\documentclass[12pt]{article}

\usepackage{amsmath,amsthm,amsfonts,amssymb,color}

\newtheorem{theorem}{Theorem}[section]

\newtheorem{proposition}[theorem]{Proposition}
\newtheorem{lemma}[theorem]{Lemma}
\newtheorem{definition}[theorem]{Definition}
\newtheorem{remark}[theorem]{Remark}

\def\cD{\mathcal{D}}
\def\cE{\mathcal{E}}
\def\cF{\mathcal{F}}
\def\cH{\mathcal{H}}

\def\cP{\mathcal{P}}
\def\cS{\mathcal{S}}

\def\bC{\mathbb{C}}
\def\bD{\mathbb{D}}
\def\bR{\mathbb{R}}

\begin{document}

\title{Hyperbolic Anderson Model with \\
space-time homogeneous Gaussian noise}

\author{Raluca M. Balan\footnote{Corresponding author. Department of Mathematics and Statistics, University of Ottawa,
585 King Edward Avenue, Ottawa, ON, K1N 6N5, Canada. E-mail
address: rbalan@uottawa.ca} \footnote{Research supported by a
grant from the Natural Sciences and Engineering Research Council
of Canada.}\and
Jian Song\footnote{Department of Mathematics. University of Hong Kong, Hong Kong. E-mail address: txjsong@hku.hk}}

\date{June 20, 2017}
\maketitle

\begin{abstract}
\noindent In this article, we study the stochastic wave equation in arbitrary spatial dimension $d$, with a multiplicative term of the form $\sigma(u)=u$, also known in the literature as the Hyperbolic Anderson Model. This equation is perturbed by a general Gaussian noise, which is homogeneous in both space and time. We prove the existence and uniqueness of a solution of this equation (in the Skorohod sense) and the H\"older continuity of its sample paths, under the same respective conditions on the spatial spectral measure of the noise as in the case of the white noise in time, regardless of the temporal covariance function of the noise.
\end{abstract}


\vspace{1mm}


\section{Introduction}

In this article, we study the stochastic wave equation with multiplicative noise:
\begin{equation}
\left\{\begin{array}{rcl}
\displaystyle \frac{\partial^2 u}{\partial t^2}(t,x) & = & \displaystyle \Delta u(t,x)+u(t,x)\dot{W}(t,x), \quad t>0, x \in \bR^d \\[2ex]
\displaystyle u(0,x) & = & 1, \quad x \in \bR^d\\[1ex]
\displaystyle \frac{\partial u}{\partial t}(0,x) & = & 0, \quad x \in \bR^d
\end{array}\right. \label{wave} 
\end{equation}

This problem is also known in the literature as the {\em Hyperbolic Anderson Model}, by analogy with the Parabolic Anderson Model in which the wave operator is replaced by the heat operator. We assume that the noise $W$ is Gaussian with covariance structure  specified by two locally integrable non-negative definitive functions $\gamma:{\bR} \to [0,\infty]$ in time and $f:\bR^d \to [0,\infty]$ in space.
Since the noise is not a martingale in time, the stochastic integral with respect to $W$ cannot be defined in the It\^o sense. To define the concept of solution we use the divergence operator from Malliavin calculus.
We refer the reader to Section \ref{section-noise} below for the precise definitions of the noise and the solution.

The Parabolic Anderson Model with the same noise $W$ as in the present article has been studied extensively in the recent years. These investigations culminated with the recent impressive article \cite{HHNT}, in which the authors have obtained a Feynman-Kac formula for the moments of the solution (for general covariance kernels $\gamma$ and $f$), as well as exponential bounds for these moments (under some quantitative conditions on $\gamma$ and $f$). The exact asymptotics for these moments were obtained in \cite{CHSX}. These extend some earlier results of \cite{hu-nualart09} and \cite{hu-nualart-song09}, in the case when the noise $W$ was fractional in space and time with index $H>1/2$ in time, and indices $H_1, \ldots,H_d>1/2$ in space.

These investigations originate in the seminal article \cite{dawson-salehi80} which studied the Parabolic Anderson Model with spatially homogenous Gaussian noise which was white in time. These authors were among the first who showed that for an equation with multiplicative noise, the (mild) solution has an explicit chaos expansion (i.e. it can be written as a series of a multiple integrals with respect to the noise), and the solution exists if and only if this series converges in $L^2(\Omega)$. Note that in \cite{dawson-salehi80} it is assumed that $f(x)=\int_{\bR^d}e^{-i \xi \cdot x}\mu(d\xi)$ for a {\em finite} measure $\mu$. In this case, the noise is a bona-fide function in the space variable $x$, whereas in the present article, the noise is a generalized function in $x$. (More precisely, if the spectral measure $\mu$ is finite, the noise $W_t(\varphi):=W(1_{[0,t]}\varphi)$ can be identified with a stationary random field $\{V_t(x)\}_{x \in \bR^d}$ via: $W_t(\varphi)=\int_{\bR^d}\varphi(x)V_t(x)dx$ for all $\varphi \in \cS(\bR^d)$.)

In contrast with its parabolic counterpart, the Hyperbolic Anderson Model with noise $W$ as above has received less attention in the literature. But there is a large amount of literature dedicated to the stochastic wave equation with spatially-homogeneous Gaussian noise which is white in time and has spectral covariance measure $\mu$ in space. (The covariance kernel $f$ is the Fourier transform of $\mu$.) We describe briefly the most important contributions in this area. In the landmark article \cite{dalang99}, Robert Dalang introduced an It\^o-type stochastic integral with respect to this noise (building upon the theory of martingale measures developed in \cite{walsh86}), and proved that the solution of the stochastic wave equation with this type of noise (and possibly a Lipschitz non-linear term $\sigma(u)$ multiplying the noise)
exists in any dimension $d=1,2,3$, provided that the measure $\mu$ satisfies what is now called {\em Dalang's condition}:
\begin{equation}
\label{Dalang-cond}
\int_{\bR^d}\frac{1}{1+|\xi|^2}\mu(d\xi)<\infty.
\end{equation}
This result was extended to any dimension $d$ in \cite{conus-dalang09}. In \cite{conus-dalang09}, it was also proved that the solution of the wave equation with
affine term $\sigma(u)=u+b$ is H\"older continuous, provided that $\mu$ satisfies:
\begin{equation}
\label{Holder-cond}
\int_{\bR^d}\left(\frac{1}{1+|\xi|^2}\right)^\beta\mu(d\xi)<\infty, \quad \mbox{for some} \quad \beta \in (0,1).
\end{equation}
A deeper study of the H\"older continuity of the solution of the wave equation in dimension $d=3$ (with general Lipschitz function $\sigma$) was carried out in \cite{dalang-sanz09} and \cite{HHN}. Exponential bounds for the moments of the solution of the Hyperbolic Anderson Model in dimension $d=3$ were obtained in \cite{dalang-mueller09}. The fact that the solution of the wave equation (with general Lipschitz function $\sigma$) has a density was proved in \cite{sanz-sus15} for any dimension $d$. In \cite{HHNS}, it was shown that this density is smooth for dimensions $d=1,2,3$.

The existence and H\"older continuity of the solution of equation \eqref{wave} with noise $W$ which is fractional in time with index $H>1/2$ and has a spatial covariance function given by the Riesz kernel $f(x)=|x|^{-\alpha},0<\alpha<d$ was proved in \cite{balan12} under the conditions $\alpha<2$, respectively $\alpha/2<\beta<1$ (which are restatements of conditions \eqref{Dalang-cond} and \eqref{Holder-cond} for the Riesz kernel). The goal of this article is to extend these results to the case of a Gaussian noise with general temporal covariance kernel $\gamma$.
In the case $d\geq 3$, the definition of solution given in \cite{balan12} is incorrect since the product between the distribution $G(t-s,x-\cdot)$ and the function $u(s,\cdot)$ is not well-defined. For this reason, we propose a new definition of the solution, and we prove its existence and uniqueness.


This article is organized as follows. In Section \ref{section-noise}, we gather some preliminary results about the space of integrands with respect to the noise $W$, and conclude with some elements of Malliavin calculus. In Section \ref{section-kernel-fn}, we study the kernels $f_n(\cdot,t,x)$ which appear in the Wiener chaos representation of the solution, and we show that the multiple Wiener integral $I_n(f_n(\cdot,t,x))$ is well-defined. In Section \ref{section-summability}, we show that the series  $\sum_{n\geq 1}I_n(f_n(\cdot,t,x))$ converges in $L^2(\Omega)$. The existence of the solution is proved in Section \ref{section-existence},  for any temporal covariance function $\gamma$ and for any spectral measure $\mu$ which satisfies \eqref{Dalang-cond}. In Section \ref{section-uniqueness}, we show that this solution is unique. In Section \ref{section-moments}, we show that the solution has uniformly bounded moments of order $p\geq 2$ and is continuous in $L^p(\Omega)$.
Section \ref{section-Holder}, we prove that this solution is H\"older continuous in time and space, provided that $\mu$ satisfies \eqref{Holder-cond}.

We specify the notation used in this article. We let $\cD_{\bC}(\bR^{d})$
be the set of complex-valued $C^{\infty}$ (i.e. infinitely differentiable) functions on $\bR^{d}$ with compact support.
We let $L_{\bC}^p(\bR^d)$ be the space of complex-valued
functions $\varphi$ on $\bR^d$ such that $|\varphi|^p$ is integrable with respect to the Lebesgue measure.
We let ${\cS}_{\bC}(\bR^d)$ be the set of complex-valued rapidly decreasing
$C^{\infty}$ functions on $\bR^d$. We denote by $\cS_{\bC}'(\bR^d)$ the space of all complex-valued tempered distributions on $\bR^d$.
Similar notations are used for spaces of real-valued elements, with the subscript $\bC$ omitted. We denote by $x \cdot y=\sum_{i=1}^{d}x_iy_i$ the inner
product in $\bR^d$ and by $|x|=(x \cdot x)^{1/2}$ the Euclidean norm
in $\bR^d$. We let $\cF \varphi(\xi)=\int_{\bR^d}e^{-i \xi \cdot x} \varphi(x)dx$ be the Fourier transform of $\varphi \in L^1(\bR^d)$. The inverse Fourier transform of $\varphi \in \cS(\bR^d)$ is
$\cF^{-1} \varphi=(2\pi)^{-d}\overline{\cF \varphi}$. We
use $\cF$ to denote the Fourier transform of functions on $\bR,\bR^d$ or $\bR^{d+1}$, but whenever there is a risk of confusion, the notation will be clearly
specified.

We end the introduction with some basic facts about distributions (see e.g. \cite{rudin91}).
The Fourier transform of $F \in \cS_{\bC}'(\bR^d)$ is a distribution $\cF F$ defined by $(\cF F,\phi)=(F,\cF \phi)$ for all $\phi \in \cS(\bR^d)$. For any $F \in \cS'(\bR^d)$ and $x \in \bR^d$, $F(x-\cdot)$ is the distribution in $\cS'(\bR^d)$ defined by
$\big(F(x-\cdot),\phi\big)=\big(F,\phi(x-\cdot)\big)$ for all $\phi \in \cS(\bR^d)$.
The product between a distribution $F \in \cS'(\bR^d)$ and a function $k \in \cS(\bR^d)$ is a distribution $Fk \in \cS'(\bR^d)$ defined by
$(Fk,\phi)=(F,k\phi)$ for all $\phi \in \cS(\bR^d)$.
The convolution between a distribution $F \in \cS'(\bR^d)$ and a function $\phi \in \cS(\bR^d)$ is a $C^{\infty}$ function $F*\phi \in \cS(\bR^d)$ with polynomial growth, defined by
$(F*\phi)(x)=(F(x-\cdot),\phi)$; its Fourier transform in $\cS'(\bR^d)$ is $\cF(F*\phi)=(\cF F) (\cF \phi)$.

\section{Preliminaries}
\label{section-noise}

In this section, we give a rigorous definition of the noise $W$, we establish a criterion for integrability with respect to $W$, and we introduce the basic elements of Malliavin calculus.

We assume that $W=\{W(\varphi);\varphi \in \cD(\bR \times \bR^d)\}$ is a zero-mean Gaussian process, defined on a complete probability space $(\Omega,\cF,P)$, with covariance
$$E[W(\varphi_1)W(\varphi_2)]=\int_{\bR^2 \times \bR^{2d}}\gamma(t-s)f(x-y)\varphi_1(t,x)\varphi_2(s,y)dxdydtds=:J(\varphi_1,\varphi_2),$$
where $\gamma:\bR \to [0,\infty]$ and $f:\bR^d \to [0,\infty]$ are continuous, symmetric, locally integrable functions, such that
$\gamma(t) <\infty$ a.e and $f(x) <\infty$ a.e.

We denote by $\cH$  the completion of $\cD(\bR \times \bR^d)$ with respect to $\langle \cdot,\cdot\rangle_{\cH}$ defined by
$$\langle \varphi_1,\varphi_2\rangle_{\cH}=J(\varphi_1,\varphi_2).$$
We are mostly interested in variables $W(\varphi)$ with $\varphi \in \cD(\bR_{+}\times \bR^d)$.

We assume that the functions $\gamma$ and $f$ are \emph{non-negative definite} (in the sense of distributions), i.e. for any $\phi \in \cS(\bR)$ and
$\varphi \in \cS(\bR^d)$
$$\int_{\bR}(\phi*\widetilde{\phi})(t)\gamma(t)dt \geq 0 \quad \mbox{and} \quad \int_{\bR^d}(\varphi*\widetilde{\varphi})(x)f(x)dx \geq 0,$$
where $\widetilde \phi(t)=\phi(-t)$ and $\widetilde \varphi(x)=\varphi(-x).$

By the Bochner-Schwartz Theorem, there exists a tempered measure $\nu$ on $\bR$ such that $\gamma$ is the Fourier transform of $\nu$ in $\cS_{\bC}'(\bR)$, i.e.
$$\int_{\bR^d}\phi(t)\gamma(t)dt=\int_{\bR}\cF \phi(\tau)\nu(d\tau) \quad \mbox{for all} \quad \phi \in \cS_{\bC}(\bR).$$
We identify two functions $\gamma_1$ and $\gamma_2$ such that $\gamma_1=\gamma_2$ a.e.
Similarly, there exists a tempered measure $\mu$ on $\bR^d$ such that $f$ is the Fourier transform of $\mu$ in $\cS_{\bC}'(\bR^d)$, i.e.
\begin{equation}
\label{Parseval-formula}
\int_{\bR^d}\varphi(x)f(x)dx=\int_{\bR^d}\cF \varphi(\xi)\mu(d\xi) \quad \mbox{for all} \quad \varphi \in \cS_{\bC}(\bR^d).
\end{equation}
We identify two functions $f_1$ and $f_2$ such that $f_1=f_2$ a.e.

It follows that for any functions $\phi_1,\phi_2 \in \cS_{\bC}(\bR)$ and  $\varphi_1,\varphi_2 \in \cS_{\bC}(\bR^d)$
\begin{equation}
\label{Parseval-t}\int_{\bR} \int_{\bR}\gamma(t-s)\phi_1(t)\overline{\phi_2(s)}dtds=\int_{\bR}\cF \phi_1(\tau)\overline{\cF \phi_2(\tau)}\nu(d\tau)
\end{equation}
and
\begin{equation}
\label{Parseval-x}
\int_{\bR^d} \int_{\bR^d}f(x-y)\varphi_1(x)\overline{\varphi_2(y)}dxdy=\int_{\bR^d}\cF \varphi_1(\xi)\overline{\cF \varphi_2(\xi)}\mu(d\xi).
\end{equation}

The next result shows that the functional $J$ has an alternative expression, in terms of Fourier transforms. In particular, this shows that $J$ is non-negative definite.

\begin{lemma}
For any $\varphi_1, \varphi_2 \in \cD(\bR \times \bR^d)$, we have:
\begin{equation}
\label{alt-expres-J}
J(\varphi_1,\varphi_2)=\int_{\bR^{d+1}}
\cF \varphi_1(\tau,\xi)\overline{\cF \varphi_2(\tau,\xi)}\nu(d\tau)\mu(d\xi),
\end{equation}
where $\cF$ denotes the Fourier transform in both variables $t$ and $x$. Moreover, $J$ is non-negative definite.
\end{lemma}

\noindent {\bf Proof:} Since $\varphi_k(t,\cdot) \in \cD(\bR^d)$ for any $t \in \bR$ and $k=1,2$, by \eqref{Parseval-x} we have:
$$J(\varphi_1,\varphi_2)=\int_{\bR}\int_{\bR}\gamma(t-s)
\left(\int_{\bR^d}\cF \varphi_1(t,\cdot)(\xi)\overline{\cF \varphi_2(s,\cdot)(\xi)}\mu(d\xi)\right)dtds.$$
For any $\xi \in \bR^d$ fixed, we denote
$\phi_{\xi}^{(k)}(t)=\cF \varphi_k (t,\cdot)(\xi)=\int_{\bR^d}e^{-i \xi \cdot x}\varphi_k(t,x)dx$. Note that $\phi_{\xi}^{(k)} \in \cD_{\bC}(\bR)$ for any $\xi \in \bR^d$ and $k=1,2$. Hence, by Fubini's theorem and \eqref{Parseval-t}, we have
\begin{eqnarray}
\nonumber
J(\varphi_1,\varphi_2)&=&
\int_{\bR^d} \left(\int_{\bR}\int_{\bR}\gamma(t-s) \phi_{\xi}^{(1)}(t)\overline{\phi_{\xi}^{(2)}(\xi)}dtds \right) \mu(d\xi)\\
\label{alt-def-J}
&=& \int_{\bR^d}\int_{\bR}\cF \phi_{\xi}^{(1)}(\tau)\overline{\phi_{\xi}^{(2)}(\tau)} \nu(d\tau) \mu(d\xi),
\end{eqnarray}
where for any $\tau \in \bR$ and $k=1,2$, we denote $$\cF \phi_{\xi}^{(k)}(\tau)=\int_{\bR}e^{-i \tau \cdot t}\phi_{\xi}^{(k)}(t)dt=\int_{\bR}e^{-i \tau \cdot t} \left(\int_{\bR^d}e^{-i \xi \cdot x}\varphi_k(t,x)dx\right)dt=\cF \varphi_k(\tau,\xi).$$

\noindent This proves \eqref{alt-expres-J}. Consequently, for any $a_1, \ldots,a_n \in \bC$ and $\varphi_1, \ldots,\varphi_n \in \cD_{\bC}(\bR_+ \times \bR^d)$,
$$\sum_{j,k=1}^{n}a_j\overline{a}_k J(\varphi_j,\varphi_k) =\int_{\bR^{d}}\int_{\bR}\left|\sum_{j=1}^{n}a_j\cF \varphi_j(\tau,\xi)\right|^2\nu(d\tau)\mu(d\xi) \geq 0.$$
This proves that $J$ is non-negative definite. $\Box$

\vspace{3mm}

The map $\varphi \mapsto W(\varphi)$ is an isometry which can be extended to $\cH$. For any $\varphi \in \cH$, we say that $W(\varphi)$ is the Wiener integral of $\varphi$ with respect to $W$ and we denote
$$W(\varphi)=\int_{\bR}\int_{\bR^d}\varphi(t,x)W(dt,dx).$$
We note that the space $\cH$ may contain distributions in $\cS'(\bR^{d+1})$ (see Theorem 3.5. of \cite{BGP12} with $F=\nu \times \mu$).

Let $|\cH|$ be the set of measurable functions $\varphi:\bR_{+} \times \bR^d \to \bR$ such that $\|\varphi\|_{|\cH|}^2:=J(|\varphi|,|\varphi|)<\infty$. Since $|\cH|$ is a complete with respect to $\|\cdot\|_{|\cH|}$ and $\|\cdot\|_{\cH} \leq \|\cdot\|_{|\cH|}$,
\begin{equation}
\label{inclusion-H}
|\cH| \subset \cH.
\end{equation}

To obtain a criterion for integrability, we need the following approximation result.

\begin{lemma}
\label{jolis-lemma}
If $\mu$ is a tempered measure on $\bR^d$, then $\cF(\cD(\bR^d))$ is dense in $\widetilde{L}_{\bC}^2(\bR^d,\mu)$, where
$$\widetilde{L}_{\bC}^2(\bR^d,\mu)=\{\varphi \in L_{\bC}^2(\bR^d,\mu); \varphi(\xi)=\overline{\varphi(-\xi)} \ \mbox{for all} \ \xi \in \bR^d\}.$$
\end{lemma}

\noindent {\bf Proof:} We refer the reader to the proof of Theorem 3.2 of \cite{jolis10} for the case $d=1$. The same argument can be used for higher dimensions $d$. $\Box$

\vspace{3mm}

We also need the following result on the ``energy'' of a complex-valued function $\varphi$ with respect to a kernel $\kappa$.

\begin{lemma}
\label{energy-lemma}
Let $m$ be a tempered measure on $\bR^d$ whose Fourier transform in $\cS_{\bC}'(\bR^d)$ is a locally integrable function $\kappa:\bR^d \to [0,\infty]$ such that $\kappa(x)<\infty$ a.e. Then for any bounded function $\varphi :\bR^d \to \bC$ with bounded support, which is continuous a.e., we have:
\begin{equation}
\label{Parseval-x-0}
\cE_{\kappa}(\varphi):=\int_{\bR^d} \kappa(x-y)\varphi(x)\overline{\varphi(y)}dxdy=
\int_{\bR^d}|\cF \varphi(\xi)|^2 m(d\xi).
\end{equation}
\end{lemma}

\begin{remark}
{\rm If we assume that $\kappa$ is a kernel of positive type (i.e. the measure $m$ is absolutely continuous with respect to the Lebesgue measure), relation \eqref{Parseval-x-0} can be deduced from Lemma 5.6 of \cite{KX09} for any function $\varphi \in L_{\bC}^1(\bR^d)$ with $\cE_{\kappa}(|\varphi|)<\infty$. In the proof of Theorem \ref{phi-in-H} below, we will use relation \eqref{Parseval-x-0} for the kernel $\kappa=\gamma$ on $\bR$ and the measure $m=\nu$. We do not use the result of \cite{KX09} since we do not assume that $\nu$ is absolutely continuous with respect to the Lebesque measure. (Relation \eqref{Parseval-x-0} will also be used in the proof of Theorem \ref{phi-in-Hn} below for the kernel $\kappa=\gamma_n$ on $\bR^n$, with
$\gamma_n(t_1,\ldots,t_n)=\prod_{i=1}^{n}\gamma(t_i)$.)
}
\end{remark}

\noindent {\bf Proof of Lemma \ref{energy-lemma}:} Suppose that $\varphi=\varphi_1+i\varphi_2$, $|\varphi(x)|\leq K$ for all $x \in \bR^n$ and the support of $\varphi$ is contained in the set $\{x \in \bR^n; |x|\leq M\}$.
We proceed by approximation. Let $p\in \cD(\bR^n)$ be such that $p \geq 0$, $\int_{\bR^n}p(x)dx=1$ and the support of $p$ is contained in $\{x \in \bR^n; |x|\leq 1\}$. For any $\varepsilon>0$, we define $p_{\varepsilon}(x)=\varepsilon^{-d}p(x/\varepsilon)$ for all $x \in \bR^d$.
Let $$\varphi_{\varepsilon}=\varphi*p_{\varepsilon}=
\varphi_{\varepsilon,1}+i\varphi_{\varepsilon,2},$$ where $\varphi_{\varepsilon,1}=\varphi_1*p_{\varepsilon}$ and $\varphi_{\varepsilon,2}=\varphi_2*p_{\varepsilon}$. Then $\varphi_{\varepsilon} \in \cD_{\bC}(\bR^d)$, $|\varphi_{\varepsilon}(x)|\leq K$ for all $x \in \bR^d$, $\varphi_{\varepsilon}(x) \to \varphi(x)$ for any continuity point $x$ of $\varphi$, and the support of $\varphi_{\varepsilon}$ is contained in the set $\{x \in \bR^d; |x| \leq M+1\}$, for any $\varepsilon \in (0,1)$. Moreover, $\cF \varphi_{\varepsilon}=\cF \varphi \cF p_{\varepsilon} \to \cF \varphi$ as $\varepsilon \downarrow 0$ and $|\cF \varphi_{\varepsilon}| \leq |\cF \varphi|$.
By the definition of the Fourier transform in $\cS_{\bC}'(\bR^d)$, for any $\varepsilon>0$,
\begin{equation}
\label{Parseval-x-1}
\int_{\bR^d} \int_{\bR^d}\kappa(x-y)\varphi_{\varepsilon}(x)\overline{\varphi_{\varepsilon}(y)}dxdy=
\int_{\bR^d}|\cF \varphi_{\varepsilon}(\xi)|^2 m(d\xi).
\end{equation}

Note that
\begin{equation}
\label{Parseval-x-2}
\lim_{\varepsilon \downarrow 0}\int_{\bR^d} \int_{\bR^d}\kappa(x-y)\varphi_{\varepsilon}(x)\overline{\varphi_{\varepsilon}(y)}
dxdy=\int_{\bR^d} \int_{\bR^d}\kappa(x-y)\varphi(x)\overline{\varphi(y)}dxdy.
\end{equation}
(This follows by applying the dominated convergence theorem to the real and imaginary part of the integrals above. In fact, since the integral on the right-hand side of \eqref{Parseval-x-1} is real-valued, the term on the left-hand side has to be real-valued for any $\varepsilon>0$, and hence its limit as $\varepsilon \downarrow 0$ is real-valued.) On the other hand, by Fatou's lemma,
\begin{equation}
\label{Parseval-x-3}
\int_{\bR^d} |\cF \varphi(\xi)|^2 m(d\xi) \leq \liminf_{\varepsilon \downarrow 0} \int_{\bR^d}|\cF \varphi_{\varepsilon}(\xi)|^2 m(d\xi).
\end{equation}
From \eqref{Parseval-x-1}, \eqref{Parseval-x-2} and \eqref{Parseval-x-3}, we obtain that
$$\int_{\bR^d} |\cF \varphi(\xi)|^2 m(d\xi) \leq \int_{\bR^d} \int_{\bR^d}\kappa(x-y)\varphi(x)\overline{\varphi(y)}dxdy.$$

Hence, if the right-hand side of \eqref{Parseval-x-0} is infinite, so must be the left-hand side. If the right-hand side of \eqref{Parseval-x-0} is finite, then by the dominated convergence theorem, we have:
\begin{equation}
\label{Parseval-x-4}
\int_{\bR^d} |\cF \varphi(\xi)|^2 m(d\xi) = \lim_{\varepsilon \downarrow 0} \int_{\bR^d}|\cF \varphi_{\varepsilon}(\xi)|^2 m(d\xi).
\end{equation}
In this case, relation \eqref{Parseval-x-0} follows by \eqref{Parseval-x-1}, \eqref{Parseval-x-2} and \eqref{Parseval-x-4}. $\Box$

\begin{remark}
\label{def-version}
{\rm
Recall that the Fourier transform $\cF S$ of a distribution $S \in \cS'(\bR^d)$ is defined by $\big(\cF S,\varphi\big)=\big(S,\cF \varphi \big)$ for all $\varphi \in \cS(\bR^d)$. When $S$ is a genuine distribution and $\cF S=g$ is a function, this means that
\begin{equation}
\label{Fourier-S}
\int_{\bR^d}g(\xi)\varphi(\xi)d\xi=\big(S,\cF \varphi\big) \quad \mbox{for all} \quad \varphi \in \cS(\bR^d).
\end{equation}
In this case, $\cF S$ is understood as the equivalence class of all functions $g$ which satisfy \eqref{Fourier-S}. If $g$ is an element of this class, we say that $g$ is a {\em version} of the Fourier transform $\cF S$. If $g_1$ and $g_2$ are versions of $\cF S$, then $g_1=g_2$ a.e.}
\end{remark}

This leads us to the following hypothesis.

\vspace{3mm}

{\bf Hypothesis A.} $\mu$ is absolutely continuous with respect to the Lebesgue measure.

\vspace{3mm}

Using the alternative expression given by \eqref{alt-def-J} for the inner product $\langle \varphi_1,\varphi_2\rangle_{\cH}$ and the previous lemmas, we obtain the following criterion for integrability with respect to $W$.

\begin{theorem}
\label{phi-in-H}
Let $\bR \ni t \mapsto S(t) \in \cS'(\bR^{d})$ be a deterministic function such that $\cF S(t,\cdot)$ is a function for all $t\in \bR$. If $\cF S(t,\cdot)$ is uniquely determined only up to a set of Lebesgue measure zero, we assume that $\mu$ satisfies Hypothesis A.
Suppose that:\\
(i) for each $t \in \bR$, there exists a version of $\cF S(t,\cdot)$ such that $(t,\xi) \mapsto \cF S(t,\cdot)(\xi)=:\phi_{\xi}(t)$ is measurable on $\bR \times \bR^{d}$;\\
(ii) for all $\xi \in \bR^d$, $\int_{\bR}|\phi_{\xi}(t)|dt<\infty$.

Then the following statements hold:

a) The function $(\tau,\xi)\mapsto \cF \phi_{\xi}(\tau)$ is measurable on $\bR \times \bR^{d}$, where $\cF \phi_{\xi}$ denotes the Fourier transform of $\phi_{\xi}$, i.e.
$\cF \phi_{\xi}(\tau)=\int_{\bR}
e^{-i \tau t}\phi_{\xi}(t)dt$, $\tau \in \bR$.

b) If
\begin{equation}
\label{0-norm-finite}
\|S\|_0^2:=\int_{\bR^{d}} \int_{\bR} |\cF \phi_{\xi}(\tau)|^2 \nu(d\tau)\mu(d\xi)<\infty
\end{equation}
then $S \in \cH$ and $\|S\|_{\cH}^2=\|S\|_0^2$.

c) Assume in addition that $S(t,\cdot)=0$ for all $t \not \in [0,T]$, for some $T>0$. If for every $\xi \in \bR^d$, the function $t \mapsto \cF S(t,\cdot)(\xi)$ is continuous a.e. and bounded on $[0,T]$, and
$$I_T:=\int_{\bR^d}\int_0^T \int_0^T \gamma(t-s)\cF S(t,\cdot)(\xi)\overline{\cF S(t,\cdot)(\xi)}dtds \mu(d\xi)<\infty,$$
then \eqref{0-norm-finite} holds, $S \in \cH$ and $\|S\|_{\cH}^2=\|S\|_{0}^2=I_T$.
\end{theorem}

\noindent {\bf Proof:} a) This follows by Fubini's theorem, using the fact that $(t,\tau,\xi) \mapsto e^{-i \tau t}\phi_{\xi}(t)$ is measurable on $\bR \times \bR \times \bR^d$, by hypothesis (i).

b) 
Let $a(\tau,\xi)=\cF \phi_{\xi}(\tau)$. By \eqref{0-norm-finite} and part a), the function $a$ lies in $L_{\bC}^2(\bR^{d+1},\Pi)$, where $$\Pi(d\tau,d\xi)=\nu(d\tau)\mu(d\xi).$$

\noindent We denote $\widetilde{L}_{\bC}^2(\bR^{d+1},\Pi)=\{\varphi \in L_{\bC}^2(\bR^{d+1},\Pi); \varphi(\tau,\xi)=\overline{\varphi(-\tau,-\xi)} \ \mbox{for all} \ \tau \in \bR, \xi \in \bR^d\}$.
We observe that
$a \in \widetilde{L}_{\bC}^2(\bR^{d+1},\Pi)$, since by Lemma 3.3 of \cite{BGP12},
$$\phi_{-\xi}(t)=\cF S(t,\cdot)(-\xi)=\overline{\cF S(t,\cdot)(\xi)}=\overline{\phi_{\xi}(t)} \quad \mbox{for all} \ \xi \in \bR^d,$$
and hence
$$a(-\tau,-\xi)=\int_{\bR}e^{i\tau t}\phi_{-\xi}(t)dt=\int_{\bR}\overline{e^{-i\tau t}\phi_{\xi}(t)}dt=\overline{a(\tau,\xi)} \quad \mbox{for all} \ \tau \in \bR,\xi \in \bR^d.$$

By Lemma \ref{jolis-lemma}, $\cF(\cD(\bR^{d+1}))$ is dense in
$\widetilde{L}_{\bC}^2(\bR^{d+1},\Pi)$. Hence, for any $\varepsilon>0$, there exists a function $l=l(\varepsilon) \in \cD(\bR^{d+1})$ such that
$$\int_{\bR^{d+1}}|a(\tau,\xi)-\cF l(\tau,\xi)|^2 \Pi(d\tau,d\xi)<\varepsilon^2.$$
Note that the previous integral is
$\int_{\bR^{d+1}}|\cF \phi_{\xi}(\tau)-\cF \psi_{\xi}(\tau)|^2 \Pi(d\tau,d\xi)=:\|S-l\|_{0}^{2}$, where $\cF \psi_{\xi}$ is the Fourier transform of the function $t \mapsto \psi_{\xi}(t)=\cF l(t,\cdot)(\xi)$.
The conclusion follows using expression \eqref{alt-def-J} for the inner product in $\cH$.

c) For every $\xi \in \bR^d$ fixed, we apply Lemma \ref{energy-lemma} to the bounded function $\phi_{\xi}:\bR \to \bC$ which is continuous a.e and has support contained in $[0,T]$. We apply this lemma for the measure $m=\nu$ and the kernel $\kappa=\gamma$ on $\bR$. We obtain that, for any $\xi \in \bR^d$,
$$\int_{0}^T \int_{0}^T \gamma(t-s) \phi_{\xi}(t)\overline{\phi_{\xi}(s)}
dtds=\int_{\bR}|\cF \phi_{\xi}(\tau)|^2 \nu(d\tau).$$
We integrate with respect to $\mu(d\xi)$ and we multiply by $(2\pi)^{-d}$.
We obtain that
$$I_T=\int_{\bR^d}\int_{\bR}|\cF \phi_{\xi}(\tau)|^2 \nu(d\tau)\mu(d\xi)=:\|S\|_0^2.$$
Since $I_T<\infty$, it follows that $\|S\|_0^2<\infty$. The conclusion follows by part b).
$\Box$

\vspace{3mm}

We are interested in applying Theorem \ref{phi-in-H} to the case when $\varphi$ is related to the fundamental solution $G$ of the wave equation on $\bR_{+}\times \bR^d$. We recall that:
\begin{eqnarray*}
G(t,x)&=& \frac{1}{2}1_{\{|x| <t\}} \quad \mbox{if} \ d=1\\
G(t,x)&=& \frac{1}{2\pi}\frac{1}{\sqrt{t^2-|x|^2}}1_{\{|x|<t\}} \quad \mbox{if} \ d=2\\
G(t,\cdot)&=&\frac{1}{4\pi t}\sigma_t, \quad \mbox{if} \ d=3,\\\
\end{eqnarray*}
where $\sigma_t$ is the surface measure on the sphere $\{x \in \bR^3; |x|=t\}$.
If $d=1$ or $d=2$, $G(t,\cdot)$ is a non-negative function in $L^1(\bR^d)$, and if $d=3$, $G(t,\cdot)$ is a finite measure in $\bR^3$.

If $d \geq 4$ is even, $G(t,\cdot)$ is a distribution with compact support in $\bR^d$ given by:
$$G(t,\cdot)= \frac{1}{1 \cdot 3 \cdot \ldots \cdot (d-1)}\left(\frac{1}{t}\frac{\partial}{\partial t}
\right)^{(d-2)/2}(t^{d-1}\Upsilon_t), \quad \Upsilon_t(\varphi)=\frac{1}{\omega_{d+1}}\int_{B(0,1)}
\frac{\varphi(ty)}{\sqrt{1-|x|^2}} dx,$$
and if
$d \geq 5$ is odd, $G(t,\cdot)$ is a distribution with compact support in $\bR^d$ given by:
$$G(t,\cdot)= \frac{1}{1 \cdot 3 \cdot \ldots \cdot (d-2)}\left(\frac{1}{t}\frac{\partial}{\partial t}
 \right)^{(d-3)/2}(t^{d-2}\Sigma_t), \quad \Sigma_t(\varphi)=\frac{1}{\omega_d}\int_{\partial B(0,1)}\varphi(tz)d\sigma(z),$$
where $\omega_d$ is the surface area of the unit sphere $\partial B(0,1)$ in $\bR^d$, and $\sigma$ is the surface measure on $\partial B(0,1)$
(see e.g. Theorem (5.28), page 176 of \cite{folland95}).

It is known that for any $d \geq $1, the Fourier transform of $G(t,\cdot)$ is given by:
\begin{equation}
\label{Fourier-G}
\cF G(t,\cdot)(\xi)=\frac{\sin(t|\xi|)}{|\xi|}, \quad \xi \in \bR^d.
\end{equation}

Note that when $d=1,2,3$, the previous formula uniquely determines $\cF G(t,\cdot)$ as the Fourier transform of a function in $L^1(\bR^d)$ for $d=1,2$, or the Fourier transform of a finite measure for $d=3$. But when $d \geq 4$, \eqref{Fourier-G} is interpreted in the sense of distributions, and the definition of $\cF G(t,\cdot)$ is unique only up to a set of Lebesgue measure zero.

We have the following result about the integrability of $G$.

\begin{theorem}
\label{gtx-in-H}
For any $t>0$ and $x \in \bR^d$, we define $g_{t,x}(s,\cdot)=G(t-s,x-\cdot)1_{[0,t]}(s)$ for any $s \in \bR$. If $d \geq 4$, we assume that $\mu$ satisfies Hypothesis A. Suppose that
\begin{equation}
\label{additive-noise-cond}
I_t:=\int_{\bR^d}\int_{0}^t\int_0^t \gamma(r-s) \frac{\sin ((t-r)|\xi|)\sin((t-s)|\xi|)}{|\xi|^2}drds \mu(d\xi)<\infty
\end{equation}
for any $t>0$.
Then, for any $t>0$ and $x \in \bR^d$, $g_{t,x} \in \cH$, the stochastic integral
$$v(t,x):=\int_0^t \int_{\bR}G(t-s,x-y)W(ds,dy)$$ is well-defined and $E|v(t,x)|^2=I_t$. In particular, \eqref{additive-noise-cond} holds for any $t>0$ if the measure $\mu$ satisfies \eqref{Dalang-cond}.
(Note that $v$ is the solution of the linear wave equation
$\frac{\partial^2 v}{\partial x^2}(t,x)=\Delta v(t,x)+\dot{W}(t,x), t>0, x\in \bR^d$
with zero initial conditions.)
\end{theorem}

\noindent {\bf Proof:} By applying Theorem \ref{phi-in-H}.c) to the function $S=g_{t,x}$ we infer that $g_{t,x} \in \cH$ and $\|g_{t,x}\|_{\cH}^2=\|g_{t,x}\|_0^2=I_t$.
To see that $g_{t,x}$ satisfies the conditions of this theorem, we note that, due to \eqref{Fourier-G}, for all $s \in \bR$ and $\xi \in \bR^d$,
$$\phi_{\xi}(s):=\cF g_{t,x}(s,\cdot)(\xi)=e^{-i \xi \cdot x}\frac{\sin((t-s)|\xi|)}{|\xi|}1_{[0,t]}(s).$$

\noindent Then $|\phi_{\xi}(s)|\leq (t-s)1_{[0,t]}(s) \leq t 1_{[0,t]}(s)$ for all $s\in \bR$ and $\xi \in \bR$. It follows that
$g_{t,x}$ satisfies conditions (i) and (ii) of Theorem \ref{phi-in-H}.
By the construction of the stochastic integral, $E|v(t,x)|^2=\|g_{t,x}\|_{0}^2=I_t$.

Note that $I_t$ coincides with the term $\alpha_1(t)$ which appears in the series representation \eqref{second-moment-u} of the second moment of the solution $u(t,x)$ to equation \eqref{wave}. (See definition \eqref{def-alpha-n} of $\alpha_n(t)$ below.) By relations \eqref{bound1-alpha} and \eqref{bound-psi} below,
we see that if \eqref{Dalang-cond} holds, $\alpha_1(t)<\infty$. $\Box$

\begin{remark}

{\rm Theorem \ref{phi-in-H}.c) can also be applied to the function $S=g_{t,x}$ where
$g_{t,x}(s,\cdot)=G(t-s,x-\cdot)1_{[0,t]}(s)$ and
\begin{equation}
\label{def-G-heat}
G(t,x)=\frac{1}{(2\pi t)^{d/2}} \exp\left(-\frac{|x|^2}{2t}\right)
\end{equation}
is the fundamental solution of the heat equation $\frac{\partial u}{\partial t}=\frac{1}{2}\Delta u$ on $\bR_{+} \times \bR^d$. Since $g_{t,x}(s,\cdot) \in L^1(\bR^d)$, its Fourier transform is uniquely determined and we do not need to assume that $\mu$ is absolutely continuous with respect to the Lebesgue measure. Note that
$g_{t,x}\in \cH$ provided that, for any $t>0$,
$$I_t:=\int_{\bR^d}\int_0^t \int_0^t \gamma(r-s) \exp\left(-\frac{(t-r)|\xi|^2}{2}\right)\exp\left(-\frac{(t-s)|\xi|^2}{2}\right)drds \mu(d\xi)<\infty.$$
In this case, $v(t,x)=W(g_{t,x})$ is the solution of $\frac{\partial v}{\partial t}=\frac{1}{2}\Delta v+\dot{W}$ and $E|v(t,x)|^2=I_t$.}

\end{remark}

We conclude this section by recalling briefly some basic elements of Malliavin calculus (see \cite{nualart06} for more details).

It is known that every square-integrable random variable $F$ which is measurable with respect to $W$, has the Wiener chaos expansion: $$F=E(F)+\sum_{n \geq 1}F_n \quad \mbox{with} \quad F_n \in \cH_n,$$
where $\cH_n$ is the $n$-th Wiener chaos space associated to $W$. Moreover, each $F_n$ can be represented as $F_n=I_n(f_n)$ for some $f_n \in \cH^{\otimes n}$, where $\cH^{\otimes n}$ is the $n$-th tensor product of $\cH$ and
$I_n:\cH^{\otimes n} \to \cH_n$ is the multiple Wiener integral with respect to $W$.
By the orthogonality of the Wiener chaos spaces and an isometry-type property of $I_n$, we obtain that
$$E|F|^2=(EF)^2+\sum_{n \geq 1}E|I_n(f_n)|^2=(EF)^2+\sum_{n \geq 1}n! \|\widetilde{f}_n\|_{\cH^{\otimes n}}^{2},$$
where $\widetilde{f}_{n}$ is the symmetrization of $f_n$ in all $n$ variables.
We note that
the space $\cH^{\otimes n}$ may contain distributions in $\cS'(\bR^{n(d+1)})$.

We denote by $\delta: {\rm Dom}(\delta) \subset L^2(\Omega;\cH) \to L^2(\Omega)$ the divergence operator with respect to $W$, defined as the adjoint of the Malliavin derivative  $D$ with respect to $W$.
If $u \in \mbox{Dom} \ \delta$, we use the notation
$$\delta(u)=\int_0^{\infty} \int_{\bR^d}u(t,x) W(\delta t, \delta x),$$
and we say that
$\delta(u)$ is the {\em Skorohod integral} of $u$ with respect to $W$. In particular, $E[\delta(u)]=0$.

\section{The kernels $f_n$}
\label{section-kernel-fn}

In this section, we give the definition of the kernels $f_n(\cdot,t,x)$ which appear in the Wiener chaos representation of the solution to equation \eqref{wave}, and we prove that they are integrable with respect to the noise $W$.

Let $t_1>0, \ldots,t_n>0$ be arbitrary. If $d\leq 2$, for any $x_1, \ldots,x_n \in \bR^d$, we define
\begin{equation}
\label{def-fn-1}
f_n(t_1,x_1, \ldots,t_n,x_n,t,x)=1_{\{0<t_1<\ldots<t_n<t\}}G(t-t_n,x-x_n) \ldots G(t_2-t_1,x_2-x_1).
\end{equation}
In this case, $f_n(t_1,\cdot, \ldots,t_n,\cdot,t,x)$ is a function in $L^1(\bR^{nd})$ whose Fourier transform is
\begin{align}
\label{def-Fourier-fn}
& \cF f_n(t_1,\cdot,\ldots,t_n,\cdot,t,x)(\xi_1,\ldots,\xi_n)=
1_{\{0<t_1<\ldots<t_n<t\}}e^{-i(\xi_1+\ldots+\xi_n)\cdot x}\overline{\cF G(t_2-t_1,\cdot)(\xi_1)} \\
\nonumber
& \quad \quad \quad \quad \quad \quad  \quad \quad \quad  \overline{\cF G(t_3-t_2,\cdot)(\xi_1+\xi_2)} \ldots \overline{\cF G(t-t_n,\cdot)(\xi_1+\ldots+\xi_n)}.
\end{align}

If $d\geq 3$, we let $f_n(t_1,\cdot,\ldots,t_n,\cdot,t,x)$ be the distribution in $\cS'(\bR^{nd})$ whose Fourier transform is given by \eqref{def-Fourier-fn}.
More precisely, the action of the $f_n(t_1,\cdot,\ldots,t_n,\cdot,t,x)$ on a test function $\phi \in \cS(\bR^{nd})$ is given by:
\begin{align}
\nonumber
& \Big(f_n(t_1,\cdot,\ldots,t_n,\cdot,t,x),\phi \Big)= \Big(\cF f_n(t_1,\cdot,\ldots,t_n,\cdot,t,x),\cF^{-1}\phi \Big)\\
\nonumber
&  \quad  = 1_{\{0<t_1<\ldots<t_n<t\}}\int_{\bR^{nd}} e^{-i(\xi_1+\ldots+\xi_n)\cdot x}\overline{\cF G(t_2-t_1,\cdot)(\xi_1)} \cdot \overline{\cF G(t_3-t_2,\cdot)(\xi_1+\xi_2)} \ldots \\
\label{def-kernel-fn}
& \quad \quad \quad \quad \quad \quad \quad \quad \quad \overline{\cF G(t-t_n,\cdot)(\xi_1+\ldots+\xi_n)} \, \cF^{-1} \phi(\xi_1,\ldots,\xi_n)d\xi_1 \ldots d\xi_n.
\end{align}

\begin{lemma}
\label{fn-distribution}
$f_n(t_1,\cdot,\ldots,t_n,\cdot,t,x)$ is a well-defined distribution in $\cS'(\bR^{nd})$.
\end{lemma}

\noindent {\bf Proof:} First, note that for any $t>0$ and $\xi \in \bR^d$, letting $C_t=2(t^2 \vee 1)$, we have:
 \begin{equation}
\label{UB-Fourier-G}
|\cF G(t,\cdot)(\xi)|^2 \leq C_t\frac{1}{1+|\xi|^2},
\end{equation}
since
$\frac{\sin^2(t|\xi|)}{|\xi|^2} \leq t^2 \leq t^2 \frac{2}{1+|\xi|^2}$ if $|\xi|\leq 1$
and
$\frac{\sin^2(t|\xi|)}{|\xi|^2} \leq \frac{1}{|\xi|^2} \leq \frac{2}{1+|\xi|^2}$ if $|\xi|>1$.

The integral on the right-hand side of \eqref{def-kernel-fn} is finite, since $\cF^{-1}\phi$ is in $L^1(\bR^{nd})$ and
$|\cF G(t,\cdot)(\xi)|\leq C_t^{1/2}$
 for any $t>0,\xi \in \bR^d$.  The map $\phi \mapsto
\Big(f_n(t_1,\cdot,\ldots,t_n,\cdot,t,x),\phi \Big)$ is clearly linear. To show that this map is continuous, assume that $\phi_k \to \phi$ in $\cS(\bR^{nd})$ as $k \to \infty$. Then $\cF^{-1}\phi_k \to \cF^{-1}\phi$ in $\cS_{\bC}(\bR^{nd})$ as $k\to \infty$, and hence for any integer $m>0$,
$$T_k:=\sup_{\xi_1, \ldots,\xi_n \in \bR^{d}}(1+|\xi_1|^2)^m \ldots (1+|\xi_1+\ldots+\xi_n|^2)^m \, |\cF^{-1}(\phi_k-\phi)(\xi_1,\ldots,\xi_n)| \to 0,$$
as $k \to \infty$. Using \eqref{UB-Fourier-G}, we see that
$\left|\Big(f_n(t_1,\cdot,\ldots,t_n,\cdot,t,x),\phi_k-\phi \Big)\right|$ is smaller than
\begin{align*}
& T_k C_t^{n/2} 1_{\{0<t_1<\ldots<t_n<t\}} \int_{\bR^{nd}} \left(\frac{1}{1+|\xi_1|^2}\right)^{m+1/2}\ldots \left(\frac{1}{1+|\xi_1+\ldots+\xi_n|^2}\right)^{m+1/2} d\xi_1 \ldots d\xi_n,
\end{align*}
and hence converges to $0$ as $k \to \infty$. The last integral is finite if $2m+1>d$. $\Box$

\vspace{3mm}

We will need the following result. Recall that for any $t>0$, $G(t,\cdot)$ is a distribution with compact support (hence is in $\cS'(\bR^d)$).

\begin{lemma}
For any $\phi \in \cS(\bR^d)$, $G(t,\cdot)*\phi$ is a well-defined function in $\cS(\bR^d)$.
\end{lemma}

\noindent {\bf Proof:} We first show that $G(t,\cdot)*\phi$ is well-defined. For any $x \in \bR^d$, we have
\begin{eqnarray}
\nonumber
\lefteqn{(G(t,\cdot)*\phi)(x)=(G(t,x-\cdot),\phi)=(\cF G(t,x-\cdot),\cF^{-1}\phi)} \\
\label{def-G*phi}
& & = \int_{\bR^d}\cF G(t,x-\cdot)(\xi)\cF^{-1}\phi(\xi)d\xi= \int_{\bR^d}e^{-i \xi \cdot x} \,\frac{\sin(t|\xi|)}{|\xi|}\cF^{-1}\phi (\xi)d\xi.
\end{eqnarray}
The function $h(\xi)=\frac{\sin(t|\xi|)}{|\xi|},\xi \in \bR^d$ is analytic (hence infinitely differentiable), being the Fourier transform of a distribution with compact support (see e.g. Theorem 7.23 of \cite{rudin91}). It can be shown that the partial derivatives of $h$ are bounded. Therefore, the function $g(\xi):= h(\xi)\cF^{-1}\phi (\xi), \xi \in \bR^d$ is in $\cS(\bR^d)$, and its Fourier transform lies in $\cS_{\bC}(\bR^d)$. Since $g(-\xi)=\overline{g(\xi)}$ for any $\xi \in \bR^d$, the Fourier transform of $g$ is real-valued.
$\Box$

\vspace{3mm}

The next result identifies the action of the distribution $f_n(t_1,\cdot,\ldots,t_n,\cdot,t,x)$ on a product test function
$\phi=\phi_1 \otimes \ldots \otimes \phi_n$, i.e.
$\phi(x_1,\ldots,x_n)=\phi_1(x_1)  \ldots  \phi_n(x_n)$ for any $x_1, \ldots,x_n \in \bR^d$.

\begin{lemma}
\label{def-kn-lemma}
If $\phi=\phi_1 \otimes \ldots \otimes \phi_n$ with $\phi_1, \ldots, \phi_n \in \cS(\bR^d)$, then
$$\Big(f_n(t_1,\cdot,\ldots,t_n,\cdot,t,x),\phi \Big)=\varphi_n(t_1,\ldots,t_n,t,x),$$
where the pairs of functions $(\psi_1,\varphi_1),\ldots,(\psi_n,\varphi_n)$ are defined recursively as follows:
$$\psi_1=\phi_1, \quad \varphi_1(t_1,t_2,\cdot)=1_{(0,t_2)}(t_1)G(t_2-t_1,\cdot)*\psi_1,$$
$$\psi_2(t_1,t_2,\cdot)=\phi_2 \varphi_1(t_1,t_2,\cdot), \quad \varphi_2(t_1,t_2,t_3,\cdot)=1_{(0,t_3)}(t_2)G(t_3-t_2,\cdot)*\psi_2(t_1,t_2,\cdot),\ldots,$$
$$\psi_n(t_1,\ldots,t_n,\cdot)=\phi_{n} \varphi_{n-1}(t_1,\ldots,t_n,\cdot),$$ $$\varphi_n(t_1, \ldots,t_{n},t,\cdot)=1_{(0,t)}(t_n)G(t-t_n,\cdot)*\psi_n(t_1,\ldots,t_n,\cdot).$$
\end{lemma}

\noindent {\bf Proof:} The statement is clear for $n=1$ since by \eqref{def-G*phi},
\begin{align*}
\Big(f_1(t_1,\cdot,t,x),\phi_1\Big)&=1_{(0,t)}(t_1)\int_{\bR^d}e^{-i \xi_1 \cdot x}\overline{\cF G(t-t_1,\cdot)(\xi_1)}\,\cF^{-1}\phi_1(\xi_1)d\xi_1\\
&=1_{(0,t)}(t_1)(G(t-t_1,\cdot)*\phi_1)(x).
\end{align*}

Assume that $n\geq 2$. Note that $\cF^{-1}\phi=h_1 \otimes \ldots \otimes h_n$ where $h_i=\cF^{-1}\phi_i$ for $i=1, \ldots,n$.
Using the change of variables $\eta_k=\xi_1+\ldots+\xi_k$ with $k=1, \ldots,n$, we see that relation \eqref{def-kernel-fn} can be written as
$$\Big(f_n(t_1,\cdot,\ldots,t_n,\cdot,t,x),\phi \Big)= 1_{(0,t)}(t_n) \int_{\bR^d}e^{-i\eta_n\cdot x} \overline{\cF G(t-t_n,\cdot)(\eta_n)} \,g_n(t_1,\ldots,t_n,\eta_n)\,d\eta_n,$$
where
\begin{align*}
g_k(t_1,\ldots,t_k,\eta_k)=&1_{\{0<t_1<\ldots<t_{k-1}<t_k\}}\int_{\bR^{(k-1)d}} \overline{\cF G(t_k-t_{k-1},\cdot)(\eta_{k-1})} \ldots \overline{\cF G(t_2-t_1,\cdot)(\eta_{1})}\\
& \quad \quad \quad h_1(\eta_1)h_2(\eta_2-\eta_1)\ldots h_{k}(\eta_{k}-\eta_{k-1})d\eta_1 \ldots d\eta_{k-1}
\end{align*}
for $k=2,\ldots,n$. We show that for
\begin{equation}
\label{gk-psi-k-induction}
g_k(t_1,\ldots,t_k,\cdot)=\cF^{-1}\psi_k(t_1,\ldots,t_k,\cdot) \quad \mbox{for all} \ k=2,\ldots,n.
 \end{equation}
This is proved by induction on $k$. For $k=2$, we have:
\begin{align*}
& \cF^{-1}\psi_2(t_1,t_2,\cdot)(\eta_2)=\frac{1}{(2\pi)^d}\int_{\bR^d}e^{i\eta_2 \cdot x_2} \phi_2(x_2) \varphi_1(t_1,t_2,x_2)dx_2\\
&=\frac{1}{(2\pi)^d}\int_{\bR^d}e^{i\eta_2 \cdot x_2} \phi_2(x_2) 1_{(0,t_2)}(t_1)
\left(\int_{\bR^d} e^{-i\eta_1\cdot x_2} \overline{\cF G(t_2-t_1,\cdot)(\eta_1)}\cF^{-1}\phi_1(\eta_1)d\eta_1\right)dx_2\\
&=  1_{(0,t_2)}(t_1)\int_{\bR^d} \overline{\cF G(t_2-t_1,\cdot)(\eta_1)} \, \cF^{-1}\phi_1(\eta_1) \left(\frac{1}{(2\pi)^d}\int_{\bR^d}e^{i(\eta_2-\eta_1) \cdot x_2} \phi_2(x_2) dx_2 \right)d\eta_1\\
&=1_{(0,t_2)}(t_1)\int_{\bR^d} \overline{\cF G(t_2-t_1,\cdot)(\eta_1)} \, \cF^{-1}\phi_1(\eta_1) \cF^{-1}\phi_2(\eta_2-\eta_1)d\eta_1=g_2(t_1,t_2,\eta_2),
\end{align*}
where we used the definition of $\varphi_1$ and \eqref{def-G*phi} for the second equality. Assume that the statement is true for $k-1$. Then
\begin{align*}
& \cF^{-1}\psi_k(t_1,\ldots,t_k,\cdot)(\eta_k)=\frac{1}{(2\pi)^d}\int_{\bR^d}e^{i\eta_k \cdot x_k} \phi_k(x_k) \varphi_{k-1}(t_1,\ldots,t_k,x_k)dx_k\\
&=\frac{1}{(2\pi)^d}\int_{\bR^d}e^{i\eta_k \cdot x_k} \phi_k(x_k) 1_{(0,t_k)}(t_{k-1})
\left(\int_{\bR^d} e^{-i\eta_{k-1}\cdot x_{k}} \overline{\cF G(t_k-t_{k-1},\cdot)(\eta_{k-1})} \right.\\ & \quad \quad \quad \quad \quad
 \quad \quad \quad \quad  \cF^{-1}\psi_{k-1}(t_1,\ldots,t_{k-1},\cdot)(\eta_{k-1})
d\eta_{k-1}\Big)dx_k\\
&=  1_{(0,t_k)}(t_{k-1})\int_{\bR^d} \overline{\cF G(t_k-t_{k-1},\cdot)(\eta_{k-1})} \, \ \cF^{-1}\psi_{k-1}(t_1,\ldots,t_{k},\cdot)(\eta_{k-1}) \\ & \quad \quad \quad \quad \quad
 \quad \quad \quad \quad\left(\frac{1}{(2\pi)^d}\int_{\bR^d}e^{i(\eta_k-\eta_{k-1}) \cdot x_k} \phi_k(x_k) dx_k \right)d\eta_{k-1}\\
&=1_{(0,t_k)}(t_{k-1})\int_{\bR^d} \overline{\cF G(t_k-t_{k-1},\cdot)(\eta_{k-1})}
\,g_{k-1}(t_1,\ldots,t_{k-1},\eta_{k-1})
  \cF^{-1}\phi_k(\eta_k-\eta_{k-1})d\eta_k\\
&=g_k(t_1,\ldots,t_k,\eta_k),
\end{align*}
where we used the definition of $\varphi_{k-1}$ and \eqref{def-G*phi} for the second equality, and the induction hypothesis for the fourth equality. This concludes the proof of \eqref{gk-psi-k-induction}.

Using \eqref{gk-psi-k-induction}, it follows that
\begin{align*}
\Big(f_n(t_1,\cdot,\ldots,t_n,\cdot,t,x),\phi \Big)&= 1_{(0,t)}(t_n) \int_{\bR^d}e^{-i\eta_n\cdot x} \overline{\cF G(t-t_n,\cdot)(\eta_n)}
\cF^{-1}\psi_n(t_1,\ldots,t_n,\cdot)(\eta_n)\,d\eta_n\\
& = 1_{(0,t)}(t_n) \big(G(t-t_n,\cdot)* \psi_n(t_1,\ldots,t_n,\cdot)\big)(x)\\
&=\varphi_n(t_1,\ldots,t_n,x),
\end{align*}
where we used \eqref{def-G*phi} for the second equality. $\Box$

\begin{remark}
\label{recurrence-relation-fn}
{\rm Lemma \ref{def-kn-lemma} gives the relationship between the kernels $f_{n+1}$ and $f_n$:
\begin{align*}
& \Big(f_{n+1}(t_1,\cdot,\ldots,t_{n+1},\cdot,t,x),\phi_1 \otimes \ldots \otimes \phi_n \otimes \phi_{n+1} \Big)=\\
& \quad \quad \quad 1_{(0,t)}(t_{n+1}) \Big(G(t-t_{n+1},x-\cdot),\phi_{n+1} \big(f_{n}(t_1,\cdot,\ldots,t_n,\cdot,*),\phi_1 \otimes \ldots \otimes \phi_n\big) \Big),
\end{align*}
where on the right-hand side we have the action of the distribution $G(t-t_{n+1},x-\cdot)$ on the function $y \mapsto \phi_{n+1}(y) \big(f_{n}(t_1,\cdot,\ldots,t_n,\cdot,y),\phi_1 \otimes \ldots \otimes \phi_n\big)$. }
\end{remark}

To check that the kernel $f_n(\cdot,t,x)$ lies in $\cH^{\otimes n}$, we need the following result, which is the counterpart of Theorem \ref{phi-in-H} for multiple Wiener integrals of order $n$.

\begin{theorem}
\label{phi-in-Hn}
Let $\bR^n \ni (t_1,\ldots,t_n) \mapsto S(t_1,\cdot,\ldots,t_n,\cdot) \in \cS'(\bR^{nd})$ be a deterministic function such that $\cF S(t_1,\cdot,\ldots,t_n,\cdot)$ is a function for all $(t_1, \ldots,t_n)\in \bR^n$.
If $\cF S(t_1,\cdot, \ldots,t_n,\cdot)$ is uniquely determined only up to a set of Lebesgue measure zero, we assume that $\mu$ satisfies Hypothesis A.
Suppose that:\\
(i) for each $(t_1, \ldots,t_n) \in \bR^n$, there exists a version of $\cF S(t_1,\cdot,\ldots,t_n,\cdot)$ such that the function $(t_1,\ldots,t_n, \xi_1,\ldots,\xi_n) \mapsto \cF S(t_1,\cdot,\ldots,t_n,\cdot)(\xi_1, \ldots,\xi_n)=:\phi_{\xi_1, \ldots,\xi_n}(t_1,\ldots,t_n)$ is measurable on $\bR^n \times \bR^{nd}$;\\
(ii) for all $\xi_1, \ldots,\xi_n \in \bR^d$, $\int_{\bR^n}|\phi_{\xi_1,\ldots,\xi_n}(t_1,\ldots,t_n)|dt_1 \ldots dt_n<\infty$.

Then the following statements hold:

a) The function $(\tau_1, \ldots,\tau_n,\xi_1,\ldots,\xi_n)\mapsto \cF \phi_{\xi_1,\ldots,\xi_n}(\tau_1,\ldots,\tau_n)$ is measurable on $\bR^n \times \bR^{nd}$, where $\cF \phi_{\xi_1,\ldots,\xi_n}$ denotes the Fourier transform of $\phi_{\xi_1,\ldots,\xi_n}$, i.e.
$$\cF \phi_{\xi_1,\ldots,\xi_n}(\tau_1,\ldots,\tau_n)=\int_{\bR^{n}}
e^{-i (\tau_1 t_1+\ldots+\tau_n t_n)}\phi_{\xi_1,\ldots,\xi_n}(t_1,\ldots,t_n)dt_1 \ldots dt_n.$$

b) If
\begin{equation}
\label{S-norm-finite}
\|S\|_{0,n}^{2}:=\int_{\bR^{nd}} \int_{\bR^{n}} |\cF \phi_{\xi_1,\ldots,\xi_n}(\tau_1,\ldots,\tau_n)|^2   \nu(d\tau_1)\ldots \nu(d\tau_n)\mu(d\xi_1)\ldots \mu(d\xi_n)<\infty,
\end{equation}
then $S \in \cH^{\otimes n}$ and $\|S\|_{\cH^{\otimes n}}^2=\|S\|_{0,n}^2$.

c) Assume in addition that $S(t_1,\cdot, \ldots,t_n,\cdot)=0$ for all $(t_1, \ldots,t_n) \not \in [0,T]^n$, for some $T>0$. If for every $\xi_1,\ldots,\xi_n \in \bR^d$, the function $(t_1,\ldots,t_n) \mapsto \cF S(t_1,\cdot,\ldots,t_n,\cdot)(\xi_1,\ldots,\xi_n)$ is bounded and continuous a.e. on $[0,T]^n$, and
$$I_T(n):=\int_{\bR^{nd}}\int_{[0,T]^{2n}} \prod_{j=1}^{n}\gamma(t_j-s_j)\cF S(t_1,\cdot,\ldots,t_n,\cdot)(\xi)\overline{\cF S(s_1,\cdot,\ldots,s_n,\cdot)(\xi)}d{\bf t}d{\bf s} \mu_n(d \xi)<\infty,$$
then \eqref{S-norm-finite} holds, $S \in \cH^{\otimes n}$ and $\|S\|_{\cH^{\otimes n}}^2=I_T(n)$. In the integral $I_T(n)$ above, ${\bf t}=(t_1,\ldots,t_n)$, ${\bf s}=(s_1,\ldots,s_n)$ and $\mu_n(d\xi_1, \ldots,d\xi_n)=\prod_{j=1}^{n}\mu(d\xi_j)$ is a measure on $\bR^{nd}$.

\end{theorem}

\noindent {\bf Proof:} We argue as in the proof of Theorem \ref{phi-in-H}.
Part a) follows by Fubini's theorem and hypothesis (i).
For b), we note that $a\in \widetilde{L}_{\bC}^2(\bR^{n(d+1)},\Pi_n)$, where $a(\tau_1,\ldots,\tau_n,\xi_1, \ldots,\xi_n)=\cF \phi_{\xi_1,\ldots,\xi_n}(\tau_1,\ldots,\tau_n)$
and $$\Pi_n(d\tau_1,\ldots,d\tau_n,d\xi_1,\ldots,d\xi_n)=\nu(d\tau_1)\ldots \nu(d\tau_n) \mu(d\xi_1)\ldots \mu(d\xi_n).$$

By Lemma \ref{jolis-lemma}, $\cF(\cD(\bR^{n(d+1)}))$ is dense in $\widetilde{L}_{\bC}^2(\bR^{n(d+1)},\Pi_n)$. Hence, for any $\varepsilon>0$, there exists a function $l=l(\varepsilon) \in \cD(\bR^{n(d+1)})$ such that
$$\|\varphi-l\|_{0,n}:=\int_{\bR^{n(d+1)}} |a-\cF l|^2 d\Pi_n<\varepsilon^2.$$
The conclusion follows since $\cH^{\otimes n}$ is the completion of $\cD(\bR^{n(d+1)})$ with respect to the inner product $\langle \cdot,\cdot \rangle_{\cH^{\otimes n}}$ defined by
$$\langle \varphi_1,\varphi_2 \rangle_{\cH^{\otimes n}}=\int_{\bR^{n(d+1)}}
\cF \phi_{\xi_1, \ldots,\xi_n}^{(1)}(\tau_1,\ldots,\tau_n)
\overline{\cF \phi_{\xi_1, \ldots,\xi_n}^{(2)}(\tau_1,\ldots,\tau_n)}\Pi_n(d\tau_1,\ldots,d\tau_n,d\xi_1,\ldots,d\xi_n)$$
where $\phi_{\xi_1, \ldots,\xi_n}^{(k)}(t_1,\ldots,t_n)=
\cF\varphi(t_1,\cdot\ldots,t_n,\cdot)(\xi_1,\ldots,\xi_n)$ for $k=1,2$.

c) For every $\xi_1, \ldots,\xi_n \in \bR^d$ fixed, we apply Lemma \ref{energy-lemma} to the bounded function $\phi_{\xi_1,\ldots,\xi_n}:\bR^n \to \bC$ which is continuous a.e. and has support contained in $[0,T]^n$. We apply this lemma for the measure $m=\nu_n$ and the kernel $\kappa=\gamma_n$ on $\bR^n$, where $\nu_n(d\tau_1, \ldots,d\tau_n)=\prod_{j=1}^{n}\nu(d\tau_j)$ and $\gamma_n(t_1, \ldots,t_n)=\prod_{j=1}^{n}\gamma(t_j)$. We obtain that, for any $\xi_1, \ldots,\xi_n \in \bR^d$,
$$\int_{[0,T]^{2n}}\prod_{j=1}^{n}\gamma(t_j-s_j)\phi_{\xi_1,\ldots,\xi_n}({\bf t})\overline{\phi_{\xi_1,\ldots,\xi_n}({\bf s})}d{\bf t}d{\bf s}=\int_{\bR^n}|\cF \phi_{\xi_1,\ldots,\xi_n}(\tau_1,\ldots,\tau_n)|^2 \nu_n(d\tau),$$
where ${\bf t}=(t_1, \ldots,t_n)$ and ${\bf s}=(s_1,\ldots,s_n)$.
We integrate with respect to $\mu_n(d\xi_1, \ldots,d\xi_n)$. We obtain that
$I_T(n)=\|S\|_{0,n}^2$.
Since $I_T(n)<\infty$, it follows that \eqref{S-norm-finite} holds. The conclusion follows by part b). $\Box$

\vspace{3mm}

As a consequence of the previous theorem, we obtain the following result.

\begin{theorem}
\label{fn-in-Hn}
Suppose that $\mu$ satisfies \eqref{Dalang-cond}. If $d \geq 4$, suppose in addition that Hypothesis A holds.
Then for any $t>0$, $x \in \bR^d$ and $n \geq 1$,
$$f_n(\cdot,t,x) \in \cH^{\otimes n} \quad \mbox{and} \quad \|f_n(\cdot,t,x)\|_{\cH^{\otimes n}}^{2}=I_t(n),$$
where
$$I_t(n)=\int_{\bR^{nd}}\int_{[0,t]^{2n}} \frac{\sin((t_2-t_1)|\xi_1|)}{|\xi_1|}
\cdot \frac{\sin((t_3-t_2)|\xi_1+\xi_2|)}{|\xi_1+\xi_2|}\ldots \frac{\sin((t-t_n)|\xi_1+\ldots+\xi_n|)}{|\xi_1+\ldots+\xi_n|}$$
$$\frac{\sin((s_2-s_1)|\xi_1|)}{|\xi_1|}
\cdot \frac{\sin((s_3-s_2)|\xi_1+\xi_2|)}{|\xi_1+\xi_2|}\ldots \frac{\sin((t-s_n)|\xi_1+\ldots+\xi_n|)}{|\xi_1+\ldots+\xi_n|}$$
$$\prod_{j=1}^{n}\gamma(t_j-s_j)dt_1 \ldots dt_n ds_1 \ldots ds_n \mu(d\xi_1)\ldots \mu(d\xi_n).$$
\end{theorem}

\noindent {\bf Proof:} We apply Theorem \ref{phi-in-Hn}.c) to the function $S=f_n(\cdot,t,x)$ for fixed $t>0$ and $x \in \bR^d$, i.e. $S(t_1, \ldots,t_n)=f_n(t_1,\cdot,\ldots,t_n,\cdot,t,x)$. To see that $f_n(\cdot,t,x)$ satisfies the conditions of this theorem, recall that
for any $(t_1, \ldots,t_n) \in \bR$ and $\xi_1,\ldots,\xi_n \in \bR^d$,
\begin{eqnarray*}
\lefteqn{\phi_{\xi_1,\ldots,\xi_n}(t_1,\ldots,t_n):= \cF f_n(t_1,\cdot,\ldots,t_n,\cdot,t,x)(\xi_1,\ldots,\xi_n)} \\
& & =e^{-i (\xi_1+\ldots+\xi_n)\cdot x}  \overline{\cF G(t_2-t_1,\cdot)(\xi_1)}\,
\overline{\cF G(t_3-t_2,\cdot)(\xi_1+\xi_2)}\ldots \overline{\cF G(t-t_n,\cdot)(\xi_1+\ldots+\xi_n)}\\
& & =e^{-i (\xi_1+\ldots+\xi_n)\cdot x} \frac{\sin((t_2-t_1)|\xi_1|)}{|\xi_1|}
\cdot \frac{\sin((t_3-t_2)|\xi_1+\xi_2|)}{|\xi_1+\xi_2|}\ldots \frac{\sin((t-t_n)|\xi_1+\ldots+\xi_n|)}{|\xi_1+\ldots+\xi_n|}
\end{eqnarray*}
if $0<t_1<\ldots<t_n<t$ and $\phi_{\xi_1,\ldots,\xi_n}(t_1,\ldots,t_n)=0$ otherwise. Hence,
$$|\phi_{\xi_1,\ldots,\xi_n}(t_1, \ldots,t_n)|\leq (t_2-t_1)\ldots(t-t_n)1_{\{0<t_1<\ldots<t_n<t\}}\leq t^n 1_{[0,t]^n}.$$
It follows that $f_n(\cdot,t,x)$ satisfies conditions (i) and (ii) of Theorem \ref{phi-in-Hn}.

Similarly to the calculations done in the proof of Theorem \ref{theorem-sum-finite} below, one can prove that $I_t(n)<\infty$, under condition \eqref{Dalang-cond}. By Theorem \ref{phi-in-Hn}.c), we conclude that $f_n(\cdot,t,x)\in \cH^{\otimes n}$ and $\|f_n(\cdot,t,x)\|_{\cH^{\otimes n}}^2=I_t(n)$. $\Box$

\begin{remark}

{\rm Theorem \ref{phi-in-Hn}.c) can also be applied to the function $S=f_n(\cdot,t,x)$ where $f_n(\cdot,t,x)$ is defined by \eqref{def-fn-1} and $G$ is the fundamental solution of the heat equation, given by \eqref{def-G-heat}. Using the same argument as in the proof of Theorem \ref{fn-in-Hn}, we infer that, if $\mu$ satisfies \eqref{Dalang-cond}, then $f_n(\cdot,t,x) \in \cH^{\otimes n}$ for all $t>0$ and $x \in \bR^d$.}

\end{remark}

\section{Summability of the series}
\label{section-summability}

In this section, we show that under condition \eqref{Dalang-cond}, for any $t>0$ and $x \in \bR^d$,
\begin{equation}
\label{sum-finite}
\sum_{n\geq 1}n!\,\|\widetilde{f}_n(\cdot,t,x)\|_{\cH^{\otimes n}}^2<\infty,
\end{equation}
where $\widetilde{f}_n(\cdot,t,x)$ the symmetrization of $f_n(\cdot,t,x)$ defined as follows. Let $t_1>0,\ldots,t_n>0$ be arbitrary. If $d\leq 2$, for any $x_1,\ldots,x_n \in \bR^d$, we let
$$\widetilde{f}_n(t_1,x_1,\ldots,t_n,x_n,t,x)=\frac{1}{n!}\sum_{\rho \in {\cal P}_n} f_n(t_{\rho(1)},x_{\rho(1)},\ldots,t_{\rho(n)},x_{\rho(n)},t,x),$$
where ${\cal P}_n$ is the set of all permutations of $\{1,\ldots,n\}$. If $d \geq 3$,
we let $\widetilde{f}_n(t_1,\ldots,t_n,\cdot,t,x)$ be the distribution in $\cS'(\bR^{nd})$ defined as follows: for any $\psi \in \cS(\bR^{nd})$,
\begin{equation}
\label{def-fn-tilde}
\big(\widetilde{f}_n(t_1,\cdot,\ldots,t_n,\cdot,t,x),\psi\big)=\frac{1}{n!}\sum_{\rho \in \cP_n}\big(f_n(t_{\rho(1)},\cdot,\ldots,t_{\rho(n)},\cdot,t,x),\psi_{\rho} \big),
\end{equation}
where
\begin{equation}
\label{def-psi-rho}
\psi_{\rho}(x_1,\ldots,x_n)=\psi(x_{\rho^{-1}(1)},\ldots,x_{\rho^{-1}(n)}) \quad \mbox{for all} \quad x_1,\ldots,x_n \in \bR^d.
\end{equation}
It follows that for any $d\geq 1$, the Fourier transform of $\widetilde{f}_n(t_1,\ldots,t_n,\cdot,t,x)$ is the function
\begin{align*}
 \cF \widetilde{f}_n(t_1,\ldots,t_n,\cdot,t,x)(\xi_1,\ldots,\xi_n)&=e^{-i(\xi_1\cdot x_1+\ldots+\xi_n\cdot x_n)}\frac{1}{n!}\sum_{\rho \in {\cal P}_n}\overline{\cF G(t_{\rho(2)}-t_{\rho(1)},\cdot)(\xi_{\rho(1)})}\ldots\\
&\overline{\cF G(t-t_{\rho(n)},\cdot)(\xi_{\rho(1)}+\ldots+\xi_{\rho(n)})}
\, 1_{\{0<t_{\rho(1)}<\ldots<t_{\rho(n)}<t\}}.
\end{align*}

To prove \eqref{sum-finite}, we need some preliminary results. The first result is based on the fact that the function $f$ is {\em non-negative}. (See also relation (3.4) of \cite{dalang-mueller03} for a related result.)

\begin{lemma}
\label{lemmaA}
Let $\mu$ be a tempered measure on $\bR^d$ whose Fourier transform in $\cS_{\bC}'(\bR^d)$ is a locally-integrable function $f:\bR^d \to [0,\infty]$ such that $f(x)<\infty$ a.e.
Then for any $\beta>0$,
\begin{equation}
\label{MaximumPrin}
\sup_{\eta \in \bR^d}\int_{\bR^d}\left(\frac{1}{1+|\xi+\eta|^2} \right)^\beta\mu(d\xi)=\int_{\bR^d}\left(\frac{1}{1+|\xi|^2}\right)^\beta \mu(d\xi).
\end{equation}
\end{lemma}

\noindent {\bf Proof:} We prove the result in a similar way as in Remark 5.8 in \cite{song15}. We assume that the right hand side of \eqref{MaximumPrin} is finite, otherwise it is trivial.  Note that for $c>0$ and $\beta>0$,
\begin{equation}\label{Gamma-eqn}
c^{-\beta}=\frac1{\Gamma(\beta)} \int_0^\infty t^{\beta-1} e^{-ct}dt.
\end{equation}

Fix $\eta \in \bR^d$. We apply \eqref{Gamma-eqn} to $c=1+|\xi+\eta|^2$ and then integrate $\mu(d\xi)$. Using Fubini's theorem, we obtain:
$$\int_{\bR^d}\left(\frac{1}{1+|\xi+\eta|^2} \right)^\beta\mu(d\xi)=\frac{1}{\Gamma(\beta)}\int_0^\infty t^{\beta-1} e^{-t} \left(\int_{\bR^d} e^{-t|\xi+\eta|^2} \mu(d\xi)\right) dt.$$
Let $p_t(x)=(2\pi t)^{-d/2}e^{-|x|^2/(2t)}$. Note that for any $\xi,\eta \in \bR^d$,
$$\cF(e^{-i \eta \cdot} p_{2t})(\xi)=\int_{\bR^d}e^{-i (\xi+\eta) \cdot x}p_{2t}(x)dx=\cF p_{2t}(\xi+\eta)=e^{-t|\xi+\eta|^2}.$$
By applying Parseval's identity \eqref{Parseval-formula} to the function $\varphi=e^{-i \eta \cdot}p_{2t} \in \cS_{\bC}(\bR^d)$, we see that
$$\int_{\bR^d}e^{-i \eta \cdot x}p_{2t}(x)f(x)dx=\frac{1}{(2\pi)^d} \int_{\bR^d} e^{-t|\xi+\eta|^2}\mu(d\xi)$$
Hence, by applying Fubini's theorem,
\begin{eqnarray}
\nonumber
\int_{\bR^d}\left(\frac{1}{1+|\xi+\eta|^2} \right)^\beta\mu(d\xi)
&=&\frac{1}{\Gamma(\beta)}\int_0^\infty t^{\beta-1} e^{-t} \left(\int_{\bR^d} e^{-i\eta\cdot x}p_{2t}(x) f(x)dx \right) dt\\
\label{DM-relation}
&=&\frac{1}{\Gamma(\beta)}\int_{\bR^d}e^{-i \eta \cdot x}G_{d,\beta}(x)f(x)dx,
\end{eqnarray}
 where $G_{d,\beta}$ is the Bessel kernel:
$$G_{d,\beta}(x)=\frac{1}{\Gamma(\beta)} \int_0^{\infty}t^{\beta-1}e^{-t}p_{2t}(x)dt>0.$$
We take the modulus on both sides of \eqref{DM-relation} and we use the fact that the left-hand side of this relation is non-negative. We use the inequality $|\int \ldots| \leq \int |\ldots|$ on the right-hand side. Since $|e^{-i\eta \cdot x}|=1$ and $f$ is non-negative, we obtain that
$$\int_{\bR^d}\left(\frac{1}{1+|\xi-\eta|^2}\right)^\beta \mu(d\xi)  \leq \int_{\bR^d}G_{d,k}(x)f(x)dx=\int_{\bR^d}\left(\frac{1}{1+|\xi|^2}\right)^\beta \mu(d\xi).$$
$\Box$

\vspace{3mm}

Based on the previous lemma, we obtain the following result.

\begin{lemma}
\label{sup-lemma}
For any $t>0$,
\begin{equation}
\label{sup-G}
\sup_{\eta \in \bR^d}\int_{\bR^d}|\cF G(t,\cdot)(\xi+\eta)|^2 \mu(d\xi)\leq 4t^2\int_{\bR^d}\frac{1}{1+t^2|\xi|^2}\mu(d\xi)
\end{equation}
and
\begin{equation}
\label{sup-G-2}
\sup_{\eta \in \bR^d}\int_{\bR^d}|\cF G(t,\cdot)(\xi+\eta)|^2 \mu(d\xi)\leq 2(t^2 \vee 1)\int_{\bR^d}\frac{1}{1+|\xi|^2}\mu(d\xi).
\end{equation}

\end{lemma}

\noindent {\bf Proof:} We first prove \eqref{sup-G}. Note that $\frac{|\sin x|}{x} \leq \frac{2}{1+x}$ for any $x>0$. (This can be seen as follows: if $x \leq 1$, then $\frac{|\sin x|}{x} \leq 1 \leq \frac{2}{1+x}$, and if $x>1$, then $\frac{|\sin x|}{x} \leq \frac{1}{x}\leq \frac{2}{1+x}$.) Hence
$$|\cF G(t,\cdot)(\xi)|^2=\frac{\sin^2(t|\xi|)}{|\xi|^2}\leq \frac{4t^2}{(1+t|\xi|)^2}\leq \frac{4t^2}{1+t^2|\xi|^2}.$$
It follows that
\begin{eqnarray*}
\lefteqn{\sup_{\eta \in \bR^d}\int_{\bR^d}|\cF G(t,\cdot)(\xi+\eta)|^2 \mu(d\xi) \leq
\sup_{\eta \in \bR^d} \int_{\bR^d}\frac{4t^2}{1+t^2|\xi+\eta|^2}\mu(d\xi)= } \\
& & 4t^2 \sup_{\eta \in \bR^d}\int_{\bR^d}\frac{1}{1+|t\xi+\eta|^2} \mu(d\xi)=
4t^2 \sup_{\eta \in \bR^d} \int_{\bR^d}\frac{1}{1+|\xi+\eta|^2}\mu_t(d\xi),
\end{eqnarray*}
where $\mu_t=\mu \circ h_t^{-1}$ and $h_t(\xi)=t\xi$. We now apply Lemma \ref{lemmaA} (with $\beta=1$) to the measure $\mu_t$. To justify the application of this result, we note that the Fourier transform in $\cS'(\bR^d)$ of the measure $\mu_t$ is the non-negative definite function $f_t$ defined by $f_t(x)=f(tx),x \in \bR^d$, since for any $\varphi \in \cS(\bR^d)$ we have:
\begin{eqnarray*}
 \int_{\bR^d}\cF \varphi(\xi)\mu_t(d\xi)&=&\int_{\bR^d}\cF\varphi(t\xi)\mu(d\xi)=
\int_{\bR^d}\cF\varphi^{(t)}(\xi)\mu(d\xi)\\
&=&\int_{\bR^d}\varphi^{(t)}(x)f(x)dx=\int_{\bR^d}\varphi(x)f_t(x)dx,
\end{eqnarray*}
where $\varphi^{(t)}(x)=t^{-d}\varphi(x/t)$.
It follows that
$$\sup_{\eta \in \bR^d} \int_{\bR^d}\frac{1}{1+|\xi+\eta|^2}\mu_t(d\xi)=
\int_{\bR^d}\frac{1}{1+|\xi|^2}\mu_t(d\xi)=
\int_{\bR^d}\frac{1}{1+t^2|\xi|^2}\mu(d\xi).$$

Inequality \eqref{sup-G-2} follows similarly, using \eqref{UB-Fourier-G}.
$\Box$

\vspace{3mm}

We will need the following elementary result.

\begin{lemma}
\label{basic-ineq-lemma}
For any $n \geq 1$ and for any function $h:[0,t]^n \to \bR$ which is either non-negative or integrable,
\begin{equation}
\label{basic-ineq}
\int_{[0,t]^n}\int_{[0,t]^n}\prod_{j=1}^n\gamma(t_j-s_j)h(t_1,\ldots,t_n)dt_1 \ldots dt_n ds_1 \ldots ds_n\leq \Gamma_t^n \int_{[0,t]^n}|h(t_1,\ldots,t_n)|dt_1 \ldots dt_n,
\end{equation}
where $\Gamma_t=\int_{-t}^{t}\gamma(s)ds=2\int_0^t \gamma(s)ds$.
 \end{lemma}

\noindent {\bf Proof:} We consider only the case when $h$ is a non-negative function. The proof for an integrable function $h$ is similar. We use an induction argument on $n \geq 1$. For $n=1$,
we note that $\int_0^t \gamma(r-s)dr=\int_{-s}^{t-s}\gamma(r)dr\leq \Gamma_t$ and hence
$$\int_0^t h(s)\left(\int_0^t \gamma(r-s)dr \right) ds \leq \Gamma_t \int_0^t h(s)ds.$$
For the induction step, we assume that the inequality holds for $n-1$. Then
\begin{eqnarray*}
& &\int_0^t \int_0^t \gamma(t_n-s_n)\left(\int_{[0,t]^{2(n-1)}} h(t_1,\ldots,t_n)\prod_{j=1}^{n-1}\gamma(t_j-s_j)dt_1 ds_1 \ldots dt_{n-1}ds_{n-1} \right)dt_n ds_n \leq \\
& &  \int_0^t \int_0^t \gamma(t_n-s_n)\left(\Gamma_t^{n-1}\int_{[0,t]^{n-1}} h(t_1,\ldots,t_n)dt_1  \ldots dt_{n-1} \right)dt_n ds_n=\\
& &  \Gamma_t^{n-1}\int_{[0,t]^{n-1}} \left(\int_0^t \int_0^t \gamma(t_n-s_n)h(t_1,\ldots,t_n) dt_n ds_n\right) dt_1\ldots dt_{n-1} \leq \\
& &  \Gamma_t^{n-1}\int_{[0,t]^{n-1}} \left(\Gamma_t \int_0^t h(t_1,\ldots,t_n)dt_n \right) dt_1 \ldots dt_{n-1}
\end{eqnarray*}
where we used the induction hypothesis for the first inequality, and  inequality \eqref{basic-ineq} for the case $n=1$ for the last inequality.
For the equality above, we used Fubini's theorem whose application is justified since $h$ is non-negative.  $\Box$

\vspace{3mm}

\begin{theorem}
\label{theorem-sum-finite}
Suppose that $d\geq 1$ is arbitrary and $\mu$ satisfies \eqref{Dalang-cond}. If $d\geq 4$, suppose in addition that Hypothesis A holds. Then relation \eqref{sum-finite} holds for any $t>0$ and $x \in \bR^d$.
\end{theorem}

\noindent {\bf Proof:} We will prove that
\begin{equation}
\label{sum-alpha-finite}
\sum_{n \geq 0}\frac{1}{n!}\alpha_n(t)<\infty,
\end{equation}
where
\begin{equation}
\label{alpha-J}\alpha_n(t)=E|J_n(t,x)|^2=E|I_n(f_n(\cdot,t,x))|^2=(n!)^2 \|\widetilde{f}_n(\cdot,t,x)\|_{\cH^{\otimes n}}^2.
\end{equation}

For this, we proceed as in the proof of Theorem 3.2 of \cite{HHNT}. In the integrals below, we use the notation ${\bf t}=(t_1, \ldots,t_n)$, ${\bf s}=(s_1, \ldots,s_n)$,
${\bf x}=(x_1, \ldots,x_n)$ and ${\bf y}=(y_1, \ldots,y_n)$.
Then
\begin{equation}
\label{def-alpha-n}
\alpha_n(t)=\int_{[0,t]^{2n}}\prod_{j=1}^{n}\gamma(t_j-s_j)\psi_n({\bf t},{\bf s})d{\bf t}d{\bf s},
\end{equation}
where
\begin{eqnarray*}
\psi_n({\bf t},{\bf s})&=&\int_{\bR^{nd}}\cF g_{\bf t}^{(n)}(\cdot,t,x)(\xi_1,\ldots,\xi_n)\overline{\cF g_{\bf s}^{(n)}(\cdot,t,x)(\xi_1,\ldots,\xi_n)}\mu(d\xi_1)\ldots \mu(d\xi_n)
\end{eqnarray*}
and
\begin{equation}
\label{gtn}
g_{\bf t}^{(n)}(\cdot,t,x)=n!\widetilde{f}_n(t_1,\cdot,\ldots,t_n,\cdot,t,x).
\end{equation}

If the permutation $\rho$ of $\{1, \ldots,n\}$ is chosen such that
$t_{\rho(1)}<\ldots<t_{\rho(n)}$,
then
\begin{eqnarray}
\nonumber
\cF g_{\bf t}^{(n)}(\xi_1,\ldots,\xi_n)&=&e^{-i \sum_{j=1}^{n}\xi_j\cdot x} \overline{\cF G(t_{\rho(2)}-t_{\rho(1)},\cdot)(\xi_{\rho(1)})} \  \overline{\cF G(t_{\rho(3)}-t_{\rho(2)},\cdot)(\xi_{\rho(1)}+\xi_{\rho(2)})} \\
\label{Fourier-gn}
& & \ldots \overline{\cF G(t-t_{\rho(n)},\cdot)(\xi_{\rho(1)}+\ldots+\xi_{\rho(n)})}
\end{eqnarray}

By the Cauchy-Schwarz inequality  and the
inequality $ab \leq (a^2+b^2)/2$, we obtain:
$$\psi_n({\bf t},{\bf s})\leq \psi_n({\bf t},{\bf t})^{1/2}\psi_n({\bf s},{\bf s})^{1/2} \leq \frac{1}{2}\Big(\psi_n({\bf t},{\bf t})+\psi_n({\bf s},{\bf s})\Big).$$
Using \eqref{def-alpha-n} and the symmetry of the function $\gamma$, it follows that
\begin{eqnarray*}
\alpha_n(t) \leq  \int_{[0,t]^{2n}}\prod_{j=1}^{n}\gamma(t_j-s_j)
\frac{\psi_n({\bf t},{\bf t})+\psi_n({\bf s},{\bf s})}{2}d{\bf t}d{\bf s}
=\int_{[0,t]^{2n}}\prod_{j=1}^{n}\gamma(t_j-s_j)
\psi_n({\bf t},{\bf t})d{\bf t}d{\bf s}.
\end{eqnarray*}

\noindent Using Lemma \ref{basic-ineq-lemma} for the function $h({\bf t})=\psi_n({\bf t},{\bf t})$, we obtain:
\begin{equation}
\label{bound1-alpha}
\alpha_n(t)\leq \Gamma_t^n \int_{[0,t]^n}\psi_n({\bf t},{\bf t})d{\bf t}.
\end{equation}

We now estimate $\psi_n({\bf t},{\bf t})$. We denote $u_j=t_{\rho(j+1)}-t_{\rho(j)}$ for $j=1,\ldots,n$, where $t_{\rho(n+1)}=t$. We have:
\begin{eqnarray*}
\psi_n({\bf t},{\bf t})
&=&   \int_{\bR^{nd}} |\cF G(u_1,\cdot)(\xi_{\rho(1)})|^2
\ |\cF G(u_2,\cdot)(\xi_{\rho(1)}+\xi_{\rho(2)})|^2 \ldots \\
& & \quad \quad \quad \quad \quad |\cF G(u_n,\cdot)(\xi_{\rho(1)}+\ldots+\xi_{\rho(n)})|^2 \mu(d\xi_1) \ldots \mu(d\xi_n) \\
&=&   \int_{\bR^{d}}\mu(d\xi_1') |\cF G(u_1,\cdot)(\xi_{1}')|^2 \left(
\int_{\bR^d} \mu(d\xi_2') |\cF G(u_2,\cdot)(\xi_{1}'+\xi_{2}')|^2 \ldots \right. \\
& & \left.\quad \quad \quad \quad \quad \left(\int_{\bR^d} |\cF G(u_n,\cdot)(\xi_{1}'+\ldots+\xi_{n}')|^2  \mu(d\xi_n') \right) \ldots \right),
\end{eqnarray*}
where for the last equality we used the change of variable $\xi_j'=\xi_{\rho(j)}$ for $j=1,\ldots,n$.
Using Lemma \ref{sup-lemma}
it follows that
\begin{eqnarray}
\label{bound-psi}
\psi_n({\bf t},{\bf t})  \leq    \prod_{j=1}^{n}\left(\sup_{\eta \in \bR^d} \int_{\bR^d} |\cF G(u_j,\cdot)(\xi_j+\eta)|^2 \mu(d\xi_j)\right)
\leq   \prod_{j=1}^{n}\int_{\bR^d} \frac{4u_j^2}{1+u_j^2|\xi_j|^2} \mu(d\xi_j).
\end{eqnarray}

We now go back to the estimate \eqref{bound1-alpha} for $\alpha_n(t)$. We decompose the set $[0,t]^n$ into $n!$ disjoint regions of the form $t_{\rho(1)}<\ldots<t_{\rho(n)}$ with $\rho \in {\cal P}_n$.
Using \eqref{bound-psi}, it follows that
\begin{eqnarray}
\nonumber
\alpha_n(t)&\leq &  \Gamma_t^n  \sum_{\rho \in {\cal P}_n} \int_{t_{\rho(1)}<\ldots<t_{\rho(n)}} \int_{\bR^{nd}} \prod_{j=1}^{n}
\frac{4(t_{\rho(j+1)}-t_{\rho(j)})^2}{1+(t_{\rho(j+1)}-t_{\rho(j)})^2|\xi_j|^2}\mu(d\xi_1) \ldots \mu(d\xi_n)d{\bf t} \\
\nonumber
&= & \Gamma_t^n  n! \int_{0<t_{1}<\ldots<t_{n}<t} \int_{\bR^{nd}} \prod_{j=1}^{n}
\frac{4(t_{j+1}-t_{j})^2}{1+(t_{j+1}-t_{j})^2|\xi_j|^2}\mu(d\xi_1) \ldots \mu(d\xi_n)d{\bf t}\\
\nonumber
&=& \Gamma_t^n n! \,\int_{\bR^{nd}}  \int_{S_{t,n}} \prod_{j=1}^{n}
\frac{4w_j^2}{1+w_j^2|\xi_j|^2} d{\bf w}\mu(d\xi_1) \ldots \mu(d\xi_n)\\
\label{def-In(t)}
&=:& \Gamma_t^n n!\,I^{(n)}(t)
\end{eqnarray}
where $S_{t,n}=\{(w_1,\ldots,w_n) \in [0,t]^n;w_1+\ldots+w_n \leq t\}$ and ${\bf w}=(w_1,\ldots,w_n)$.
As in the proof of Lemma 3.3 of \cite{HHNT}, since $S_{t,n} \subset S_{t}^I \times S_t^{I^c}$, the integral $I^{(n)}(t)$ is smaller than
\begin{eqnarray}
\label{def-Jn(t)}
J^{(n)}(t) &:=& \sum_{I \subset \{1,\ldots,n\}} \int_{\bR^{d|I|}} \prod_{j \in I}1_{\{|\xi_j|\leq N\}} \left(\int_{S_t^I} \prod_{j\in I} \frac{4w_j^2}{1+w_j^2|\xi_j|^2} d{\bf w}_I\right) \prod_{j \in I}\mu(d\xi_j)\\
\nonumber
& & \int_{\bR^{d|I^c|}}\prod_{j \in I^c}1_{\{|\xi_j|>N\}} \left(\int_{S_t^{I^c}}
\prod_{j\in I^c} \frac{4w_j^2}{1+w_j^2|\xi_j|^2} d{\bf w}_{I^c}\right)\prod_{j \in I^c}\mu(d\xi_j),
\end{eqnarray}
where $S_t^I=\{{\bf w}_I=(w_j)_{j \in I}; w_j \geq 0, \sum_{j \in I}w_j \leq t\}$ and $S_t^{I^c}=\{{\bf w}_{I^c}=(w_j)_{j \in I^c}; w_j \geq 0, \sum_{j \in I^c}w_j \leq t\}$. Here $|I|$ is the cardinality of $I$ and $N>0$ is arbitrary.

For the integral over the set $S_t^I$ we use the bound
$$\frac{4w_j^2}{1+w_j^2|\xi_j|^2} \leq 4w_j^2 \leq 4t^2,$$
and so, this integral is bounded by $(4t^2)^{|I|} \int_{S_t^I}d{\bf w}_{I}=4^{|I|}t^{3|I|}/|I|!$. For the integral over $S_t^{I^c}$, we have:
$$\int_{S_t^{I^c}}
\prod_{j\in I^c}\frac{4w_j^2}{1+w_j^2|\xi_j|^2}d{\bf w}_{I^c}\leq \prod_{j \in I^c}\int_0^t \frac{4w_j^2}{1+w_j^2|\xi_j|^2} dw_j
 \leq \prod_{j \in I^c}\int_0^t \frac{4}{|\xi_j|^2} dw_j =4^{|I^c|}t^{|I^c|} \prod_{j \in I^c}\frac{1}{|\xi_j|^2}.$$

We denote
$$C_N=\int_{\{|\xi|>N\}}\frac{1}{|\xi|^2}\mu(d\xi) \quad \mbox{and} \quad
D_N=\int_{\{|\xi|\leq N\}}\mu(d\xi).$$
It follows that
\begin{eqnarray}
\nonumber
J^{(n)}(t) &\leq & 4^n \sum_{I \subset \{1,\ldots,n\}} \frac{t^{3|I|}}{|I|!} D_N^{|I|}\cdot t^{|I^c|} C_N^{|I^c|}
= 4^n \sum_{k=0}^{n}\left(\begin{array}{c}n \\
k \end{array} \right) \frac{t^{3k}}{k!}D_N^k t^{n-k}C_N^{n-k} \\
\label{def-Kn(t)}
& \leq & 4^n \sum_{k=0}^{n}2^n \frac{t^{n+2k}}{k!}D_N^k C_N^{n-k}=:K^{(n)}(t).
\end{eqnarray}
Hence
\begin{equation}
\label{alpha-bound}
\alpha_n(t)\leq  \Gamma_t^n  n! \frac{1}{(2\pi)^{nd}} 8^n \sum_{k=0}^{n}\frac{t^{n+2k}}{k!}D_N^k C_N^{n-k}
\end{equation}
and
\begin{eqnarray*}
\sum_{n \geq 0}\frac{1}{n!}\alpha_n(t)&\leq &  \sum_{n \geq 0}\Gamma_t^n 8^n \sum_{k=0}^{n}\frac{t^{n+2k}}{k!}D_N^k C_N^{n-k}
= \sum_{k \geq 0}\frac{t^{2k}}{k!}D_N^k C_{N}^{-k}\sum_{n \geq k}(8C_N \Gamma_t t)^n \\
&=&   \sum_{k \geq 0}\frac{t^{2k}}{k!}D_N^k C_{N}^{-k}(8C_N \Gamma_t  t)^{k} \sum_{n \geq 0}(8C_N \Gamma_t t)^n.
\end{eqnarray*}
Due to condition \eqref{Dalang-cond}, $C_N \to 0$ as $N \to \infty$. Hence, there exists $N_t>0$ such that $8C_N \Gamma_t t<1$ for all $N > N_t$. We choose $N>N_t$ arbitrary. We have:
$$
\sum_{n \geq 0}\frac{1}{n!}\alpha_n(t) \leq  \frac{1}{1-8C_N \Gamma_t t} \sum_{k \geq 0}\frac{1}{k!}(8D_N \Gamma_t t^3)^k
= \frac{1}{1-8C_N \Gamma_tt} \exp\Big(8D_N \Gamma_t t^3 \Big)<\infty.$$
This concludes the proof of \eqref{sum-alpha-finite}.

\begin{remark}
{\rm In the proof of Theorem \ref{theorem-sum-finite}, we expressed $\alpha_n(t)$ as an integral which depends on the measure $\mu$ (instead of the kernel $f$); see \eqref{def-alpha-n}. However, the fact that the Fourier transform  of $\mu$ is the {\em locally integrable non-negative function $f$} was used in Lemma \ref{lemmaA}.}
\end{remark}

\section{Existence of the solution}
\label{section-existence}

In this section, we show that the process $u=\{u(t,x);t\geq 0,x \in \bR^d\}$  defined by:
\begin{equation}
\label{def-u}
u(t,x)=1+\sum_{n\geq 1}I_n(f_n(\cdot,t,x)), \quad t>0,x\in \bR^d
\end{equation}
is a solution of \eqref{wave}. (Note that
the series above converges in $L^2(\Omega)$, due to \eqref{sum-finite}.)

The results of Sections 3 and 4 (in particular, Theorem 3.5 which gives criteria for integrability, and Theorem 4.4 which gives the summability of the series) play a crucial role in the proof of the existence of solution. Generally speaking, a solution of \eqref{wave} is a process $u$ which satisfies the equation:
\begin{equation}
\label{def-sol}
u(t,x) = 1 + \int_{0}^{t}\int_{\bR^d}G(t-s,x-y)u(s,y)W(\delta s,\delta y),
\end{equation}
where the integral is interpreted in the Skorohod sense. To give a rigorous meaning to this equation, we need to discuss separately the cases $d\leq 2$ and $d\geq 3$.

\subsection{The case $d\leq 2$}

In this case, the solution is defined exactly as in the case of the Parabolic Anderson Model.

\begin{definition}
\label{definition-solution}
{\rm Assume that $d\leq 2$. A square-integrable process $u=\{u(t,x);t\geq 0,x\in \bR^d\}$ with $u(0,x)=1$ for all $x \in \bR^d$ and Wiener chaos expansion
\begin{equation}
\label{def-u-kernel-k}
u(t,x)=1+\sum_{n\geq 1}I_n(k_n(\cdot,t,x)), \quad t>0,x\in \bR^d,
\end{equation}
for some symmetric non-negative functions $k_n(\cdot,t,x) \in \cH^{\otimes n}$, is a {\bf solution} to equation \eqref{wave} if the following conditions are satisfied:
\begin{description}
\item[(a)] $u$ has a jointly  measurable modification (denoted also by $u$) and \begin{equation}
    \label{moment2-u}
    \sup_{(t,x)\in [0,T] \times \bR^d}E|u(t,x)|^2<\infty \quad \mbox{for all} \ T>0;
    \end{equation}

\item[(b)] for any $t>0$ and $x \in \bR^d$, $v^{(t,x)} \in {\rm Dom} \ \delta$ and $u(t,x)=1+\delta(v^{(t,x)})$ in $L^2(\Omega)$, where
\begin{equation}
\label{def-v}
v^{(t,x)}(s,y)=1_{(0,t)}(s)G(t-s,x-y)u(s,y), \quad s \geq 0,y \in \bR^d.
\end{equation}
\end{description}
 }
\end{definition}

The proof of the existence of the solution is identical to the parabolic case. We include it for the sake of completeness.

\begin{theorem}
\label{d2-exist-sol}
Suppose that $d\leq 2$ and condition \eqref{Dalang-cond} holds. Then the process $\{u(t,x);t \geq 0,x \in \bR^d\}$ given by \eqref{def-u} is a solution of equation \eqref{wave}.
\end{theorem}

\noindent {\bf Proof:} We apply Proposition \ref{corr-prop-25-0} (Appendix A) to the process $v^{(t,x)}$. For $s\geq 0,y \in \bR^d$,
$$v^{(t,x)}(s,y)=\sum_{n\geq 0}I_n(g_n^{(t,x)}(\cdot,s,y)) \quad \mbox{in} \quad L^2(\Omega),$$
where $g_n^{(t,x)}(\cdot,s,y)=1_{(0,t)}(s)G(t-s,x-y)f_n(\cdot,s,y)$.
We use Remark \ref{corr-prop-25-remark} to verify hypothesis {\em (i)} of Proposition \ref{corr-prop-25-0}. By Theorem \ref{exist-th} below, $u$ is $L^2(\Omega)$-continuous and satisfies \eqref{moment2-u}. By Theorem 30, Chapter IV of \cite{dellacherie-meyer75}, $u$ has a jointly measurable modification. We work with this modification. It follows that $v^{(t,x)}$ is also jointly measurable. Note that
\begin{align*}
& E \left[\int_{(\bR_{+} \times \bR^d)^2} \gamma(s-r) f(y-z) |v^{(t,x)}(s,y) v^{(t,x)}(r,z)|dsdy drdz\right]=\\
& \int_{((0,t) \times \bR^d)^2} \gamma(s-r) f(y-z) G(t-s,x-y)G(t-r,x-z)E|u(s,y) u(r,z)|dsdy drdz<\infty,
\end{align*}
since by Cauchy-Schwarz inequality and \eqref{moment2-u}, for any $s \in (0,t)$, $y \in \bR^d$ and $z \in \bR^d$, $E|u(s,y) u(r,z)|\leq \big(E|u(s,y)|^2 \big)^{1/2} \big(E|u(r,z)|^2\big)^{1/2} \leq C_t$, where $C_t>0$ is a constant depending on $t$.
Hypothesis {\em (ii)} of Proposition \ref{corr-prop-25-0} holds since $g_n^{(t,x)}=f_{n+1}(\cdot,t,x)\in \cH^{\otimes(n+1)}$ by Theorem \ref{fn-in-Hn}. By Theorem \ref{theorem-sum-finite}, $V^{(t,x)}:=\sum_{n\geq 0}I_{n+1}(\widetilde{f}_{n+1}(\cdot,t,x))$ converges in $L^2(\Omega)$.

By Proposition \ref{corr-prop-25-0}, it follows that $v^{(t,x)} \in {\rm Dom} \ \delta$ and $\delta(v^{(t,x)})=V^{(t,x)}$. On the other hand, by \eqref{def-u}, $u(t,x)=1+V^{(t,x)}$ in $L^2(\Omega)$. Hence $u(t,x)=1+\delta(v^{(t,x)})$ in $L^2(\Omega)$. $\Box$

\subsection{The case $d\geq 3$}

In the case $d\geq 3$, the definition of solution proposed in \cite{balan12} is incorrect since the product between the distribution $G(t-s,x-\cdot)$ and the (random) function $u(s,\cdot)$ is not well-defined.
We give below a new definition of the solution, and we prove the existence of a solution. With this new definition, we will also be able to prove the uniqueness of the solution in Section \ref{section-uniqueness} below. In this section, we assume that
Hypothesis A holds.

If $u=\{u(\varphi);\varphi \in \cS(\bR^d)\}$ and $v=\{v(\varphi);\varphi \in \cS(\bR^d)\}$ are two collections of random variables defined on the same probability space $(\Omega,\cF,P)$, we say that $u$ is a {\em modification} of $v$ if $u(\varphi)=v(\varphi)$ a.s. for any $\varphi \in \cS(\bR^d)$.

A distribution $F \in \cS'(\bR^{nd})$ is {\em symmetric} if $\big(F,\psi\big)=\big(F,\psi_{\rho}\big)$ for any for any $\psi \in \cS(\bR^{nd})$ and $\rho \in {\cal P}_n$, where $\psi_{\rho}$ is defined by \eqref{def-psi-rho}.

\begin{definition}
\label{new-def-solution}
{\rm Assume that $d\geq 3$. A square-integrable process $u=\{u(t,x);t\geq 0,x\in \bR^d\}$ with $u(0,x)=1$ for all $x \in \bR^d$ and Wiener chaos expansion \eqref{def-u-kernel-k}
with $k_n(\cdot,t,x) \in \cH^{\otimes n}$, is a {\bf solution} to equation \eqref{wave} if the following conditions are satisfied:
\begin{description}
\item[(a)] for any $s>0$ and $t_1>0,\ldots,t_n>0$, we have: (i) $k_n(t_1,\cdot,\ldots,t_n,\cdot,s,x)$ is a symmetric distribution in $\cS'(\bR^{nd})$ for any $x \in \bR^d$; (ii) for any $\psi \in \cS(\bR^{nd})$ the function $$\bR^d \ni x \mapsto \big(k_n(t_1,\cdot,\ldots,t_n,\cdot,s,x),\psi \big):=h_{\psi}(x) \quad \mbox{is in $\cS(\bR^d)$};$$
 (iii) the map $\psi \mapsto h_{\psi}$ is continuous from $\cS(\bR^{nd})$ to $\cS(\bR^d)$;

\item[(b)] for any $t>0$, $x \in \bR^d$, $s>0$ and $\varphi \in \cS(\bR^d)$, $S_{n}^{(t,x,s,\varphi)}\in \cH^{\otimes n}$ and the series
$\sum_{n \geq 1} I_n(S_{n}^{(t,x,s,\varphi)})$
converges in $L^2(\Omega)$, where $S_{n}^{(t,x,s,\varphi)}$ is defined in Remark \ref{def-gn-tx} below; for any $t_1, \ldots,t_n$, the Fourier transform of  $S_{n}^{(t,x,s,\varphi)}(t_1,\cdot,\ldots,t_n,\cdot)$ is a function which has a version such that
$(t_1,\ldots,t_n,\xi_1,\ldots,\xi_n)\mapsto \cF S_{n}^{(t,x,s,\varphi)}(t_1,\cdot,\ldots,t_n,\cdot)(\xi_1,\ldots,\xi_n)$ is measurable (see Remark \ref{def-version} for the definition of a version of a Fourier transform);

\item[(c)] for any $t>0$, $x \in \bR^d$ and $s>0$, the process $\{(v^{(t,x)}(s,\cdot),\varphi);\varphi \in \cS(\bR^d)\}$ defined by $$(v^{(t,x)}(s,\cdot),\varphi)=(G(t-s,x-\cdot),\varphi)+\sum_{n \geq 1} I_n(S_{n}^{(t,x,s,\varphi)}),$$
    has a modification with values in $\cS'(\bR^d)$ which satisfies hypotheses {\em (i)}-{\em (iv)} of Proposition \ref{corr-prop-25-2} (Appendix A);

\item[(d)] for any $t>0,x \in \bR^d$, $v^{(t,x)} \in {\rm Dom} \ \delta$ and  $u(t,x)=1+\delta(v^{(t,x)})$ in $L^2(\Omega)$.
\end{description}
}
\end{definition}

\begin{remark}
\label{def-gn-tx}
{\rm For any $t_1>0,\ldots,t_n>0,s>0$ and $\psi \in \cS(\bR^{nd})$, we consider the product between the distribution
$1_{(0,t)}(s)G(t-s,x-\cdot)$ and the function $y \mapsto \big(k_n(t_1,\cdot,\ldots,t_n,\cdot,s,y),\psi \big)$. This product, which we denote by $\big(g_n^{(t,x)}(t_1,\cdot,\ldots,t_n,\cdot,s,*),\psi\big)$, is a distribution in $\cS'(\bR^d)$ whose action on a test function $\varphi \in \cS(\bR^d)$ is given by:
$$\Big(\big(g_n^{(t,x)}(t_1,\cdot,\ldots,t_n,\cdot,s,*),\psi\big),\varphi\Big)= 1_{(0,t)}(s)\Big(G(t-s,x-\cdot),\varphi \big(k_n(t_1,\cdot,\ldots,t_n,\cdot,s,*),\psi \big)\Big).$$
Due to condition (a).(iii) in Definition \ref{new-def-solution}, for each $\varphi \in \cS(\bR^d)$ fixed,
the map
$$\cS(\bR^{nd})\ni \psi \mapsto \Big(\big(g_n^{(t,x)}(t_1,\cdot,\ldots,t_n,\cdot,s,*),\psi\big),\varphi\Big)$$
is a distribution in $\cS'(\bR^{nd})$ which
we denote by $S_{n}^{(t,x,s,\varphi)}(t_1,\cdot,\ldots,t_n,\cdot)$. This distribution is symmetric since $k_n(t_1,\cdot,\ldots,t_n,\cdot,s,x)$ is symmetric for any $x \in \bR^d$.}
\end{remark}

\begin{remark}
{\rm
Intuitively, the process $v^{(t,x)}$ given by Definition \ref{new-def-solution}.(c)
should satisfy $v^{(t,x)}(s,\cdot)=G(t-s,x-\cdot)u(s,\cdot)$ where $u$ is a solution. Since $u(s,\cdot)$ may not be a smooth function, this product is not well-defined. For this reason, we define the process $v^{(t,x)}$ using its Wiener chaos expansion. }
\end{remark}

The following theorem establishes the existence of the solution. This result is a correction to Theorem 2.8 of \cite{balan12} whose proof is incorrect since the claim on  page 14, line 15 (that the convolution of the distribution $G(t-s,\cdot)$ with the infinite series $\sum_{n\geq 0}\phi J_n(s,\cdot)$ is equal to the series $\sum_{n \geq 0}(\phi J_n(s,\cdot)*G(t-s,\cdot))(x)$) cannot be justified.

\begin{theorem}
Suppose that $d\geq 3$, $\mu$ satisfies \eqref{Dalang-cond} and Hypothesis A holds. Then the process $\{u(t,x);t \geq 0,x \in \bR^d\}$ given by \eqref{def-u} is a solution of equation \eqref{wave}.
\end{theorem}

\noindent {\bf Proof:} We show that $u$ satisfies the conditions of Definition \ref{new-def-solution} with $k_n(\cdot,t,x)=\widetilde{f}_n(\cdot,t,x)$.

{\em Step 1. (Verification of condition (a))} From the definition of $\widetilde{f}_n(\cdot,t,x)$, we see that it is enough to show that (a) holds when $k_n=f_n$.
Property (a).(i) is clear. To prove (a).(ii), note that $h_{\psi}(x)=H_{\psi}(x,\ldots,x)$ for all $x\in \bR^d$, where
\begin{align*}
H_{\psi}(x_1,\ldots,x_n)&:=1_{\{0<t_1<\ldots<t_n<s\}}\int_{\bR^{nd}}e^{-i(\xi_1\cdot x_1+\ldots+\xi_n\cdot x_n)} \overline{\cF G(t_2-t_1,\cdot)(\xi_1)} \ldots \\
&  \quad \quad \quad \quad \quad \quad \quad \quad \overline{\cF G(s-t_n,\cdot)(\xi_1+\ldots+\xi_n)} \, \cF^{-1}\psi(\xi_1,\ldots,\xi_n)d\xi_1, \ldots d\xi_n,
\end{align*}
for any $x_1, \ldots,x_n \in \bR^d$ and $\psi \in \cS(\bR^{nd})$. The function $$F(\xi_1, \ldots,\xi_n):=\frac{\sin((t_2-t_1)|\xi_1|)}{|\xi_1|}\ldots \frac{\sin((s-t_n)|\xi_1+\ldots+\xi_n|)}{|\xi_1+\ldots+\xi_n|}, \quad \xi_1,\ldots,\xi_n \in \bR^d$$
is infinitely differentiable on $\bR^{nd}$ and has bounded partial derivatives of any order. Hence the product $F \cF^{-1}\psi$ is a function in $\cS(\bR^{nd})$, and therefore $H_{\psi}=1_{\{0<t_1<\ldots<t_n<s\}}\cF(F \cF^{-1}\psi)$ is in $\cS(\bR^{nd})$. Consequently, $h_{\psi} \in \cS(\bR^d)$. 
To prove (a).(iii), note that if $\psi_k \to \psi$ in $\cS(\bR^{nd})$, then $F \cF^{-1}\psi_k \to F \cF^{-1}\psi$ in $\cS(\bR^{nd})$
and $H_{\psi_k} \to H_{\psi}$ in $\cS(\bR^{nd})$, by the continuity of the Fourier transform. From this, we deduce that $h_{\psi_k} \to h_{\psi}$ in $\cS(\bR^d)$. This proves (a).(iii).

In the argument for (b) below, we will need the following fact. Using the change of variables $\eta_k=\xi_1+\ldots+\xi_k,k=1, \ldots,n$, we see that
\begin{equation}
\label{h-psi-Fourier}
h_{\psi}=1_{\{0<t_1<\ldots<t_n<s\}}\cF(\overline{\cF G(s-t_n,\cdot)} V_{\psi}),
\end{equation}
 where
\begin{align}
\nonumber
V_{\psi}(\eta_n)&=\int_{\bR^{n(d-1)}} \overline{\cF G(t_2-t_1,\cdot)(\eta_1)}\ldots \overline{\cF G(t_n-t_{n-1},\cdot)(\eta_{n-1})} \\
\label{def-V-psi}
 & \quad \quad \quad \cF^{-1}\psi(\eta_1,\eta_2-\eta_1,\ldots,\eta_n-\eta_{n-1}) d\eta_1 \ldots d\eta_{n-1}, \quad \quad \quad \eta_n \in \bR^d.
\end{align}
Since $h_{\psi} \in \cS(\bR^d)$, it follows that $\overline{\cF G(s-t_n,\cdot)}V_{\psi} \in \cS(\bR^d)$.

\vspace{3mm}

{\em Step 2. (Verification of condition (b)).}
We fix $t>0,x \in \bR^d,s \in (0,t), \varphi \in \cS(\bR^d)$.
We first show that condition (b) holds when $k_n=f_n$. Since $t,x,s$ and $\varphi$ are fixed in this step, we denote $S_n^{(t,x,s,\varphi)}$ by $S_n$ for simplicity. We prove that $S_n \in \cH^{\otimes n}$ using Theorem \ref{phi-in-Hn}.b).

By definition, for any $\psi \in \cS(\bR^{nd})$, we have:
\begin{equation}
\label{def-Sn}
\big(S_n(t_1,\cdot,\ldots,t_n,\cdot),\psi\big)=
\Big(G(t-s,x-\cdot),\varphi \big(f_n(t_1,\cdot,\ldots,t_n,\cdot,s,*),\psi\big)\Big).
\end{equation}

We treat first the case $n=1$. By definition, for any $\psi \in \cS(\bR^d)$,
$$\big(S_1(t_1,\cdot),\psi\big)=\big(G(t-s,\cdot)*\varphi h_{\psi}\big)(x),$$
where $h_{\psi}(y)=\big(f_1(t_1,\cdot,s,y),\psi \big)$. We first show that the Fourier transform of $S_1(t_1,\cdot)$ in $\cS'(\bR^d)$ is a function and we identify this function. By \eqref{def-G*phi}, for any $\psi \in \cS(\bR^d)$, we have:
\begin{align*}
\big(\cF S_1(t_1,\cdot), \psi\big)&= \big( S_1(t_1,\cdot), \cF \psi\big)=\big(G(t-s,\cdot)*\varphi h_{\cF \psi} \big)(x)\\
&=\int_{\bR^d} e^{-i \xi \cdot x} \, \frac{\sin((t-s)|\xi|)}{|\xi|} \, \cF^{-1}(\varphi h_{\cF \psi})(\xi)\,d\xi\\
&=\int_{\bR^d} e^{-i \xi \cdot x} \, \frac{\sin((t-s)|\xi|)}{|\xi|} \, \big( \cF^{-1}\varphi * \cF^{-1} h_{\cF \psi} \big)(\xi) \,d\xi.
\end{align*}
Note that by the definition of $f_1(t_1,\cdot,s,y)$ and \eqref{def-G*phi},
\begin{align*}
& h_{\cF \psi}(y)=\big(f_1(t_1,\cdot,s,y),\cF \psi \big)=1_{(0,s)}(t_1) \big(G(s-t_1,\cdot)*\cF \psi \big)(y)=\\
& \quad \quad \quad 1_{(0,s)}(t_1)\int_{\bR^d}e^{-i \xi \cdot y} \, \frac{\sin((s-t_1)|\xi|)}{|\xi|}\, \psi(\xi)d\xi=1_{(0,s)}(t_1) \cF \big(\overline{\cF G(s-t_1,\cdot)}\,\psi\big),
\end{align*}
and hence $\cF^{-1} h_{\cF \psi}=1_{(0,s)}(t_1) \overline{\cF G(s-t_1,\cdot)}\,\psi$.
Therefore,
\begin{align*}
\big(\cF S_1(t_1,\cdot), \psi\big)&=1_{(0,s)}(t_1) \int_{\bR^d}e^{-i \xi \cdot x}
\,\frac{\sin((t-s)|\xi|)}{|\xi|}\left(\int_{\bR^d}\cF^{-1}\varphi(\xi-\xi_1) \frac{\sin((s-t_1)|\xi_1|)}{|\xi_1|}\psi(\xi_1)d\xi_1\right) d\xi\\
&=1_{(0,s)}(t_1) \int_{\bR^d}\psi(\xi_1) \frac{\sin((s-t_1)|\xi_1|)}{|\xi_1|}
\left(\int_{\bR^d}e^{-i\xi \cdot x}\frac{\sin((t-s)|\xi|)}{|\xi|}\cF^{-1} (\varphi e^{-i\xi_1 \cdot})d\xi\right)d\xi_1\\
&=1_{(0,s)}(t_1) \int_{\bR^d}\psi(\xi_1) \frac{\sin((s-t_1)|\xi_1|)}{|\xi_1|} \big( G(t-s,\cdot)*\varphi e^{-i \xi_1 \cdot}\big)(x)d\xi_1,
\end{align*}
where for the second equality we used Fubini's theorem and the fact that
\begin{equation}
\label{Fourier-phi-a}
\cF^{-1}\varphi(\xi-a)=\frac{1}{(2\pi)^d}\int_{\bR^d}e^{i(\xi-a)\cdot x}\varphi(x)dx=\cF^{-1}(\varphi e^{-i a \cdot})(\xi),
\end{equation}
and for the last equality we used
\eqref{def-G*phi}. This proves that (a version of) the Fourier transform of $S_1(t_1,\cdot)$ is the function
$$\cF S(t_1,\cdot)(\xi_1)=1_{(0,s)}(t_1) \frac{\sin((s-t_1)|\xi_1|)}{|\xi_1|} \big(G(t-s,\cdot)*\varphi e^{-i \xi_1 \cdot}\big)(x), \quad \xi_1 \in \bR^d.$$

To prove that $S_1 \in \cH$, we apply Theorem \ref{phi-in-H}.b).
The function $(t_1,\xi_1)\mapsto \phi_{\xi_1}(t_1):=\cF S_1(t_1,\cdot)(\xi_1)$ is measurable by Fubini's theorem. The function $t_1 \mapsto \phi_{\xi_1}(t_1)$ is integrable on $\bR$ and
its Fourier transform is
$$\cF \phi_{\xi_1}(\tau_1)=\int_{\bR^d}e^{-i\tau_1 t_1}\phi_{\xi_1}(t_1)dt_1=\big(G(t-s,\cdot)*\varphi e^{-i \xi_1 \cdot}\big)(x)\int_0^s e^{-i \tau_1 t_1} \frac{\sin((s-t_1)|\xi_1|)}{|\xi_1|}dt_1.$$

We denote
$$\|S_1\|_0^2 :=\frac{1}{(2\pi)^{d+1}}\int_{\bR^d}\int_{\bR}|\cF \phi_{\xi_1}(\tau_1)|^2 \nu(d\tau_1)\mu(d\xi_1).$$
By \eqref{UB-Fourier-G}, for any $\xi_1 \in \bR^d$,
\begin{equation}
\label{bound-convolution-G}
\left|\big(G(t-s,\cdot)*\varphi e^{-i \xi_1 \cdot}\big)(x)\right|=\left|\int_{\bR^d} e^{-i \xi \cdot x} \frac{\sin((t-s)|\xi|)}{|\xi|}\cF^{-1}\varphi(\xi-\xi_1)d\xi \right|\leq  t\,\|\cF^{-1}\varphi\|_1=:C_{t,\varphi},
\end{equation}
where $\|\cdot\|_1$ denotes the norm in $L^1(\bR^d)$. Hence,
$$|\cF \phi_{\xi_1}(\tau_1)|\leq C_{t,\varphi} \left|\int_0^s e^{-i \tau_1 t_1} \frac{\sin((s-t_1)|\xi_1|)}{|\xi_1|}dt_1 \right|,$$
and
\begin{align*}
\|S_1\|_0^2 & \leq C_{t,\varphi}^2 \frac{1}{(2\pi)^{d+1}} \int_{\bR^d}\int_{\bR} \left|\int_0^s e^{-i \tau_1 t_1} \frac{\sin((s-t_1)|\xi_1|)}{|\xi_1|}dt_1 \right|^2 \nu(d\tau_1)\mu(d\xi_1)\\
&= C_{t,\varphi}^2 \, \|f_1(\cdot,s,y)\|_{\cH}^2<\infty \quad \mbox{for any $y \in \bR^d$},
\end{align*}
where the last equality is due to the fact that $f_1(t_1,\cdot,s,y)=1_{(0,s)}(t_1)G(s-t_1,y-\cdot)=:g_{s,y}(t_1,\cdot)$ and in the proof of Theorem \ref{gtx-in-H}, we showed that $\|g_{s,y}\|_{\cH}=\|g_{s,y}\|_{0}=\alpha_1(s)$.
This proves that $S_1 \in \cH$ and
$\|S_1\|_{\cH}^2=\|S_{1}\|_{0}^2 \leq C_{t,\varphi}^2  \|f_1(\cdot,s,y)\|_{\cH}^2$.

Next, we treat the case $n\geq 2$.
By definition, for any $\psi \in \cS(\bR^{nd})$,
$$\big(S_n(t_1,\cdot,\ldots,t_n,\cdot),\psi\big)=\big(G(t-s,\cdot)*\varphi h_{\psi}\big)(x),$$
where $h_{\psi}(y)=\big(f_n(t_1,\cdot,\ldots,t_n,\cdot,s,y),\psi \big)$. First, we show that the Fourier transform of $S_n(t_1,\cdot,\ldots,t_n,\cdot)$ in $\cS'(\bR^{nd})$ is a function and we identify this function.  By \eqref{def-G*phi}, for any $\psi \in \cS(\bR^{nd})$, we have:
\begin{align*}
\big(\cF S_n(t_1,\cdot,\ldots,t_n,\cdot), \psi\big)&= \big( S_n(t_1,\cdot,\ldots,t_n,\cdot), \cF \psi\big)=\big(G(t-s,\cdot)*\varphi h_{\cF \psi} \big)(x)\\
&=\int_{\bR^d} e^{-i \xi \cdot x} \, \frac{\sin((t-s)|\xi|)}{|\xi|} \, \cF^{-1}(\varphi h_{\cF \psi})\,d\xi\\
&=\int_{\bR^d} e^{-i \xi \cdot x} \, \frac{\sin((t-s)|\xi|)}{|\xi|} \, \big( \cF^{-1}\varphi * \cF^{-1} h_{\cF \psi} \big)(\xi) \,d\xi.
\end{align*}

By \eqref{h-psi-Fourier}, $\cF^{-1}h_{\cF \psi}(\eta_n)=1_{\{0<t_1<\ldots<t_n<s\}} \overline{\cF G(s-t_n,\cdot)(\eta_n)}V_{\psi}(\eta_n)$. Therefore,
\begin{align*}
& \big(\cF S_n(t_1,\cdot,\ldots,t_n), \psi\big)=1_{\{0<t_1<\ldots<t_n<s\}} \int_{\bR^d}e^{-i \xi \cdot x}
\,\frac{\sin((t-s)|\xi|)}{|\xi|} \\
& \quad \quad \quad \quad \quad \quad \quad \quad \quad \left(\int_{\bR^d}\cF^{-1}\varphi(\xi-\eta_n) \frac{\sin((s-t_n)|\eta_n|)}{|\eta_n|}V_{\psi}(\eta_n)d\eta_n\right)d\xi \\
&=1_{\{0<t_1<\ldots<t_n<s\}}  \int_{\bR^d}V_{\psi}(\eta_n) \frac{\sin((s-t_n)|\eta_n|)}{|\eta_n|}
\left(\int_{\bR^d}e^{-i\xi \cdot x}\frac{\sin((t-s)|\xi|)}{|\xi|}\cF^{-1} \varphi(\xi-\eta_n) d\xi\right)d\eta_n,
\end{align*}
where the second equality is due to Fubini's theorem. Using definition \eqref{def-V-psi} of $V_{\psi}(\eta_n)$, followed by the change of variables $\xi_1=\eta_1, \xi_k=\eta_k-\eta_{k-1}$ for $k=2,\ldots,n$, we obtain:
\begin{align*}
& \big(\cF S_n(t_1,\cdot,\ldots,t_n,\cdot),\psi\big)=1_{\{0<t_1<\ldots<t_n<s\}}
\int_{\bR^{nd}}\psi(\xi_1,\ldots,\xi_n)
\frac{\sin((t_2-t_1)|\xi_1|)}{|\xi_1|} \ldots \\
& \frac{\sin((s-t_n)|\xi_1+\ldots+\xi_n|)}{|\xi_1+\ldots+\xi_n|}
\left(\int_{\bR^d}e^{-i\xi \cdot x}\frac{\sin((t-s)|\xi|)}{|\xi|}
\cF^{-1} \varphi(\xi-\xi_1-\ldots-\xi_n)
d\xi\right)d\xi_1 \ldots d\xi_n.
\end{align*}
By \eqref{Fourier-phi-a} and \eqref{def-G*phi}, (a version of) the Fourier transform of $S_n(t_1,\cdot,\ldots,t_n,\cdot)$ is the function:
\begin{align}
\nonumber
&\cF S_n(t_1,\cdot,\ldots,t_n,\cdot)(\xi_1, \ldots,\xi_n)=1_{\{0<t_1<\ldots<t_n<s\}}
\frac{\sin((t_2-t_1)|\xi_1|)}{|\xi_1|} \ldots \\
\label{Fourier-Sn}
&  \quad \quad \quad \frac{\sin((s-t_n)|\xi_1+\ldots+\xi_n|)}{|\xi_1+\ldots+\xi_n|}
\big(G(t-s,\cdot)*\varphi e^{-i (\xi_1+\ldots+\xi_n) \cdot}\big)(x).
\end{align}

To prove that $S_n \in \cH^{\otimes n}$, we apply Theorem \ref{phi-in-Hn}.b). The function $(t_1,\ldots,t_n)\mapsto \phi_{\xi_1,\ldots,\xi_n}(t_1,\ldots,t_n):=\cF S_n(t_1,\cdot,\ldots,t_n,\cdot)(\xi_1, \ldots,\xi_n)$ is integrable on $\bR^n$ and its Fourier transform is
\begin{align*}
& \cF \phi_{\xi_1,\ldots,\xi_n}(\tau_1,\ldots,\tau_n)=\int_{\bR^d}e^{-i(\tau_1 t_1+\ldots+\tau_n t_n)}\phi_{\xi_1,\ldots,\xi_n}(t_1,\ldots,t_n)dt_1 \ldots dt_n\\
&=\big(G(t-s,\cdot)*\varphi e^{-i (\xi_1+\ldots+\xi_n) \cdot}\big)(x)
\int_{\{0<t_1<\ldots<t_n<s\}}e^{-i(\tau_1 t_1+\ldots+\tau_n t_n)}
\frac{\sin((t_2-t_1)|\xi_1|)}{|\xi_1|} \ldots \\ & \quad \quad \quad
\quad \quad \quad \quad \quad \quad \quad \quad \quad \quad \quad \quad
\quad \quad \quad \quad\frac{\sin((s-t_n)|\xi_1+\ldots+\xi_n|)}{|\xi_1+\ldots+\xi_n|}dt_1\ldots dt_n.
\end{align*}
We consider the quantity $\|\cdot\|_{0,n}$ defined in Theorem \ref{phi-in-Hn}. By \eqref{bound-convolution-G} and the calculation of $\|f_n(\cdot,s,y)\|_{\cH^{\otimes n}}^2$ given in the proof of Theorem \ref{fn-in-Hn}, we see that
$\|S_n\|_{0,n}^2 \leq C_{t,\varphi}^2 \|f_n(\cdot,s,y)\|_{\cH^{\otimes n}}^2\linebreak <\infty$ for any $y \in \bR^d$.
This proves that $S_n \in \cH^{\otimes n}$ and
$\|S_n\|_{\cH^{\otimes n}}^2 \leq C_{t,\varphi}^2 \, \|f_n(\cdot,s,y)\|_{\cH^{\otimes n}}^2$.

To show that condition (b) holds when $k_n=\widetilde{f}_n$, we denote by $\widetilde{S}_n(t_1,\cdot\ldots,t_n,\cdot)$ the distribution in $\cS'(\bR^{nd})$ given by: for any $\psi \in \cS(\bR^{nd})$,
$$\big(\widetilde{S}_n(t_1,\cdot,\ldots,t_n,\cdot),\psi\big)=
\Big(G(t-s,x-\cdot),\varphi \big(\widetilde{f}_n(t_1,\cdot,\ldots,t_n,\cdot,s,*),\psi\big)\Big).$$
Note that $\widetilde{S}_n$ is the symmetrization of $S_n$, i.e. for any $\psi \in \cS(\bR^d)$,
$$\big(\widetilde{S}_n(t_1,\cdot,\ldots,t_n,\cdot,t,x),\psi\big)=\frac{1}{n!}\sum_{\rho \in \cP_n}\big(S_n(t_{\rho(1)},\cdot,\ldots,t_{\rho(n)},\cdot,t,x),\psi_{\rho} \big),$$
where $\psi_{\rho}$ is defined by \eqref{def-psi-rho}.
Similarly to the case $k_n=f_n$, it can be proved that $\widetilde{S}_n \in \cH^{\otimes n}$ and for any $y \in \bR^d$,
\begin{equation}
\label{bound-norm-Sn}
\|\widetilde{S}_n\|_{\cH^{\otimes n}}^2 \leq C_{t,\varphi}^2 \, \|\widetilde{f}_n(\cdot,s,y)\|_{\cH^{\otimes n}}^2.
\end{equation}

By \eqref{bound-norm-Sn} and Theorem \ref{theorem-sum-finite},
$\sum_{n \geq 1}n! \, \|\widetilde{S}_n\|_{\cH^{\otimes n}}^2 \leq C_{t,\varphi}^2 \sum_{n \geq 1}n! \, \|\widetilde{f}_n(\cdot,s,y)\|_{\cH^{\otimes n}}^2<\infty$, and therefore, the series $\sum_{n\geq 1} I_n(S_n)$ converges in $L^2(\Omega)$.

\vspace{3mm}

{\em Step 3. (Verification of condition (c))}. In this step, we denote $S_{n}^{(t,x,s,\varphi)}$ by $S_n^{(s,\varphi)}$, since $t$ and $x$ are fixed.
By the linearity of the multiple integrals $I_n$, the map $\varphi \mapsto \big(v^{(t,x)}(s,\cdot),\varphi\big)$ is linear from $\cS(\bR^d)$ to $L^2(\Omega)$. This map is $L^2(\Omega)$-continuous, since if $\varphi_k \to \varphi$ in $\cS(\bR^d)$,
\begin{align*}
E\left|\big(v^{(t,x)}(s,\cdot),\varphi_k-\varphi\big)\right|^2&= E\left|\sum_{n\geq 0}I_n(S_n^{(s,\varphi_k-\varphi)}) \right|^2=\sum_{n\geq 0}n!\,\|\widetilde{S}_n^{(s,\varphi_k-\varphi)} \|_{\cH^{\otimes n}}^2\\
&\leq C_{t,\varphi_k-\varphi}^2 \sum_{n\geq 1}\|\widetilde{f}_n(\cdot,s,y)\|_{\cH^{\otimes n}}^2 \to 0 \quad \mbox{as} \ k\to \infty,
\end{align*}
where the inequality is due to \eqref{bound-norm-Sn}, and we recall that $C_{t,\varphi}=t\,\|\cF^{-1}\varphi\|_1$. By Corollary 4.2 of \cite{walsh86}, the process $\{\big(v^{(t,x)}(s,\cdot),\varphi\big);\varphi \in \cS(\bR^d)\}$ has a modification with values in $\cS'(\bR^d)$, which we denote also by $v^{(t,x)}(s,\cdot)$.
We prove that this modification satisfies hypotheses {\em (i)}-{\em (iv)} of Proposition \ref{corr-prop-25-2} (Appendix A).

Hypotheses {\em (ii)}-{\em (iv)} are easily verified. {\em (ii)} holds since by Remark \ref{recurrence-relation-fn} and definition \eqref{def-Sn} of $S_n^{(s,\varphi)}$, for any $s\in (0,t)$, there is a distribution in $\cS'(\bR^{(n+1)d})$, namely $f_{n+1}(t_1,\cdot,\ldots,t_n,\cdot,s,\cdot,t,x)$, which satisfies
$$\big(f_{n+1}(t_1,\cdot,\ldots,t_n,\cdot,s,\cdot,t,x),\psi \otimes \varphi \big)=\big( S_{n}^{(s,\varphi)}(t_1,\cdot,\ldots,t_n,\cdot),\psi\big).$$
{\em (iii)} is clear since $\cF f_{n+1}(t_1,\cdot,\ldots,t_n,\cdot,s,\cdot,t,x)$ is measurable in $(t_1,\ldots,t_{n+1},\xi_1,\ldots,\xi_{n+1})$, bounded, continuous a.e. in $(t_1,\ldots,t_{n+1})$ and $\|f_{n+1}(\cdot,t,x)\|_{\cH^{\otimes (n+1)}}<\infty$ (see the proof of Theorem \ref{fn-in-Hn}). To check {\em (iv)}, we need to prove that the function
$$h(\xi):= \cF f_{n+1}(t_1,\cdot,\ldots,t_n,\cdot,s,\cdot,t,x)(\xi_1,\ldots,\ldots,\xi_n,\xi)\phi(\xi),\quad \xi \in \bR^d$$
is in $\cS(\bR^d)$, for any $t_1, \ldots,t_n,s \in (0,t)$, $\xi_1, \ldots,\xi_n \in \bR^d$ and $\phi \in \cS(\bR^d)$. This is clear since $h(\xi)=Cg(\xi)\phi(\xi)$, where $g(\xi)=e^{-i \xi \cdot x} \frac{\sin((t-s)|\xi_1+\ldots+\xi_n+\xi|)}{|\xi_1+\ldots+\xi_n+\xi|}$ and $C$ is a constant depending on $\xi_1,\ldots,\xi_n$. (Since $g$ is a $C^{\infty}$ function with bounded derivatives, $g \phi \in \cS(\bR^d)$.)

It remains to check that the process $v^{(t,x)}$ has a modification which satisfies hypothesis {\em (i)} of Proposition \ref{corr-prop-25-2}. For this modification, we will show that $\cF v^{(t,x)}(\omega,s,\cdot)$ is a function and we will identify this function. We know that the Fourier transform of $v^{(t,x)}(\omega,s,\cdot)$ is a distribution in $\cS'(\bR^d)$ which satisfies: for any $\varphi \in \cS(\bR^d)$,
\begin{equation}
\label{Fourier-vtx}
\big(\cF v^{(t,x)}(s,\cdot),\varphi\big)=\big(v^{(t,x)}(s,\cdot),\cF \varphi \big)=\sum_{n\geq 0}I_n(S_n^{(s,\cF \varphi)}) \quad \mbox{in} \ L^2(\Omega).
\end{equation}

For any $t_1,\ldots,t_n$, $s \in (0,t)$ and $\xi \in \bR^d$, we consider the distribution $F_n^{(s,\xi)}(t_1,\cdot,\ldots,t_n,\cdot)$ in $\cS'(\bR^{nd})$ whose Fourier transform is the function
\begin{align}
\nonumber
& \cF F_n^{(s,\xi)}(t_1, \cdot,\ldots,t_n,\cdot)(\xi_1,\ldots,\xi_n)= 1_{\{0<t_1<\ldots<t_n<s\}} \, e^{-i (\xi_1+\ldots+\xi_n+\xi) \cdot x}  \, \frac{\sin((t_2-t_1)|\xi_1|)}{|\xi_1|}\ldots \\
\nonumber
& \quad \quad \quad \quad \quad \quad \frac{\sin((s-t_n)|\xi_1+\ldots+\xi_n|)}{|\xi_1+\ldots+\xi_n|}\cdot
\frac{\sin((t-s)|\xi_1+\ldots+\xi_n+\xi|)}{|\xi_1+\ldots+\xi_n+\xi|}\\
\label{Fourier-F-f}
&= \cF f_n(t_1,\cdot,\ldots,t_n,\cdot,s,x)(\xi_1,\ldots,\xi_n)e^{-i \xi \cdot x} \, \frac{\sin((t-s)|\xi_1+\ldots+\xi_n+\xi|)}{|\xi_1+\ldots+\xi_n+\xi|}.
\end{align}

Similarly to Lemma \ref{fn-distribution}, it can be shown that $F_n^{(s,\xi)}(t_1,\cdot,\ldots,t_n,\cdot)$ is a well-defined distribution in $\cS'(\bR^{nd})$.
Using inequality $|\sin x|\leq |x|$ for the last term in \eqref{Fourier-F-f} and Theorem \ref{phi-in-Hn}.b), we see that $F_{n}^{(s,\xi)} \in \cH^{\otimes n}$ and $\|F_n^{(s,\xi)}\|_{\cH^{\otimes n}}^2 \leq t^2 \|f_n(\cdot,s,x)\|_{\cH^{\otimes n}}^2$. The same argument can be used to show that $\|\widetilde{F}_n^{(s,\xi)}\|_{\cH^{\otimes n}}^2 \leq t^2 \|\widetilde{f}_n(\cdot,s,x)\|_{\cH^{\otimes n}}^2$, where
$\widetilde{F}_n^{(s,\xi)}$ is the symmetrization of $F_n^{(s,\xi)}$ defined similarly to \eqref{def-fn-tilde}.

Recall that $I_n(F_n^{(s,\xi)})$ is a square-integrable random variable which is identified with any other random variable that is equal to it in $L^2(\Omega)$.
It is possible to chose the random variable $I_n(F_n^{(s,\xi)})$ such that $(\omega,s,\xi) \mapsto I_n(F_n^{(s,\xi)})$ is measurable.

We claim that
\begin{equation}
\label{relation-2-Fourier}
\cF S_n^{(s,\cF \varphi)}(t_1,\ldots,t_n)(\xi_1,\ldots,\xi_n)=\int_{\bR^d}\cF F_n^{(s,\xi)}(t_1,\cdot,\ldots,t_n,\cdot)(\xi_1,\ldots,\xi_n)\varphi(\xi)d\xi.
\end{equation}
To see this, note that by \eqref{Fourier-Sn},
\begin{align*}
\cF S_n^{(s,\cF \varphi)}(t_1,\cdot,\ldots,t_n,\cdot)(\xi_1, \ldots,\xi_n)&=1_{\{0<t_1<\ldots<t_n<s\}}
\frac{\sin((t_2-t_1)|\xi_1|)}{|\xi_1|} \ldots \\
\frac{\sin((s-t_n)|\xi_1+\ldots+\xi_n|)}{|\xi_1+\ldots+\xi_n|}
& \big(G(t-s,\cdot)*(\cF \varphi )e^{-i (\xi_1+\ldots+\xi_n) \cdot}\big)(x),
\end{align*}
and by \eqref{def-G*phi},
\begin{align*}
\big(G(t-s,\cdot)* (\cF \varphi) e^{-i (\xi_1+\ldots+\xi_n) \cdot} \big)(x)&=
\int_{\bR^d}e^{-i \xi \cdot x} \,\frac{\sin((t-s)|\xi|)}{|\xi|}\cF^{-1}\big((\cF \varphi) e^{-i (\xi_1+\ldots+\xi_n) \cdot} \big)(\xi)d\xi\\
&= \int_{\bR^d}e^{-i \xi \cdot x} \, \frac{\sin((t-s)|\xi|)}{|\xi|}\varphi(\xi-\xi_1 -\ldots-\xi_n)d\xi
\end{align*}
since $\cF^{-1}(\phi \, e^{-i a \cdot})(\xi)=(\cF^{-1} \phi)(\xi-a)$ for any $\phi \in \cS(\bR^d)$ and $a \in \bR^d$. This proves \eqref{relation-2-Fourier}.

For any $\psi \in \cS(\bR^{nd})$, define
\begin{align*}
&\big(T_n^{(s,\varphi)}(t_1,\cdot,\ldots,t_n,\cdot),\psi\big)=\int_{\bR^d} \big( F_n^{(s,\xi)}(t_1,\cdot,\ldots,t_n,\cdot),\psi\big)\varphi(\xi)d\xi\\
& =\int_{\bR^{(n+1)d}}\cF F_n^{(s,\xi)}(t_1,\cdot,\ldots,t_n,\cdot)(\xi_1,\ldots,\xi_n)
\cF^{-1}\psi(\xi_1,\ldots,\xi_n)\varphi(\xi)d\xi_1 \ldots d\xi_n d\xi.
\end{align*}
Similarly to Lemma \ref{fn-distribution}, it can be shown that $T_n^{(s,\varphi)}(t_1,\cdot,\ldots,t_n,\cdot)$ is a well-defined distribution in $\cS'(\bR^{nd})$. Its Fourier transform is a distribution in $\cS'(\bR^{nd})$ given by:
\begin{align*}
&\big(\cF T_n^{(s,\varphi)}(t_1,\cdot,\ldots,t_n,\cdot),\psi\big)=\big( T_n^{(s,\varphi)}(t_1,\cdot,\ldots,t_n,\cdot),\cF \psi\big)\\
&=\int_{\bR^{nd}}\left(\int_{\bR^d}\cF F_n^{(s,\xi)}(t_1,\cdot,\ldots,t_n,\cdot)(\xi_1,\ldots,\xi_n)\varphi(\xi)d\xi\right)
\psi(\xi_1,\ldots,\xi_n)d\xi_1 \ldots d\xi_n\\
&= \int_{\bR^{nd}} \cF S_n^{(s,\cF \varphi)}(t_1,\cdot,\ldots,t_n,\cdot)(\xi_1,\ldots,\xi_n) \psi(\xi_1,\ldots,\xi_n)d\xi_1 \ldots d\xi_,
\end{align*}
for any $\psi \in \cS(\bR^{nd})$, where we used \eqref{relation-2-Fourier} for the last equality. This shows that the Fourier transform of $T_n^{(s,\varphi)}(t_1,\ldots,t_n)$ is a function, namely
$\cF T_n^{(s,\varphi)}(t_1,\ldots,t_n)=\cF S_n^{(s,\cF \varphi)}(t_1,\ldots,t_n)$.
By Step 2 above, $S_{n}^{(s,\cF \varphi)} \in \cH^{\otimes n}$. Hence, $T_n^{(s,\varphi)} \in \cH^{\otimes n}$ and $\|T_{n}^{(s,\varphi)}-S_n^{(s,\cF \varphi)}\|_{\cH^{\otimes n}}^2=0$, which implies that $I_n(T_n^{(s,\varphi)})=I_n(S_n^{(s,\cF \varphi)})$ in $L^2(\Omega)$. We claim that:
\begin{equation}
\label{Fourier-vtx-2}
I_n(S_n^{(s,\cF \varphi)})=I_n(T_n^{(s,\varphi)})=\int_{\bR^d} I_n(F_n^{(s,\xi)})\varphi(\xi)d\xi \quad \mbox{in} \quad L^2(\Omega).
\end{equation}
The second equality follows since $E|I_n(T_n^{(s,\varphi)})|^2=\|T_n^{(s,\varphi)}\|_{\cH^{\otimes n}}^2=\|S_n^{(s,\cF \varphi)}\|_{\cH^{\otimes n}}^2$, by \eqref{relation-2-Fourier},
\begin{align*}
& E\left|\int_{\bR^d}I_n(F_n^{(s,\xi)})\varphi(\xi)d\xi \right|^2=\int_{\bR^{2d}}E[
I_n(F_n^{(s,\xi)}) I_n(F_n^{(s,\xi')})]\varphi(\xi)\varphi(\xi')d\xi d\xi'\\
&= \int_{(\bR_{+}^{2}\times \bR^d)^n} \prod_{i=1}^{n}\gamma(t_i-s_i) \left(\int_{\bR^d}\cF F_n^{(s,\xi)}(t_1,\ldots,t_n)(\xi_1,\ldots,\xi_n)\varphi(\xi)d\xi\right)\\
& \quad  \left(\int_{\bR^d}\cF F_n^{(s,\xi')}(t_1,\ldots,t_n)(\xi_1,\ldots,\xi_n)\varphi(\xi')d\xi'\right)\mu(d\xi_1)\ldots \mu(d\xi_n)d{\bf t}d{\bf s}=\|S_n^{(s,\cF \varphi)}\|_{\cH^{\otimes n}}^2,
\end{align*}
where we used \eqref{relation-2-Fourier} for the second equality. Similarly,
$E\left[ I_n(T_n^{(s,\varphi)})\int_{\bR^d}I_n(F_n^{(s,\xi)})\varphi(\xi)d\xi\right]=\|S_n^{(s,\cF \varphi)}\|_{\cH^{\otimes n}}^2$. By \eqref{Fourier-vtx} and \eqref{Fourier-vtx-2},
\begin{equation}
\label{Fourier-vtx-3}
\big(\cF v^{t,x}(s,\cdot),\varphi\big)=\sum_{n\geq 0}\int_{\bR^d} I_n(F_n^{(s,\xi)})\varphi(\xi)d\xi \quad \mbox{in} \quad L^2(\Omega).
\end{equation}

We claim that:
\begin{equation}
\label{Fourier-vtx-4}
\sum_{n\geq 0}\int_{\bR^d} I_n(F_n^{(s,\xi)})\varphi(\xi)d\xi=\int_{\bR^d}\sum_{n\geq 0} I_n(F_n^{(s,\xi)})\varphi(\xi)d\xi  \quad \mbox{in} \quad L^2(\Omega).
\end{equation}

To see this, we denote $X=\int_{\bR^d}\sum_{n\geq 0} I_n(F_n^{(s,\xi)})\varphi(\xi)d\xi$ and $X_n=\int_{\bR^d} I_n(F_n^{(s,\xi)})\varphi(\xi)d\xi$. We first show that $X$ is a random variable in $L^2(\Omega)$. Note that $\sum_{n\geq 0}I_n(F_n^{(s,\xi)})$ converges in $L^2(\Omega)$ since $\sum_{n\geq 0}n!\,\|\widetilde{F}_n^{(s,\xi)}\|_{\cH^{\otimes n}}^2\leq t^2 \sum_{n\geq 0}n!\,\|\widetilde{f}_n(\cdot,s,x)\|_{\cH^{\otimes n}}^2<\infty$ by Theorem \ref{theorem-sum-finite}. By Minkowski's inequality, $\|X\|_2 \leq \int_{\bR^d}\|\sum_{n\geq 0}I_n(F_n^{(s,\xi)})\|_2 |\varphi(\xi)|d\xi \leq t^2 C_s \int_{\bR^d}|\varphi(\xi)|d\xi<\infty$, where $C_s=\sum_{n\geq 0}n!\, \|\widetilde{f}_n(\cdot,s,x)\|_{\cH^{\otimes n}}^{2}$.
To show that $X=\sum_{n\geq 0}X_n$ in $L^2(\Omega)$, it suffices to prove that
$E(GX)=E(G X_n)$ for any $G=I_n(g)$ with $g \in \cH^{\otimes n}$ symmetric. This is true since $E(GX_n)=\int_{\bR^d}E[I_n(g) I_n(F_n^{(s,\xi)})]\varphi(\xi)d\xi=\int_{\bR^d}n! \, \langle g,\widetilde{F}_n\rangle_{\cH^{\otimes n}}\varphi(\xi)d\xi$ and
\begin{align*}
E(GX)&=E\left(I_n(g) \int_{\bR^d}\sum_{m\geq 0}I_m(F_m^{(s,\xi)}) \varphi(\xi)d\xi\right)=\int_{\bR^d}E\left(\sum_{m\geq 0}I_n(g)I_m(F_{m}^{(s,\xi)})\right)\varphi(\xi)d\xi\\
&=\int_{\bR^d}n! \, \langle g,\widetilde{F}_n\rangle_{\cH^{\otimes n}}\varphi(\xi)d\xi.
\end{align*}

From \eqref{Fourier-vtx-3} and \eqref{Fourier-vtx-4}, we obtain that for any $\varphi \in \cS(\bR^d)$,
\begin{equation}
\label{Fourier-vtx-5}
\big(\cF v^{t,x}(s,\cdot),\varphi\big)=\int_{\bR^d}\sum_{n\geq 0} I_n(F_n^{(s,\xi)})\varphi(\xi)d\xi  \quad \mbox{in} \quad L^2(\Omega).
\end{equation}
For any $s \in (0,t)$ and $\xi \in \bR^d$, consider the random variable $U(s,\xi):=\sum_{n\geq 0} I_n(F_n^{(s,\xi)})$, and define
$$\big(\overline{v}^{(t,x)}(s,\cdot),\cF \varphi \big)=\int_{\bR^d}U(s,\xi)\varphi(\xi)d\xi, \quad \mbox{for all} \ \varphi \in \cS(\bR^d).$$
Due to \eqref{Fourier-vtx-5}, for any $\varphi \in \cS(\bR^d)$, $\big(\overline{v}^{(t,x)}(s,\cdot),\cF \varphi \big)=\big(v^{(t,x)}(s,\cdot),\cF \varphi \big)$ in $L^2(\Omega)$ (hence a.s.). This proves that $\{\big(\overline{v}^{(t,x)}(s,\cdot),\varphi);\varphi \in \cS(\bR^d)\}$ is a modification of
$\{\big(v^{(t,x)}(s,\cdot),\varphi\big);\varphi \in \cS(\bR^d)\}$ which
we will denote also by $v^{(t,x)}$. For any $\omega \in \Omega$ and $s\in (0,t)$, the Fourier transform of $v^{(t,x)}(\omega,s,\cdot)$ is the function $U(\omega,s,\cdot)$. We denote
$$\phi_{\xi}^{(t,x)}(s):=\cF v^{(t,x)}(s,\cdot)(\xi)=\sum_{n\geq 0}I_n(F_n^{(s,\xi)}).$$

The map $s \mapsto \cF v^{(t,x)}(\omega,s,\cdot)(\xi)$ is square-integrable (hence integrable) on $(0,t)$, for almost all $(\omega,\xi) \in \Omega \times \bR^d$. To see this, let $\varphi\geq 0$ be an arbitrary function in $\cS(\bR^d)$. Then
\begin{align*}
&E \int_{\bR^d}  \varphi(\xi) \left(\int_0^t |\cF v^{(t,x)}(s,\cdot)(\xi)|^2 ds\right) d\xi =\int_0^t \int_{\bR^d}\varphi(\xi)\sum_{n\geq 0}n!\, \|\widetilde{F}_n^{(s,\xi)}\|_{\cH^{\otimes n}}^{2}d\xi ds\\
& \leq \int_0^t (t-s)^2 \int_{\bR^d}\varphi(\xi)\sum_{n\geq 0}n! \, \|\widetilde{f}_n(\cdot,s,x)\|_{\cH^{\otimes n}}^2 d\xi ds\\
&\leq C_t \int_0^t (t-s)^2 \int_{\bR^d}\varphi(\xi)d\xi ds<\infty,
\end{align*}
where for the first inequality we used \eqref{Fourier-F-f} and the fact that $\frac{\sin^2((t-s)|\xi_1+\ldots+\xi_n+\xi|)}{|\xi_1+\ldots+\xi_n+\xi|^2} \leq (t-s)^2$ for all $\xi_1,\ldots,\xi_n,\xi \in \bR^d$, and for the second inequality we used
$\sum_{n\geq 0}n! \, \|\widetilde{f}_n(\cdot,s,x)\|_{\cH^{\otimes n}}^2 \leq C_s$ where $C_s>0$ is a non-decreasing function of $s$ (see the proof of Theorem \ref{theorem-sum-finite}).

For any $\tau \in \bR$, we let $\cF \phi_{\xi}^{(t,x)}(\tau)=\int_{0}^t e^{-i \tau s} \cF v^{(t,x)}(s,\cdot)(\xi) ds$. Similarly to \eqref{Fourier-vtx-4}, it can be proved that for almost all $(\tau,\xi)\in \bR^{d+1}$,
$$\cF \phi_{\xi}^{(t,x)}(\tau)=\sum_{n\geq 0}I_n(G_n^{(\tau,\xi)}) \quad \mbox{in} \ L^2(\Omega),$$
where $G_n^{(\tau,\xi)} \in \cH^{\otimes n}$ is such that for any $t_1,\ldots,t_n \in (0,s)$, $G_n^{(\tau,\xi)}(t_1,\cdot,\ldots,t_n,\cdot)$ is a distribution in $\cS'(\bR^{nd})$ whose Fourier transform is the function:
\begin{equation}
\label{Fourier-Gn}\cF G_n^{(\tau,\xi)}(t_1,\cdot,\ldots,t_n,\cdot)(\xi_1,\ldots,\xi_n)=\int_0^t e^{-i \tau s} \cF F_n^{(s,\xi)}(t_1,\cdot,\ldots,t_n,\cdot)(\xi_1,\ldots,\xi_n)ds.
\end{equation}
To see that condition \eqref{u-cond} holds for the process $v^{(t,x)}$, we note that
\begin{equation}
\label{v-hyp-i}
I:=E\left[\int_{\bR^d}\int_{\bR}|\cF \phi_{\xi}^{(t,x)}(\tau)|^2  \nu(d\tau)\mu(d\xi)\right]= \int_{\bR^d}\int_{\bR} \sum_{n\geq 0}n! \, \|\widetilde{G}_n^{(\tau,\xi)}\|_{\cH^{\otimes n}}^2 \nu(d\tau)\mu(d\xi),
\end{equation}
where $\widetilde{G}_n^{(\tau,\xi)}$ is the symmetrization of $G_n^{(\tau,\xi)}$. We need to compute $\|\widetilde{G}_n^{(\tau,\xi)}\|_{\cH^{\otimes n}}^2$. For this, we define $\phi_{\xi_1,\ldots,\xi_n}^{(\tau,\xi)}(t_1,\ldots,t_n):=\cF \widetilde{G}_n^{(\tau,\xi)}(t_1,\cdot,\ldots,t_n,\cdot)(\xi_1,\ldots,\xi_n)$. By \eqref{Fourier-F-f} and \eqref{Fourier-Gn},
\begin{align*}
& \phi_{\xi_1,\ldots,\xi_n}^{(\tau,\xi)}(t_1,\ldots,t_n)=\int_0^t e^{-i \tau s} \cF \widetilde{F}_n^{(s,\xi)}(t_1,\cdot,\ldots,t_n,\cdot)(\xi_1,\ldots,\xi_n)ds \\
&= \int_0^t e^{-i \tau s}\cF \widetilde{f}_n(t_1,\cdot,\ldots,t_n,\cdot,s,x)(\xi_1,\ldots,\xi_n)e^{-i \xi \cdot x} \, \frac{\sin((t-s)|\xi_1+\ldots+\xi_n+\xi|)}{|\xi_1+\ldots+\xi_n+\xi|}ds.
\end{align*}
For any $t_1,\ldots, t_n,s>0$, let $f_{n+1}^{*}(t_1,\cdot,\ldots,t_n,\cdot,s,\cdot,t,x)$ be the distribution in $\cS'(\bR^{(n+1)d})$ whose Fourier transform is the function
\begin{align*}
\cF f_{n+1}^{*}(t_1,\cdot,\ldots,t_n,\cdot,s,\cdot,t,x)(\xi_1, \ldots,\xi_n,\xi)=& \cF \widetilde{f}_n(t_1,\cdot,\ldots,t_n,\cdot,s,x)(\xi_1,\ldots,\xi_n) e^{-i \xi \cdot x} \\
& \frac{\sin((t-s)|\xi_1+\ldots+\xi_n+\xi|)}{|\xi_1+\ldots+\xi_n+\xi|}\,  1_{\{s<t\}}.
\end{align*}
Note that $f_{n+1}^*(\cdot,t,x)$ is similar to $\widetilde{f}_{n+1}(\cdot,t,x)$.
Then
$$\phi_{\xi_1,\ldots,\xi_n}^{(\tau,\xi)}(t_1,\ldots,t_n)=\int_0^t e^{-i \tau s}\phi_{\xi_1,\ldots,\xi_n,\xi}^{(n+1),*}(t_1,\ldots,t_n,s) ds,$$
where $\phi_{\xi_1,\ldots,\xi_n,\xi}^{(n+1),*}(t_1,\ldots,t_n,s)=\cF f_{n+1}^{*}(t_1,\cdot,\ldots,t_n,\cdot,s,\cdot,t,x)(\xi_1, \ldots,\xi_n,\xi)$. It follows that the Fourier transform of the function $(t_1,\ldots,t_n) \mapsto \phi_{\xi_1,\ldots,\xi_n}^{(\tau,\xi)}(t_1,\ldots,t_n)$ is
\begin{align*}
&\cF \phi_{\xi_1,\ldots,\xi_n}^{(\tau,\xi)}(\tau_1,\ldots,\tau_n)=\int_{(0,s)^n}e^{-i (\tau_1 t_1+\ldots+\tau_n t_n)} \phi_{\xi_1,\ldots,\xi_n}^{(\tau,\xi)}(t_1,\ldots,t_n)d{\bf t}\\
&=\int_0^t e^{-i \tau s} \int_{(0,s)^n}e^{-i (\tau_1 t_1+\ldots+\tau_n t_n)} \phi_{\xi_1,\ldots,\xi_n,\xi}^{(n+1),*}(t_1,\ldots,t_n)d{\bf t}ds=\cF \phi_{\xi_1,\ldots,\xi_n,\xi}^{(n+1),*}(\tau_1,\ldots,\tau_n).
\end{align*}
Coming back to \eqref{v-hyp-i}, we see that
\begin{align*}
I&=\sum_{n\geq 0}n!\,\int_{\bR^{(n+1)d}}\int_{\bR^{n+1}}|\cF \phi_{\xi_1,\ldots,\xi_n}^{(\tau,\xi)}(\tau_1,\ldots,\tau_n)|^2 \nu(d\tau_1)\ldots \nu(d\tau_n)\nu(d\tau)\mu(d\xi_1)\ldots \mu(d\xi_n)\mu(d\xi)\\
&=\sum_{n\geq 0}n!\,\int_{\bR^{(n+1)d}}\int_{\bR^{n+1}}|\cF \phi_{\xi_1,\ldots,\xi_n,\xi}^{(n+1),*}(\tau_1,\ldots,\tau_n,\tau)|^2 \nu(d\tau_1)\ldots \nu(d\tau_n)\nu(d\tau)\mu(d\xi_1)\ldots \mu(d\xi_n)\mu(d\xi)\\
&=\sum_{n\geq 0}n!\, \|f_{n+1}^*(\cdot,t,x)\|_{\cH^{\otimes (n+1)}}^2<
\infty.
\end{align*}
(The last series converges by the same argument as in the proof of Theorem \ref{theorem-sum-finite}.) This shows that the process $v^{(t,x)}$ satisfies hypothesis {\em (i)} of Proposition \ref{corr-prop-25-2}.

\vspace{3mm}

{\em Step 4. (Verification of condition (d)).} We apply Proposition Proposition \ref{corr-prop-25-2} to the process $v^{(t,x)}$. Recall that hypothesis {\em (ii)} of this proposition holds for the distribution $f_{n+1}(t_1,\cdot,\ldots,t_n,\cdot,s,\cdot,t,x)$. The symmetrization of $f_{n+1}(\cdot,t,x)$ in all $n+1$ variables is $\widetilde{f}_{n+1}(\cdot,t,x)$.
As in the last part of the proof of Theorem \ref{d2-exist-sol}, we conclude that
$v^{(t,x)} \in {\rm Dom}\ \delta$ and $u(t,x)=1+\delta(v^{(t,x)})$ in $L^2(\Omega)$. $\Box$

\section{Uniqueness of the solution}
\label{section-uniqueness}

In this section, we establish the uniqueness of the solution. We discuss separately the cases $d\leq2$ and $d\geq 3$.

In the case $d\leq 2$, the proof of the uniqueness of the solution is the same as for the Parabolic Anderson Model (see Section 4.1 of \cite{hu-nualart09}).
We include this proof for the sake of completeness.

\begin{theorem}
If $d\leq 2$ and $\mu$ satisfies \eqref{Dalang-cond}, then the unique solution (up to a modification) of equation \eqref{wave} is the process $u$ given by \eqref{def-u}.
\end{theorem}

\noindent {\bf Proof:} Let $u$ be a solution of equation \eqref{wave} with Wiener chaos expansion \eqref{def-u-kernel-k} for some symmetric non-negative functions $k_n(\cdot,t,x) \in \cH^{\otimes n}$. We fix $t>0$ and $x \in \bR^d$. Let $k_0(t,x)=1$. We will show that $k_n(\cdot,t,x)=\widetilde{f}_n(\cdot,t,x)$ for any $n\geq 1$. Let $v^{(t,x)}$ be the process defined by \eqref{def-v}.

For any $s>0$ and $y \in \bR^d$, $v^{(t,x)}(s,y)=\sum_{n \geq 0}I_{n}(g_n^{(t,x)}(\cdot,s,y))$, where
\begin{equation}
\label{def-gn-allvar}
g_{n}^{(t,x)}(t_1,x_1,\ldots,t_n,x_n,s,y)=1_{(0,t)}(s)
G(t-s,x-y)k_n(t_1,x_1,\ldots,t_n,x_n,s,y).
\end{equation}

We apply Proposition \ref{corr-prop-25-0} to the process $v^{(t,x)}$. Hypothesis {\em (i)} of this proposition is verified as in the proof of Theorem \ref{d2-exist-sol}. Hypothesis {\em (ii)} holds i.e. $g_n^{(t,x)} \in \cH^{\otimes (n+1)}$, since
\begin{align*}
& \sum_{n\geq 0}n! \, \|g_n^{(t,x)}\|_{\cH^{\otimes (n+1)}}^2=\int_{((0,t) \times \bR^d)^2} \gamma(s-r)f(y-z) G(t-s,x-y)G(t-r,x-z)\\
& \quad \quad \quad \quad \quad \quad \quad \quad \quad \quad \quad \quad \sum_{n\geq 0}n!\, \langle k_n(\cdot,s,y),k_n(\cdot,s,z)\rangle_{\cH^{\otimes n}}dsdydrdz\\
&\quad =\int_{((0,t) \times \bR^d)^2} \gamma(s-r)f(y-z) G(t-s,x-y)G(t-r,x-z) E[u(s,y)u(r,z)]dsdydrdz\\
& \quad \leq C_t \int_{((0,t) \times \bR^d)^2} \gamma(s-r)f(y-z) G(t-s,x-y)G(t-r,x-z) dsdydrdz<\infty,
\end{align*}
using Cauchy-Schwarz inequality and \eqref{moment2-u}.
(Since $g_{n}^{(t,x)}$ is non-negative, it follows that $g_n^{(t,x)} \in |\cH^{\otimes (n+1)}| \subset \cH^{\otimes (n+1)}$, where the space $|\cH^{\otimes n}|$ is defined similarly to $|\cH|$; see \eqref{inclusion-H}.)

Since $u$ is a solution, $v^{(t,x)} \in {\rm Dom}\ \delta$ and $u(t,x)=1+\delta(v^{t,x})$. On the other hand, by Proposition \ref{corr-prop-25-0},
$\delta(v^{(t,x)})=\sum_{n \geq 0}I_{n+1}(\widetilde{g_{n}^{(t,x)}})$,
where $\widetilde{g_0^{(t,x)}}=g_0^{(t,x)}$ and for $n\geq 1$, $\widetilde{g_n^{(t,x)}}$ is the symmetrization of $g_n^{(t,x)}$ in
all $n+1$ variables, defined by:
\begin{align}
\label{def-gn-tilde}
& \widetilde{g_n^{(t,x)}}(t_1,x_1,\ldots,t_n,x_n,s,y) =
\frac{1}{n+1}\Big[
g_n^{(t,x)}(t_1,x_1,\ldots,t_n,x_n,s,y)+  \\
\nonumber
& \quad \quad \quad \sum_{i=1}^{n}
g_n^{(t,x)}(t_1,x_1,\ldots,t_{i-1},x_{i-1},s,y,t_{i+1},x_{i+1},\ldots,t_n,x_n,t_i,x_i)
\Big].
\end{align}

This shows that
$\sum_{n \geq 0}I_{n+1}(k_{n+1}(\cdot,t,x))=u(t,x)-1=\sum_{n \geq 0}I_{n+1}(\widetilde{g_n^{(t,x)}})$ in $L^2(\Omega)$.
By the uniqueness of the Wiener chaos expansion with symmetric
kernels (see e.g. Theorem 1.1.2 of \cite{nualart06}), it follows that for any $n \geq 0$,
$k_{n+1}(\cdot,t,x)=\widetilde{g_n^{(t,x)}}$, i.e.
\begin{equation}
\label{k-g-relation}
k_{n+1}(t_1,x_1,\ldots,t_n,x_n,t_{n+1},x_{n+1},t,x)=
\widetilde{g_n^{(t,x)}}(t_1,x_1,\ldots,t_n,x_n,t_{n+1},x_{n+1}).
\end{equation}

The functions $k_n(\cdot,t,x)$ can now be found recursively. By \eqref{k-g-relation} with $n=0$, we obtain:
$$k_1(t_1,x_1,t,x)=\widetilde{g_0^{(t,x)}}(t_1,x_1)=g_0^{(t,x)}(t_1,x_1)=1_{(0,t)}(t_1) G(t-t_1,x-x_1).$$
Using \eqref{k-g-relation} with $n=1$, followed by the definition of $\widetilde{g_1^{(t,x)}}$ and relation \eqref{def-gn-allvar}, we obtain:
\begin{eqnarray*}
\lefteqn{k_2(t_1,x_1,t_2,x_2,t,x)=\widetilde{g_1^{(t,x)}}(t_1,x_1,t_2,x_2)
= \frac{1}{2}\left[g_1^{(t,x)}(t_1,x_1,t_2,x_2)+g_1^{(t,x)}(t_2,x_2,t_1,x_1)\right]= } \\
& & \frac{1}{2}\Big[1_{(0,t)}(t_2)
G(t-t_2,x-x_2)k_1(t_1,x_1,t_2,x_2)+1_{(0,t)}(t_1)
G(t-t_1,x-x_1) k_1(t_2,x_2,t_1,x_1) \Big].
\end{eqnarray*}
We now use the formula for $k_1$ that we found above. This leads us to conclude that
\begin{eqnarray*}
k_2(t_1,x_1,t_2,x_2,t,x)&=& \frac{1}{2}\Big[1_{(0,t)}(t_2)
G(t-t_2,x-x_2) 1_{(0,t_2)}(t_1)G(t_2-t_1,x_2-x_1)+\\
& & 1_{(0,t)}(t_1) G(t-t_1,x-x_1)
1_{(0,t_1)}(t_2)G(t_1-t_2,x_1-x_2) \Big],
\end{eqnarray*}
that is $k_2(\cdot,t,x)=\widetilde{f}_2(\cdot,t,x)$. Iterating this procedure, we infer that $k_n(\cdot,t,x)=\widetilde{f}_n(\cdot,t,x)$ for any $n \geq 1$. $\Box$

\vspace{3mm}

The next result gives the uniqueness of the solution in the case $d\geq 3$. For this, we use Lemma \ref{def-kn-lemma}.

\begin{theorem}
If $d\geq 3$, $\mu$ satisfies \eqref{Dalang-cond} and Hypothesis A holds, then the unique solution (up to a modification) of equation \eqref{wave} is the process $u$ given by \eqref{def-u}.
\end{theorem}

\noindent {\bf Proof:} Let $u$ be a solution of equation \eqref{wave} with Wiener chaos expansion \eqref{def-u-kernel-k} for some elements $k_n(\cdot,t,x) \in \cH^{\otimes n}$ as in Definition \ref{new-def-solution}. We fix $t>0$ and $x \in \bR^d$. Let $k_0(t,x)=1$. We will show that $k_n(\cdot,t,x)=\widetilde{f}_n(\cdot,t,x)$ for any $n\geq 1$.

By definition (see Remark \ref{def-gn-tx}), for any $s>0$ and $\varphi \in \cS(\bR^d)$,
$\big(v^{(t,x)}(s,\cdot),\varphi\big)=\sum_{n\geq 0}I_n(S_n^{(t,x,s,\varphi)})$, where $S_0^{(t,x,s,\varphi)}=1_{(0,t)}(s)\big(G(t-s,x-\cdot),\varphi\big)$ and $S_n^{(t,x,s,\varphi)}(t_1,\cdot,\ldots,t_n,\cdot)$ is a symmetric distribution in $\cS'(\bR^{nd})$ which satisfies, for any $\psi \in \cS(\bR^{nd})$,
\begin{equation}
\label{def-Sn-2}
\big(S_n^{(t,x,s,\varphi)}(t_1,\cdot,\ldots,t_n,\cdot),\psi\big)=1_{(0,t)}(s)
\Big(G(t-s,x-\cdot),\varphi \big(k_n(t_1,\cdot,\ldots,t_n,\cdot,s,*),\psi \big)\Big),
\end{equation}
where $*$ denotes the missing argument of the function $h_{\psi}(y)=\big(k_n(t_1,\cdot,\ldots,t_n,\cdot,s,y),\psi \big)$.

We apply Proposition \ref{corr-prop-25-2} to the process $v^{(t,x)}$. Hypotheses {\em (i)}-{\em (iv)} of this proposition are verified due to condition (c) in Definition \ref{new-def-solution}. In particular, by hypothesis {\em (ii)}, we know that for any $t_1, \ldots,t_n,s \in (0,t)$, there exists a distribution in $\cS'(\bR^{(n+1)d})$, which we denote by $g_{n}^{(t,x)}(t_1,\cdot,\ldots,t_n,\cdot,s,\cdot)$, such that for any $\psi \in \cS(\bR^{nd})$ and $\varphi \in \cS(\bR^d)$,
$\big(g_{n}^{(t,x)}(t_1,\cdot,\ldots,t_n,\cdot,s,\cdot), \psi \otimes \varphi\big)=
\big(S_n^{(t,x,s,\varphi)}(t_1,\cdot,\ldots,t_n,\cdot),\psi\big)$. Due to \eqref{def-Sn-2}, this means that for any $\psi \in \cS(\bR^d)$ and $\varphi \in \cS(\bR^d)$,
\begin{align}
\label{def-gn-allvar3}
\big(g_{n}^{(t,x)}(t_1,\cdot,\ldots,t_n,\cdot,s,\cdot), \psi \otimes \varphi\big)=1_{(0,t)}(s) \Big(G(t-s,\cdot)* \varphi \big(k_n(t_1,\cdot,\ldots,t_n,\cdot,s,*),\psi \big)\Big)(x).
\end{align}

Since $u$ is a solution, $v^{(t,x)} \in {\rm Dom}\ \delta$ and
$u(t,x)=1+\delta(v^{(t,x)})$.
On the other hand, by Proposition \ref{corr-prop-25-2},
$\delta(v^{(t,x)})=\sum_{n \geq 0}I_{n+1}(\widetilde{g_{n}^{(t,x)}})$,
where $\widetilde{g_0^{(t,x)}}=g_0^{(t,x)}$ and for $n\geq 1$, $\widetilde{g_n^{(t,x)}}$ is the symmetrization of $g_n^{(t,x)}$ in
all $n+1$ variables. In particular, for $\phi_1, \ldots,\phi_n,\phi \in \cS(\bR^d)$,
\begin{align}
\label{def-gn-tilde3}
& \Big(\widetilde{g_n^{(t,x)}}(t_1,\cdot,\ldots,t_n,\cdot,s,\cdot),\phi_1 \otimes \ldots \otimes \phi_n \otimes \phi \Big)=\\
\nonumber
&\frac{1}{n+1}\Big[
\Big(g_n^{(t,x)}(t_1,\cdot,\ldots,t_n,\cdot,s,\cdot),\phi_1 \otimes \ldots \otimes \phi_n \otimes \phi \Big)+  \\
\nonumber
& \sum_{i=1}^{n}
\Big(g_n^{(t,x)}(t_1,\cdot,\ldots,t_{i-1},\cdot,s,\cdot,t_{i+1},\cdot,\ldots,t_n,\cdot,t_i,\cdot),
\phi_1 \otimes \ldots \otimes \phi_{i-1}\otimes \phi \otimes \phi_{i+1}\otimes \ldots \phi_n \otimes \phi_i \Big)
\Big].
\end{align}

This shows that
$\sum_{n \geq 0}I_{n+1}(k_{n+1}(\cdot,t,x))=u(t,x)-1=\sum_{n \geq 0}I_{n+1}(\widetilde{g_n^{(t,x)}})$ in $L^2(\Omega)$.
By the uniqueness of the Wiener chaos expansion, for any $n \geq 0$
$k_{n+1}(\cdot,t,x)=\widetilde{g_n^{(t,x)}}$, i.e.
\begin{equation}
\label{k-g-relation3}
k_{n+1}(t_1,\cdot,\ldots,t_n,\cdot,t_{n+1},\cdot,t,x)=
\widetilde{g_n^{(t,x)}}(t_1,\cdot,\ldots,t_n,\cdot,t_{n+1},\cdot).
\end{equation}

The elements $k_n(\cdot,t,x)$ can now be found recursively. By \eqref{k-g-relation3} with $n=0$, we obtain:
$$k_1(t_1,\cdot,t,x)=\widetilde{g_0^{(t,x)}}(t_1,\cdot)=g_0^{(t,x)}(t_1,\cdot)=1_{[0,t]}(t_1) G(t-t_1,x-\cdot).$$
Using \eqref{k-g-relation3} with $n=1$, followed by the definition of $\widetilde{g_1^{(t,x)}}$ and relation \eqref{def-gn-allvar3}, we obtain:
\begin{align*}
& \big(k_2(t_1,\cdot,t_2,\cdot,t,x),\phi_1 \otimes \phi_2 \big)=\big(\widetilde{g_1^{(t,x)}}(t_1,\cdot,t_2,\cdot),\phi_1 \otimes \phi_2 \big)\\
&  =\frac{1}{2}\left[\big(g_1^{(t,x)}(t_1,\cdot,t_2,\cdot),\phi_1 \otimes \phi_2 \big)
+\big(g_1^{(t,x)}(t_2,\cdot,t_1,\cdot),\phi_2 \otimes \phi_1 \big)\right]\\
&  =\frac{1}{2}\Big[1_{(0,)]}(t_2)
\Big(G(t-t_2,\cdot)*\phi_2 \big(k_1(t_1,\cdot,t_2,*),\phi_1\big)\Big)(x)+\\
&  \quad \quad \quad 1_{(0,t)}(t_1)
\Big(G(t-t_1,\cdot)*\phi_1 \big(k_1(t_2,\cdot,t_1,*),\phi_2\big) \Big)(x)\Big].
\end{align*}
We now use the formula for $k_1$ that we found above.
Note that
$$\big(k_1(t_1,\cdot,t_2,x_2),\phi_1\big)=1_{(0,t_2)}(t_1)\big(G(t_2-t_1,x_2-\cdot),\phi_1 \big)=1_{(0,t_2)}(t_1) \big(G(t_2-t_1,\cdot)*\phi_1\big)(x_2).$$
A similar formula holds for $\big(k_1(t_2,\cdot,t_1,*),\phi_2\big)$.
This leads us to conclude that
\begin{align*}
\big(k_2(t_1,\cdot,t_2,\cdot,t,x),\phi_1 \otimes \phi_2\big)& =
\frac{1}{2}\left[1_{(0,t)}(t_2)1_{(0,t_2)}(t_1)
\Big(G(t-t_2,\cdot)*\phi_2 \big(G(t_2-t_1,\cdot)*\phi_1\big)\Big)(x)+ \right.\\
&  \left.1_{(0,t)}(t_1)1_{(0,t_1)}(t_2)
\Big(G(t-t_1,\cdot)*\phi_1 \big(G(t_1-t_2,\cdot)*\phi_2\big) \Big)(x)\right].
\end{align*}
By Lemma \ref{def-kn-lemma}, the last term above is exactly $\big(\widetilde{f}_2(t_1,\cdot,t_2,\cdot,t,x),\phi_1 \otimes \phi_2\big)$. (Recall  definition \eqref{def-fn-tilde} of $\widetilde{f}_2(\cdot,t,x)$.)
Hence, $k_2(t_1,\cdot,t_2,\cdot,t,x)=\widetilde{f}_2(t_1,\cdot,t_2,\cdot,t,x)$ for any $t_1>0,t_2>0$. Iterating this procedure, we infer that $k_n(\cdot,t,x)=\widetilde{f}_n(\cdot,t,x)$ for any $n \geq 1$. $\Box$

\section{Moments of the solution}
\label{section-moments}

In this section, we show that the solution to equation \eqref{wave} is $L^p(\Omega)$-continuous and has uniformly bounded moments of order $p$, for any $p\geq 2$. Similar results exist for parabolic equations (see for instance, Theorem 3.2 of \cite{HHNT} and Theorem 5.2 of \cite{song15}). Recall that
\begin{equation}
\label{second-moment-u}
E|u(t,x)|^2=\sum_{n \geq 0}\frac{1}{n!}\alpha_n(t),
\end{equation}
where $\alpha_n(t)$ is given by \eqref{alpha-J}.

\begin{theorem}
\label{exist-th}
Suppose that the measure $\mu$ satisfies condition
\eqref{Dalang-cond}. If $d \geq 4$, suppose in addition that $\mu$ satisfies Hypothesis A. Then for any $p\geq 2$, the solution $u$ to equation \eqref{wave} is $L^p(\Omega)$-continuous, and
$$\sup_{(t,x)\in [0,T] \times \bR^d}E|u(t,x)|^p<\infty \quad \mbox{for all} \quad T>0.$$
\end{theorem}

\noindent {\bf Proof:} {\em Step 1.} We show that the $p$-th moments of $u$ are uniformly bounded.

We proceed as in the proof of Theorem 4.2 of \cite{balan12}. We denote by $\|\cdot\|_p$ the $L^p(\Omega)$-norm. We use the fact that for any $F \in \cH_n$ and $p \geq 2$,
\begin{equation}
\label{p-norm-F}
\|F\|_{p} \leq (p-1)^{n/2}\|F\|_2
\end{equation}
(see last line of page 62 of \cite{nualart06}). Using Minkowski's inequality, applying \eqref{p-norm-F} for $F=J_n(t,x)$, and invoking \eqref{alpha-J} and \eqref{alpha-bound}, we see that:
\begin{eqnarray*}
\|u(t,x)\|_p &\leq & \sum_{n \geq 0}\|J_n(t,x)\|_p \leq  \sum_{n \geq 0}(p-1)^{n/2}\|J_n(t,x)\|_2=
 \sum_{n \geq 0}(p-1)^{n/2}\left( \frac{1}{n!}\alpha_n(t)\right)^{1/2} \\
& \leq &  \sum_{n \geq 0}(p-1)^{n/2} \Gamma_t^{n/2}  8^{n/2} \sum_{k=0}^{n}\frac{t^{n/2+k}}{(k!)^{1/2}}D_N^{k/2} C_N^{(n-k)/2}.
\end{eqnarray*}
The last term is uniformly bounded for $(t,x) \in [0,T] \times \bR^d$ (using the same argument as in the proof of Theorem \ref{theorem-sum-finite}).

{\em Step 2.} We show that $u$ is $L^p(\Omega)$-continuous.

The argument in Step 1 above shows that for any $T>0$ and $p \geq 2$,
$$\sum_{n \geq 0}\sup_{(t,x) \in [0,T] \times \bR^d}\|J_n(t,x)\|_p \leq C_{T,p}<\infty.$$
Hence $\{u_n(t,x)=\sum_{k=0}^nJ_k(t,x)\}_{n \geq 1}$ converges to $u(t,x)$ in $L^p(\Omega)$, uniformly in $(t,x) \in [0,T] \times \bR^d$. By Lemma \ref{Jn-cont-lemma} below, $J_n$ is $L^p(\Omega)$-continuous, and hence $u_n$
is $L^p(\Omega)$-continuous. Therefore, $u$  is $L^p(\Omega)$-continuous. $\Box$

\vspace{3mm}

The following result is an extension of Lemma 4.2 of \cite{balan12} to the case of an arbitrary covariance function $\gamma(t)$.

\begin{lemma}
Under the conditions of the Theorem \ref{exist-th}, we have:\\
\label{Jn-cont-lemma}
a) for any $p \geq 2$, $n \geq 1$ and $t>0$,
$$E|J_n(t+h,x)-J_n(t,x)|^p \to 0 \quad as \ h\to 0, \ uniformly \ in \ x \in \bR^d;$$
b) for any $p \geq2$, $n \geq 1$, $t>0$ and $x \in \bR^d$
$$E|J_n(t,x+z)-J_n(t,x)|^p \to 0 \quad as \ |z|\to 0,z \in \bR^d.$$
\end{lemma}

\noindent {\bf Proof:} a) We assume that $|h| \leq 1$ and $h>0$. (The case $h<0$ is similar.) By \eqref{p-norm-F},
\begin{eqnarray}
\nonumber
\|J_n(t+h,x)-J_n(t,x)\|_p^2 &\leq &(p-1)^n E|J_n(t+h,x)-J_n(t,x)|^2 \\
\nonumber
&=& (p-1)^n n! \, \|\widetilde{f}_n(\cdot,t+h,x)-\widetilde{f}_n(\cdot,t,x)\|_{\cH^{\otimes n}}^{2} \\
\label{bound-A-B}
&\leq & (p-1)^n \frac{2}{n!} \left(A_n(t,h)+B_n(t,h)\right),
\end{eqnarray}
where
\begin{eqnarray}
\label{def-A}
A_n(t,h)&=& (n!)^2 \|\widetilde{f}_n(\cdot,t+h,x)1_{[0,t]^{n}}-\widetilde{f}_n(\cdot,t,x) \|_{\cH^{\otimes n}}^{2}\\
\label{def-B}
B_n(t,h)&=& (n!)^2 \|\widetilde{f}_n(\cdot,t+h,x)1_{[0,t+h]^{n}\verb2\2 [0,t]^n} \|_{\cH^{\otimes n}}^{2}.
\end{eqnarray}

We evaluate $A_n(t,h)$ first. We have:
$$A_n(t,h) = \int_{[0,t]^{2n}}
\prod_{j=1}^{n}\gamma(t_j-s_j) \psi_{t,h}^{(n)}({\bf t}, {\bf
s})d{\bf t} d{\bf s},$$
where
\begin{eqnarray*}
\psi_{t,h}^{(n)}({\bf t},{\bf s}) &=&\int_{\bR^{nd}}
\cF\left[g_{\bf t}^{(n)}(\cdot,t+h,x)-g_{\bf
t}^{(n)}(\cdot,t,x)\right](\xi_1,\ldots,\xi_n) \\
& & \overline{\cF\left[g_{\bf
s}^{(n)}(\cdot,t+h,x)-g_{\bf s}^{(n)}(\cdot,t,x)\right](\xi_1, \ldots,\xi_n)}\mu(d\xi_1)
\ldots \mu(d\xi_n)
\end{eqnarray*}
and $g_{\bf t}^{(n)}(\cdot,t,x)$ is given by \eqref{gtn}. By the Cauchy-Schwarz inequality and the inequality $ab \leq (a^2+b^2)/2$,
$$\psi_{t,h}^{(n)}({\bf t},{\bf s}) \leq \left(\psi_{t,h}^{(n)}({\bf t},{\bf t})\right)^{1/2}
\left(\psi_{t,h}^{(n)}({\bf s},{\bf s})\right)^{1/2} \leq \frac{1}{2}\left(
\psi_{t,h}^{(n)}({\bf t},{\bf t})+\psi_{h}^{(n)}({\bf s},{\bf s})\right).$$
Using the symmetry of the function $\gamma$ and Lemma \ref{basic-ineq-lemma}, it follows that
\begin{equation}
\label{bound-A}
A_n(t,h) \leq \int_{[0,t]^{2n}} \prod_{j=1}^{n}\gamma(t_j-s_j)\psi_{t,h}^{(n)}({\bf t},{\bf t})d{\bf t}d{\bf s} \leq \Gamma_t^n \int_{[0,t]^{n}} \psi_{t,h}^{(n)}({\bf t},{\bf t})d{\bf t}.
\end{equation}

Using definition \eqref{Fourier-gn}) of the Fourier transform of $g_{\bf t}^{(n)}(\cdot,t,x)$, we see that
\begin{eqnarray}
\nonumber
\psi_{t,h}^{(n)}({\bf t},{\bf t})
&=& \int_{\bR^{nd}} |\cF G(u_1, \cdot \ )(\xi_{\rho(1)})|^2 \ldots
|\cF G(u_{n-1}, \cdot \ )(\xi_{\rho(1)}+ \ldots +\xi_{\rho(n-1)})|^2
\\
\nonumber
& & |\cF [G(u_n+h, \cdot \ )-G(u_n , \cdot \ )](\xi_{\rho(1)}+
\ldots +\xi_{\rho(n)})|^2 \mu(d\xi_1) \ldots \mu(d\xi_n)\\\
\nonumber
&=& \int_{\bR^{nd}} |\cF G(u_1, \cdot \ )(\xi_1')|^2 \ldots
|\cF G(u_{n-1}, \cdot \ )(\xi_1'+ \ldots +\xi_{n-1}')|^2
\\
\label{estimate-psi-n}
& & |\cF [G(u_n+h, \cdot \ )-G(u_n , \cdot \ )](\xi_{1}'+
\ldots +\xi_{n}')|^2 \mu(d\xi_1') \ldots \mu(d\xi_n'),
\end{eqnarray}
where $u_j=t_{\rho(j+1)}-t_{\rho(j)}$, $\xi_j'=\xi_{\rho(j)}$ and $0<t_{\rho(1)}< \ldots<t_{\rho(n)}<t=t_{\rho(n+1)}$. By the continuity of the function $t \mapsto \cF G(t,\cdot)(\xi)$ and the dominated convergence theorem, $\psi_{t,h}^{(n)}({\bf t},{\bf t})\to 0$ as $h \to 0$. To justify the application of this theorem, note that by \eqref{UB-Fourier-G}, $|\cF G(u_j,\cdot)(\xi_1'+\ldots+\xi_j')|^2 \leq C_t \frac{1}{1+|\xi_1'+\ldots+\xi_j'|^2}$ for $j=1,\ldots, n-1$,
$$|\cF [G(u_n+h, \cdot \ )-G(u_n , \cdot \ )](\xi_{1}'+
\ldots +\xi_{n}')|^2 \leq  4C_t\frac{1}{1+|\xi_1'+\ldots+\xi_n'|^2},$$
and by Lemma \ref{lemmaA} and condition \eqref{Dalang-cond},
$$\int_{\bR^{d}} \frac{1}{1+|\xi_1'|^2}\ldots\left(\int_{\bR^d}\frac{1}{1+|\xi_1'+\ldots+\xi_n'|^2}\mu(d\xi_n')
\right)\ldots\mu(d\xi_1')\leq \left(\int_{\bR^d}\frac{1}{1+|\xi|^2}\mu(d\xi)\right)^n<\infty.$$

By applying the dominated convergence theorem again, we infer that
\begin{equation}
\label{integral-psi-th}
\int_{[0,t]^n}\psi_{t,h}^{(n)}({\bf t},{\bf t})d{\bf t} \to 0, \quad \mbox{as} \ h\to 0,
\end{equation}
and hence, $A_n(t,h) \to 0$ as $h \to 0$, due to \eqref{bound-A}. To justify the application of this theorem, we use some estimates borrowed from the proof of Theorem \ref{theorem-sum-finite}. First, note that from \eqref{estimate-psi-n}, we infer that:
\begin{eqnarray}
\label{bound-psi-th}
\psi_{t,h}^{(n)}({\bf t},{\bf t})&\leq &
\prod_{j=1}^{n-1}\left(\sup_{\eta \in \bR^d} \int_{\bR^d}|\cF G(u_j,\cdot)(\xi_j+\eta)|^2 \mu(d\xi_j)\right)\\
\nonumber
& & \cdot \sup_{\eta \in \bR^d}\int_{\bR^d} |\cF G(u_n+h,\cdot)(\xi_n+\eta)-\cF G(u_n,\cdot)(\xi_n+\eta)|^2 \mu(d\xi_n).
\end{eqnarray}
Using Lemma \ref{sup-lemma} and relation \eqref{UB-Fourier-G}, it follows that
\begin{align*}
\psi_{t,h}^{(n)}({\bf t},{\bf t})&\leq \prod_{j=1}^{n-1}\left(\int_{\bR^d} \frac{4u_j^2}{1+u_j^2|\xi_j|^2}\mu(d\xi_j)\right) 4C_t \sup_{\eta \in \bR^d} \int_{\bR^d}\frac{1}{1+|\xi_n+\eta|^2}\mu(d\xi_n) \\
&=  \prod_{j=1}^{n-1}\left(\int_{\bR^d} \frac{4u_j^2}{1+u_j^2|\xi_j|^2}\mu(d\xi_j)\right)4C_t C, \quad \mbox{with} \ C=\int_{\bR^d}\frac{1}{1+|\xi|^2}\mu(d\xi).
\end{align*}
The $d{\bf t}$ integral of the last term on $[0,t]^n$ is equal to
\begin{eqnarray*}
& & 4C_t C \,n!\int_0^t \left(\int_{0<t_{1}<\ldots<t_{n-1}<t_n} \prod_{j=1}^{n-1}\left(\int_{\bR^d} \frac{4(t_{j+1}-t_{j})^2}{1+(t_{j+1}-t_j)^2|\xi_j|^2}\mu(d\xi_j) \right)
dt_1 \ldots dt_{n-1}\right) dt_n\\
& & =4C_t C \,n! \int_0^t I^{(n-1)}(t_n)dt_n, \quad \mbox{where $I^{(n)}(t)$ is defined by \eqref{def-In(t)}}.
\end{eqnarray*}
To see that the last integral is finite, we recall that $I^{(n-1)}(t_n) \leq J^{(n-1)}(t_n) \leq K^{(n-1)}(t_n)$, where $J^{(n)}(t)$ and $K^{(n)}(t)$ are defined by \eqref{def-Jn(t)}, respectively \eqref{def-Kn(t)}. This shows that the application of the dominated convergence theorem is justified to prove \eqref{integral-psi-th}.

As for the term $B_n(t,h)$, note that
$$B_n(t,h) =\int_{[0,t+h]^{2n}}
\prod_{j=1}^{n}\gamma(t_j-s_j)
\gamma_{t,h}^{(n)}({\bf t}, {\bf s})1_{D_{t,h}}({\bf t}) 1_{D_{t,h}}({\bf s})d{\bf t} d{\bf s},
$$
where $D_{t,h}=[0,t+h]^n \verb2\2 [0,t]^n$ and
$$\gamma_{t,h}^{(n)}({\bf t},{\bf s}) = \int_{\bR^{nd}} \cF g_{\bf
t}^{(n)}(\cdot \ , t+h,x)(\xi_1, \ldots,\xi_n)\overline{\cF g_{\bf s}^{(n)}(\cdot \
, t+h,x)(\xi_1, \ldots,\xi_n)} \mu(d\xi_1) \ldots \mu(d\xi_n).$$
By the Cauchy-Schwarz inequality and the inequality $ab \leq (a^2+b^2)/2$,
$$\gamma_{t,h}^{(n)}({\bf t},{\bf s}) \leq \left(\gamma_{t,h}^{(n)}({\bf t},{\bf t})\right)^{1/2}
\left(\gamma_{t,h}^{(n)}({\bf s},{\bf s})\right)^{1/2} \leq \frac{1}{2}\left(
\gamma_{t,h}^{(n)}({\bf t},{\bf t})+\gamma_{h}^{(n)}({\bf s},{\bf s})\right).$$
Using the symmetry of the function $\gamma$ and Lemma \ref{basic-ineq-lemma}, it follows that:
\begin{eqnarray}
\nonumber
B_n(t,h) &\leq &\int_{[0,t+h]^{2n}} \prod_{j=1}^{n}\gamma(t_j-s_j)\gamma_{t,h}^{(n)}({\bf t},{\bf t})1_{D_{t,h}}({\bf t})1_{D_{t,h}}({\bf s}) d{\bf t}d{\bf s} \\
\nonumber
& \leq & \int_{[0,t+h]^{2n}} \prod_{j=1}^{n}\gamma(t_j-s_j)\gamma_{t,h}^{(n)}({\bf t},{\bf t})1_{D_{t,h}}({\bf t}) d{\bf t}d{\bf s} \\
\label{bound-B}
& \leq & \Gamma_{t+h}^n \int_{[0,t+h]^{n}} \gamma_{t,h}^{(n)}({\bf t},{\bf t}) 1_{D_{t,h}}({\bf t}) d{\bf t}.
\end{eqnarray}

We observe that for any ${\bf t}=(t_1, \ldots,t_n)\in [0,t+h]^n$, if we denote $u_j=t_{\rho(j+1)}-t_{\rho(j)}$ for $j=1,\ldots,n-1$ and $u_n=t-t_{\rho(n)}$, where $\rho \in {\cal P}_n$ is such that $0<t_{\rho(1)}<\ldots<t_{\rho(n)}<t+h$, then
\begin{eqnarray}
\nonumber
\gamma_{t,h}^{(n)}({\bf t},{\bf t})&=&\int_{\bR^{nd}} |\cF G(u_1, \cdot)(\xi_{\rho(1)})|^2 \ldots
|\cF G(u_{n-1}, \cdot)(\xi_{\rho(1)}+ \ldots +\xi_{\rho(n-1)})|^2
\\
\nonumber
& & |\cF G(u_n+h, \cdot)(\xi_{\rho(1)}+
\ldots +\xi_{\rho(n)})|^2 \mu(d\xi_1) \ldots \mu(d\xi_n)\\
\nonumber
&\leq & \prod_{j=1}^{n-1}\left(\sup_{\eta \in \bR^d} \int_{\bR^d}|\cF G(u_j,\cdot)(\xi_j+\eta)|^2 \mu(d\xi_j)\right)\\
\label{bound-gamma-th}
& & \left(\sup_{\eta \in \bR^d} \int_{\bR^d}|\cF G(u_n+h,\cdot)(\xi_n+\eta)|^2 \mu(d\xi_n)\right)
\end{eqnarray}
which is bounded by a constant of the form $C_t^n$ for any $h \in [0,1]$, due to \eqref{sup-G-2}. The fact that $B_n(t,h) \to 0$ as $h \to 0$ follows from \eqref{bound-B} by the dominated convergence theorem, since $D_{t,h} \to \emptyset$ as $h \to 0$.

b) By \eqref{p-norm-F}, we have:
\begin{equation}
\label{bound-C2}
\|J_n(t,x+z)-J_n(x,z)\|_p^2 \leq (p-1)^n E|J_n(t,x+z)-J_n(t,x)|^2=(p-1)^n \frac{1}{n!}C_n(t,z),
\end{equation} where
\begin{eqnarray}
\nonumber
C_n(t,z)&=&(n!)^2 \|\widetilde{f}_n(\cdot,t,x+z)-\widetilde{f}_n(\cdot,t,x)\|_{\cH^{\otimes n}}^2 \\
\label{def-C}
&=& \int_{[0,t]^{2n}} \prod_{j=1}^{n}\gamma(t_j-s_j) \psi_{t,z}^{(n)}({\bf t},{\bf s})d{\bf t}d{\bf s}
\end{eqnarray}
and
\begin{eqnarray*}
\psi_{t,z}^{(n)}({\bf t},{\bf s})&=&   \int_{\bR^d} \cF \left[g_{\bf t}^{(n)}(\cdot,t,x+z)-g_{\bf t}^{(n)}(\cdot,t,x) \right](\xi_1,\ldots,\xi_n)\\
& & \, \, \, \, \, \, \overline{\cF \left[g_{\bf s}^{(n)}(\cdot,t,x+z)-g_{\bf s}^{(n)}(\cdot,t,x) \right](\xi_1,\ldots,\xi_n)} \mu(d\xi_1) \ldots \mu(d\xi_n).
\end{eqnarray*}
By the Cauchy-Schwarz inequality and the inequality $ab \leq (a^2+b^2)/2$,
$$\psi_{t,z}^{(n)}({\bf t},{\bf s}) \leq \left(\psi_{t,z}^{(n)}({\bf t},{\bf t}) \right)^{1/2} \left(
\psi_{t,z}^{(n)}({\bf s},{\bf s}) \right)^{1/2} \leq \frac{1}{2}\left( \psi_{t,z}^{(n)}({\bf t},{\bf t}) +\psi_{t,z}^{(n)}({\bf s},{\bf s}) \right).$$
Using the symmetry of $\gamma$ and Lemma \ref{basic-ineq-lemma}, it follows that
\begin{equation}
\label{bound-C}
C_n(t,z)\leq \int_{[0,t]^{2n}}\prod_{j=1}^{n}\gamma(t_j-s_j)\psi_{t,z}^{(n)}({
\bf t},{\bf t})d{\bf t}d{\bf s} \leq \Gamma_t^n \int_{[0,t]^n}\psi_{t,z}^{(n)}({\bf t},{\bf t})d{\bf t}.
\end{equation}
Using the definition \eqref{Fourier-gn} of the Fourier transform of $g_{\bf t}^{(n)}(\cdot,t,x)$, we see that
\begin{eqnarray}
\nonumber
\psi_{t,z}^{(n)}({\bf t},{\bf t})&=& \int_{\bR^{nd}}
 |\cF G(u_1,\cdot)(\xi_{\rho(1)})|^2 \ldots |\cF G(u_{n-1},\cdot)(\xi_{\rho(1)}+\ldots \xi_{\rho(n-1)})|^2\\
\label{bound-psi-z}
 & & |\cF G(u_{n},\cdot)(\xi_{\rho(1)}+\ldots \xi_{\rho(n)})|^2
 |1-e^{-i (\xi_1+\ldots+\xi_n)\cdot z}|^2 \mu(d\xi_1)\ldots \mu(d\xi_n),
\end{eqnarray}
where $u_j=t_{\rho(j+1)}-t_{\rho(j)}$ and $0<t_{\rho(1)}<\ldots<t_{\rho(n)}<t=t_{\rho(n+1)}$. By applying twice the dominated convergence theorem, we conclude first that $\psi_{t,z}^{(n)}({\bf t},{\bf t}) \to 0$ when $|z|\to 0$, and then that $C_n(t,z) \to 0$ when $|z| \to 0$. $\Box$

\section{H\"older continuity}
\label{section-Holder}

In this section, we assume that the spectral measure $\mu$ satisfies \eqref{Holder-cond} and we show that the solution of equation \eqref{wave} has a H\"older continuous modification. {Note that \eqref{Holder-cond} implies \eqref{Dalang-cond}.

We will need the following results.

\begin{proposition}[Proposition 7.4 of \cite{conus-dalang09}]
Let $G$ be the fundamental solution of the wave equation in dimension $d \geq 1$. If $\mu$ satisfies \eqref{Holder-cond}, then: \\
(i) for any $T>0$ and $M>0$, there exists a constant $C>0$ depending on $T,d,M,\beta$ such that for any $h \in \bR$ with $|h|\leq M$
\begin{equation}
\label{H3}
\sup_{t \in [0,T \wedge (T-h)]}\sup_{\eta \in \bR^d}  \int_{\bR^d} |\cF G(t+h,\cdot)(\xi+\eta)-\cF G(t,\cdot)(\xi+\eta)|^2\mu(d\xi) \leq C |h|^{2-2\beta};
\end{equation}

\noindent (ii) for any $T>0$, there exists a constant $C>0$ depending on $T,d,\beta$ such that for any $t \in [0,T]$
\begin{equation}
\label{H4}
\sup_{\eta \in \bR^d}\int_{\bR^d} |\cF G(t,\cdot)(\xi+\eta)|^2 \mu(d\xi) \leq C t^{2-2\beta};
\end{equation}

\noindent (iii) for any $T>0$ and for any compact set $K \subset \bR^d$, there exists a constant $C>0$ depending on $T,K,d,\beta$ such that for any $z \in K$,
\begin{equation}
\label{H5}
\sup_{t \in [0,T]}\sup_{\eta \in \bR^d}\int_{\bR^d} |\cF G(t,\cdot)(\xi+\eta)|^2 |1-e^{-i(\xi+\eta) \cdot z}|^2 \mu(d\xi)\leq C|z|^{2-2\beta}.
\end{equation}
\end{proposition}

\begin{lemma}
\label{lemma-BT}
For any $t>0$ and $h>-1$
$${\cal I}_n(t,h):=\int_{0<t_1<\ldots<t_n<t} \prod_{j=1}^{n-1}(t_{j+1}-t_j)^h (t-t_n)^hd{\bf t}=
\frac{\Gamma(1+h)^{n+1}}{\Gamma(n(1+h)+1)}t^{n(1+h)}.$$
\end{lemma}

We are now ready to state our result about the H\"older continuity of the solution.

\begin{theorem}
\label{holder-th}
Suppose that $\mu$ satisfies \eqref{Holder-cond}. If $d \geq 4$, suppose in addition that $\mu$ satisfies Hypothesis A. Let $u$ be the solution of equation \eqref{wave}.
Then:\\
 a) for any $p \geq 2$ and $T>0$ there exists a constant $C>0$ depending on $p,T,d$ and $\beta$ such that for any $t,t' \in [0,T]$ and for any $x \in \bR^d$,
\begin{equation}
\label{t-increm}
\|u(t,x)-u(t',x)\|_p \leq C |t-t'|^{1-\beta};
\end{equation}
b) for any $p \geq 2$, $T>0$ and compact set $K \subset \bR^d$, there exists a constant $C>0$ depending on $p,T,K,d$ and $\beta$ such that for any $t \in [0,T]$ and for any $x,x' \in K$,
\begin{equation}
\label{x-increm}
\|u(t,x)-u(t,x')\|_p \leq C |x-x'|^{1-\beta}.
\end{equation}

Consequently, for any $T>0$ and for any compact set $K \subset \bR^d$, the solution $\{u(t,x); t \in [0,T], x\in  K\}$ to equation \eqref{wave} has a modification which is jointly $\theta$-H\"older continuous in time and space, for any $\theta\in (0, 1-\beta)$.
\end{theorem}

\begin{remark}
 {\rm
 If $f(x)=|x|^{-\alpha}$ is the Riesz kernel for some $0<\alpha<d$, then the spectral measure $\mu$ is given by $\mu(d\xi)=C_{\alpha,d}|\xi|^{-(d-\alpha)}d\xi$, where $C_{\alpha,d}>0$ is a constant which depends on $\alpha$ and $d$. In this case, condition \eqref{Dalang-cond} holds for any $0<\alpha<2$ and
condition \eqref{Holder-cond} holds for any
 $\beta$ with $\alpha/2<\beta<1$. Therefore, for any $T>0$ and for any compact set $K \subset \bR^d$, the solution $u=\{u(t,x);t \in [0,T],x \in K\}$ has a modification which  is jointly $\theta$-H\"older continuous in time and space, for any $\theta\in(0, \frac{2-\alpha}{2})$. This result coincides with Theorem 5.1 of \cite{balan12}.
}
\end{remark}

\noindent {\bf Proof of Theorem \ref{holder-th}:}
a) Let $t,t' \in [0,T]$ and $x \in \bR^d$ be arbitrary. Assume that $h:=t'-t>0$. (The case $h<0$ is similar.) By Minkowski's inequality, \eqref{p-norm-F} and \eqref{bound-A-B},
\begin{eqnarray}
\nonumber
\|u(t+h,x)-u(t,x)\|_p &\leq & \sum_{n\geq 0}(p-1)^{n/2} \|J_n(t+h,x))-J_n(t,x)\|_2 \\
\label{moment-time}
&\leq &\sum_{n\geq 0} (p-1)^{n/2} \left(\frac2{n!}\left[ A_n(t,h)+B_n(t,h)\right] \right)^{1/2},
\end{eqnarray}
where $A_n(t,h)$ and $B_n(t,h)$ are given by \eqref{def-A}, respectively \eqref{def-B}.

To estimate $A_n(t,h)$, we use \eqref{bound-A}. Note that by \eqref{bound-psi-th}, \eqref{H3} and \eqref{H4}, we have
$$\psi_{t,h}^{(n)}({\bf t},{\bf t}) \leq C^n (u_1 \ldots u_{n-1}h)^{2-2\beta},$$
where $u_j=t_{\rho(j+1)}-t_{\rho(j)}$ and $0<t_{\rho(1)}<\ldots<t_{\rho(n)}<t=t_{\rho(n+1)}$. By invoking Lemma \ref{lemma-BT}, it follows that
\begin{eqnarray*}
A_n(t,h) &\leq & h^{2-2\beta} \Gamma_t^n C^n  n! \int_{0<t_1<\ldots<t_n<t} \prod_{j=1}^{n-1}(t_{j+1}-t_j)^{2-2\beta}dt_1 \ldots dt_n\\
&=& h^{2-2\beta} \Gamma_t^n C^n  n! \int_0^t {\cal I}_{n-1}(t_n,2-2\beta)dt_n
\\
&=&  h^{2-2\beta} \Gamma_t^n C^n n! \frac{\Gamma(3-2\beta)^n}{\Gamma((n-1)(3-2\beta)+1)}\int_0^t t_n^{(n-1)(3-2\beta)}dt_n.
\end{eqnarray*}
We now use the fact that for all $a>1$ there exists a constant $C>0$ such that
\begin{equation}
\label{Gamma-LB}
\Gamma(an+1) \geq C (n!)^a \quad \mbox{for all} \ n\geq 1.
\end{equation}
It follows that
\begin{equation}
\label{bound-A2}
A_n(t,h) \leq h^{2-2\beta} \Gamma_t^n C^n \frac{1}{(n!)^{2-2\beta}} t^{(n-1)(3-2\beta)+1}.
\end{equation}

To estimate $B_n(t,h)$, we use \eqref{bound-B}. First note that by \eqref{bound-gamma-th} and \eqref{H4},
$$\gamma_{t,h}^{(n)}({\bf t},{\bf t}) \leq C^n [u_1 \ldots u_{n-1}(u_n+h)]^{2-2\beta},$$
where $u_j=t_{\rho(j+1)}-t_{\rho(j)}$ and $0<t_{\rho(1)}<\ldots<t_{\rho(n)}<t=t_{\rho(n+1)}$. We observe that if $(t_1, \ldots,t_n) \in D_{t,h}=[0,t+h]^n \verb2\2 [0,t]^n$ then there exists at least one index $i$ with $t_i>t$. So,
$$D_{t,h}=\bigcup_{\rho \in S_n}\{(t_1, \ldots,t_n): 0\leq t_{\rho(1)}\leq \ldots\leq t_{\rho(n-1)} \leq t_{\rho(n)}, t<t_{\rho(n)} \leq t+h\}.$$
By applying Lemma \ref{lemma-BT}, it follows that
\begin{eqnarray*}
B_n(t,h) & \leq & \Gamma_{t+h}^n C^n \sum_{\rho \in S_n} \int_t^{t+h}
\int_{0<t_{\rho(1)}<\ldots<t_{\rho(n-1)}
<t_{\rho(n)}} \prod_{j=1}^{n-1}(t_{\rho(j+1)}-t_{\rho(j)})^{2-2\beta}(t+h-t_{\rho(n)})^{2-2\beta}d{\bf t} \\
&=& \Gamma_{t+h}^n C^n n! \int_t^{t+h} {\cal I}_{n-1}(t_n,2-2\beta) \, (t+h-t_n)^{2-2\beta}dt_n\\
&=& \Gamma_{t+h}^n C^n n! \frac{\Gamma(3-2\beta)^n}{\Gamma((n-1)(3-2\beta)+1)}
\int_t^{t+h}t_n^{(n-1)(3-2\beta)}(t+h-t_n)^{2-2\beta}dt_n \\
&=&  \Gamma_{t+h}^n C^n n! \frac{\Gamma(3-2\beta)^n}{\Gamma((n-1)(3-2\beta)+1)}
\int_0^{h}(t+h-u)^{(n-1)(3-2\beta)}\, u^{2-2\beta}du \\
& \leq & \Gamma_{T}^n C^n n! \frac{\Gamma(3-2\beta)^n}{\Gamma((n-1)(3-2\beta)+1)} T^{(n-1)(3-2\beta)} \frac{1}{3-2\beta} h^{3-2\beta}.
\end{eqnarray*}
Using \eqref{Gamma-LB}, it follows that
\begin{equation}
\label{bound-B2}
B_n(t,h)\leq  h^{2-2\beta} \Gamma_{T}^n C^n \frac{1}{(n!)^{2-2\beta}} T^{(n-1)(3-2\beta)}.
\end{equation}

Relation \eqref{t-increm} follows from \eqref{moment-time}, \eqref{bound-A2} and \eqref{bound-B2}.

b) Let $t \in [0,T]$ and $x,x' \in K$ be arbitrary. We denote $z=x'-x$. By Minkowski's inequality, \eqref{p-norm-F} and \eqref{bound-C2}, we have:
$$\|u(t,x+z)-u(t,x)\|_p \leq \sum_{n \geq 0}(p-1)^{n/2} \|J_n(t,x+z)-J_n(t,x)\|_2=\sum_{n \geq 0}(p-1)^{n/2} \left(\frac{1}{n!}C_{n}(t,z) \right),$$
where $C_n(t,z)$ is defined by \eqref{def-C}. To estimate $C_n(t,z)$ we use \eqref{bound-C}. Note that by \eqref{bound-psi-z}, \eqref{H4} and \eqref{H5},
$$\psi_{t,z}^{(n)}({\bf t},{\bf t}) \leq C^n |z|^{2-2\beta}(u_1 \ldots u_{n-1})^{2-2\beta},$$
where $u_j=t_{\rho(j+1)}-t_{\rho(j)}$ and $0<t_{\rho(1)}<\ldots<t_{\rho(n)}<t=t_{\rho(n+1)}$. Hence
$$
C_n(t,z)  \leq |z|^{2-2\beta}C^n \Gamma_t^n n!
\int_{0<t_{1}<\ldots<t_{n}<t} \prod_{j=1}^{n-1}(t_{j+1}-t_{j})^{2-2\beta}d{\bf t}.$$
Using the same estimate for the last integral as above, we infer that
$$C_n(t,z) \leq |z|^{2-2\beta}C^n \Gamma_t^n \frac{1}{(n!)^{2-2\beta}}t^{(n-1)(3-2\beta)+1}.$$
Relation \eqref{x-increm} follows. The final statement is a consequence of Kolmogorov's continuity theorem.
$\Box$

{\bf Acknowledgement:} The authors would like to thank an anonymous referee for reading the paper very carefully, and for pointing out several typos.

\appendix
\label{appendixA}

\section{Malliavin calculus results}

In this section, we give some results which allow us to compute the Skorohod integral of a process $u$ (which may be a distribution in the space variable) based on its Wiener chaos expansion.
The first result is the simplest one (when the kernels are functions in all variables) and is used for proving the existence and uniqueness of the solution to equation \eqref{wave} in the case $d\leq 2$. The last result is used in the case $d\geq 3$.

\begin{proposition}
\label{corr-prop-25-0}
Let $u=\{u(t,x);t\geq 0,x \in \bR^d\}$ be a process defined on a probability space $(\Omega,\cF,P)$, such that for any $t\geq 0$ and $x\in \bR^d$, $E|u(t,x)|^2<\infty$ and
$$u(t,x)=\sum_{n\geq 0}I_n(f_n(\cdot,t,x)) \quad \mbox{in} \ L^2(\Omega),$$
where $f_0(t,x)=E[u(t,x)]$ and for $n \geq 1$, $f_n(\cdot,t,x)$ is a symmetric function in $\cH^{\otimes n}$.

 Suppose that:\\
(i) $u\in \cH$ a.s., the map $\omega \mapsto \|u(\omega)\|_{\cH}$ is measurable and $E\|u\|_{\cH}^2<\infty$;\\
(ii) $f_n \in \cH^{\otimes (n+1)}$ for any $n\geq 0$.

We denote by $\widetilde{f}_n$ the symmetrization of $f_n$ in all $n+1$ variables.
Then $u \in {\rm Dom} \ \delta$ if and only if $V:=\sum_{n\geq 0}I_{n+1}(\widetilde{f}_n)$ converges in $L^2(\Omega)$. In this case, $\delta(u)=V$.
\end{proposition}

\begin{remark}
\label{corr-prop-25-remark}
{\rm
Due to \eqref{inclusion-H}, condition {\em (i)} above holds if $u$ is jointly measurable and
\begin{equation}
\label{barH-norm}
E\left[\int_{(\bR_{+} \times \bR^d)^2}\gamma(t-s)f(x-y)|u(t,x)u(s,y)|dtdx ds dy\right]<\infty.
\end{equation}
}
\end{remark}

\noindent {\bf Proof of Proposition \ref{corr-prop-25-0}:} We use the same argument as in white noise case (see Proposition 1.3.7 of \cite{nualart06}). We include the details for the sake of completeness.

{\em Step 1.} We prove that for any $G=I_n(g)$ with $g$ a symmetric function in $\cD((\bR_{+} \times \bR^d)^n)$,
\begin{equation}
\label{MalliavinG-u}
E[\langle DG,u \rangle_{\cH}]=E[G I_{n}(\widetilde{f}_{n-1})].
\end{equation}

Since $D_{s,y}G=nI_{n-1}(g(\cdot,s,y))$, 
by the orthogonality of the Wiener chaos spaces,
$$E[u(t,x)D_{s,y}G]=n (n-1)!\,\langle f_{n-1}(\cdot,t,x),g(\cdot,s,y)\rangle_{\cH^{\otimes (n-1)}}.$$
Hence,
\begin{align}
\nonumber
&E[\langle DG,u \rangle_{\cH}]=\int_{(\bR_{+}\times
\bR^d)^2}\gamma(t-s)f(x-y)E[u(t,x)D_{s,y}G]dxdydtds\\
\label{MalliavinG-u-step1}
&= n!\int_{(\bR_{+}\times \bR^d)^2}\gamma(t-s)f(x-y)\,\langle f_{n-1}(\cdot,t,x),g(\cdot,s,y)\rangle_{\cH^{\otimes (n-1)}}dxdydtds\\
\nonumber
&=n! \int_{(\bR_{+} \times \bR^d)^{2n}}\gamma(t-s)f(x-y)
\prod_{i=1}^{n-1}\gamma(t_i-s_i)\prod_{i=1}^{n-1}f(x_i-y_i)\\
\nonumber
& \quad \quad \quad f_{n-1}(t_1,x_1,\ldots,t_{n-1},
x_{n-1},t,x)g(s_1,y_1,\ldots,s_{n-1},y_{n-1},s,y)d{\bf t}d{\bf x}d{\bf s}d{\bf y}\\
\nonumber
&=n!\, \langle g,f_{n-1}\rangle_{\cH^{\otimes n}}=n!\, \langle g,\widetilde{f}_{n-1}\rangle_{\cH^{\otimes n}}=E[GI_{n}(\widetilde{f}_{n-1})].
\end{align}

{\em Step 2.} We prove that relation \eqref{MalliavinG-u} holds also for $G=I_n(g)$ with $g \in \cH^{\otimes n}$ arbitrary.
Since $\cD((\bR_{+} \times \bR^d)^n)$ is dense in $\cH^{\otimes n}$, there exists a sequence $(g_k)_k$ of functions in $\cD((\bR_{+}\times \bR)^n)$ such that $g_k \to g$ in $\cH^{\otimes n}$. Hence $G_k=I_n(g_k) \to G=I_n(g)$ in $L^2(\Omega)$. By Step 1, relation \eqref{MalliavinG-u} holds for $G_k$ for any $k$. Letting $k \to \infty$, we infer that this relation also holds for $G$. On the left-hand side, we use the fact $DG_k\to DG$ in $L^2(\Omega;\cH)$. To see this, note that $D$ is a closable operator from $L^2(\Omega)$ to $L^2(\Omega;\cH)$ and $(DG_k)_k$ converges in $L^2(\Omega;\cH)$
(since $E\|DG_k-DG_l\|_{\cH}^2=n n!\|g_k-g_l\|_{\cH^{\otimes n}}^2 \to 0$ as $k,l\to \infty$).

\vspace{2mm}

{\em Step 3.} Suppose that the series $V$ converges in $L^2(\Omega)$. We show that $u \in {\rm Dom} \ \delta$ and $\delta(u)=V$. Let $F \in \bD^{1,2}$ be arbitrary. Say $F=\sum_{n\geq 0}I_n(g_n)$ with $g_n \in \cH^{\otimes n}$ symmetric. Let $F_N=\sum_{n=0}^{N}I_n(g_n)$. Note that $E\|D F_N-D F_{M}\|_{\cH}^2=\sum_{n=N+1}^{M}n n! \|g_n\|_{\cH^{\otimes n}} \to 0$ as $N,M \to \infty$ since $\sum_{n\geq 1}nE|I_n(g_n)|^2<\infty$ by Proposition 1.2.2 of \cite{nualart06}. Hence $(DF_N)_N$ converges in $L^2(\Omega;\cH)$. Since $D$ is a closable operator from $L^2(\Omega)$ to $L^2(\Omega;\cH)$, $DF_N \to DF \in L^2(\Omega;\cH)$. Let $G_k=I_k(g_k)$ for $k\geq 1$. For any $k=1, \ldots,N$, by \eqref{MalliavinG-u}, we have
$$E[VI_k(g_k)]=E[I_k(\widetilde{f}_{k-1}) I_k(g_k)]=E[\langle DG_k,u \rangle_{\cH}].$$
For $k=0$, $E[V g_0]=g_0E(V)=0$. The sum of these equations for $k=0,\ldots,N$, leads to $E[V F_N]=E[\langle DF_N,u \rangle_{\cH}]$. Letting $N\to \infty$, we obtain:
$$E[V F]=E[\langle DF,u \rangle_{\cH}],$$
using on the right-hand side, the fact that $E\|u\|_{\cH}^2<\infty$, which is hypothesis {\em (i)}. This shows that
$|E[\langle DF,u \rangle_{\cH}]|\leq \|V\|_2 \|F\|_2$ for any $F \in \bD^{1,2}$. Hence, $u \in {\rm Dom}\ \delta$ and $\delta(u)=V$.

\vspace{2mm}

{\em Step 4.} Suppose that $u \in {\rm Dom}\ \delta$. By \eqref{MalliavinG-u}, for any $G=I_n(g)$ with $g \in \cH^{\otimes n}$ symmetric,
$$E[G\delta(u)]=E[\langle DG,u \rangle_{\cH}]=E[G I_n(\widetilde{f}_{n-1})].$$
This shows that $I_n(\widetilde{f}_{n-1})$ is the projection of $\delta(u)$ on $\cH_n$, i.e. $\delta(u)=\sum_{n\geq 1}I_{n}(\widetilde{f}_{n-1})$ and the series converges in $L^2(\Omega)$. $\Box$

\vspace{3mm}

The next result will be used in the proof of Proposition \ref{corr-prop-25-2} below, where it will be applied to a regularization $u_{\varepsilon}$ of the process $u$.

\begin{proposition}
\label{corr-prop-25-1}
Let $u=\{u(t,x);t\geq 0,x \in \bR^d\}$ be a process defined on a probability space $(\Omega,\cF,P)$, such that for any $t\geq 0$ and $x\in \bR^d$, $E|u(t,x)|^2<\infty$,
$$u(t,x)=\sum_{n\geq 0}I_n(f_n(\cdot,t,x)) \quad \mbox{in} \ L^2(\Omega),$$
where $f_0(t,x)=E[u(t,x)]$ and for $n \geq 1$, $f_n(\cdot,t,x)\in \cH^{\otimes n}$ is such that $f_n(t_1,\cdot,\ldots,t_n,\cdot,t,x)$ is a symmetric distribution in $\cS'(\bR^{nd})$ whose Fourier transform is a function which has a version such that $(t_1,\ldots,t_n, \xi_1,\ldots,\xi_n)\mapsto \cF f_n(t_1,\cdot,\ldots,t_n,\cdot,t,x)(\xi_1,\ldots,\xi_n)$ is measurable.

Suppose that $u(t,\cdot)=0$ if $t>T$ and:

(i) $u\in \cH$ a.s., the map $\omega \mapsto \|u(\omega)\|_{\cH}$ is measurable and $E\|u\|_{\cH}^2<\infty$;

(ii) for any $n\geq 1$ and $t_1, \ldots,t_n,t \in [0,T]$, there is a distribution $f_n(t_1,\cdot,\ldots,t_n,\cdot,t,\cdot)$ in $\cS'(\bR^{(n+1)d})$ such that for any
$\psi \in \cS(\bR^{nd})$ and $\varphi \in \cS(\bR^d)$,
$$\big(f_n(t_1,\cdot,\ldots,t_n,\cdot,t,\cdot),\psi \otimes \varphi\big) = \int_{\bR^d}\big(f_n(t_1,\cdot,\ldots,t_n,\cdot,t,x),\psi\big) \varphi(x) dx;$$

(iii) for any $t_1,\ldots,t_n,t \in [0,T]$, the Fourier transform of $f_n(t_1,\cdot,\ldots,t_n,\cdot,t,\cdot)$ is a function which has a version such that $(t_1,\ldots,t_n,t,\xi_1,\ldots,\xi_n,\xi)\mapsto \cF f_n(t_1,\cdot,\ldots,t_n,\cdot,t,\cdot)\linebreak(\xi_1, \ldots,\xi_n,\xi)$ is measurable; for every $\xi_1, \ldots,\xi_n,\xi \in \bR^d$, the function $(t_1,\ldots,t_n,t) \mapsto  \cF f_n(t_1,\cdot,\ldots,t_n,\cdot,t,\cdot)(\xi_1, \ldots,\xi_n,\xi)$ is bounded, continuous a.e. on $[0,T]^{n+1}$ and satisfies
\begin{align*}
\int_{([0,T]^2\times \bR^d)^{n+1}} & \prod_{i=1}^{n}\gamma(t_i-s_i)\gamma(t-s) \cF f_n(t_1,\cdot,\ldots,t_n,\cdot,t,\cdot)(\xi_1,\ldots,\xi_n,\xi)\\
&\overline{\cF f_n(s_1,\cdot,\ldots,s_n,\cdot,s,\cdot)(\xi_1,\ldots,\xi_n,\xi)}\mu(d\xi_1)\ldots \mu(d\xi_n)\mu(d\xi) d{\bf t} d{\bf s}dtds<\infty.
\end{align*}

(iv) for any ${\bf t}=(t_1,\ldots,t_n)\in [0,T]^n$ and $\underline{\xi}=(\xi_1,\ldots,\xi_n) \in \bR^{nd}$, the map $x \mapsto \cF f_n(t_1,\cdot,\ldots,t_n,\cdot,t,x )=:\varphi_{{\bf t},\underline{\xi}}(x)$ is in $\cS(\bR^d)$.

We denote by $\widetilde{f}_n$ the symmetrization of $f_n$ in all $n+1$ variables.
Then $u \in {\rm Dom} \ \delta$ if and only if $V:=\sum_{n\geq 0}I_{n+1}(\widetilde{f}_n)$ converges in $L^2(\Omega)$. In this case, $\delta(u)=V$.
\end{proposition}

\noindent {\bf Proof:} We use the same argument as for Proposition \ref{corr-prop-25-0}. We only need to show the statement in Step 1, since Steps 2, 3 and 4 remain valid without any modification.

Hypothesis {\em (iii)} guarantees that $f_n \in \cH^{\otimes (n+1)}$ for any $n \geq 0$, by Theorem \ref{phi-in-Hn}.c).  Hypothesis {\em (ii)} implies that, for any $\psi \in \cS(\bR^{nd})$ and $\varphi \in \cS(\bR^d)$,
\begin{align*}
&(\cF f_n(t_1,\cdot,\ldots,t_n,\cdot), \psi \otimes \varphi)=\int_{\bR^d}\big(f_n(t_1,\cdot,\ldots,t_n,\cdot,t,x),\cF \psi\big) \cF \varphi(x)dx\\
&\quad \quad \quad =\int_{\bR^d}\int_{\bR^{nd}}\cF f_n(t_1,\cdot,\ldots,t_n,\cdot,t,x)(\xi_1,\ldots,\xi_n)
 \psi(\xi_1,\ldots,\xi_n) \cF \varphi(x) d\xi_1 \ldots d\xi_ndx\\
&\quad \quad \quad =\int_{\bR^{nd}} \int_{\bR^d}\cF \varphi_{{\bf t},\underline{\xi}}(\xi)\varphi(\xi)\psi(\xi_1,\ldots,\xi_n)d\xi_1 \ldots d\xi_n d\xi,
\end{align*}
using Plancherel theorem and hypothesis {\em (iv)} for the last equality. This shows that for almost all $\xi_1,\ldots,\xi_n,\xi$ in $\bR^d$,
\begin{equation}
\label{2-Fourier-equal}
\cF \varphi_{{\bf t},\underline{\xi}}(\xi)=\cF f_n(t_1,\cdot,\ldots,t_n,\cdot,t,\cdot)(\xi_1,\ldots,\xi_n).
 \end{equation}

Note that $y \mapsto \cF g(s_1,\cdot,\ldots,s_n,\cdot,s,y)(\xi_1,\ldots,\xi_n)$ is a function in $\cS(\bR^d)$ whose Fourier transform evaluated at $\xi$ is equal to $\cF g(s_1,\cdot,\ldots,s_{n-1},\cdot,s,\cdot)(\xi_1,\ldots,\xi_{n-1},\xi)$.

We use \eqref{MalliavinG-u-step1}, but we express differently $\langle f_{n-1}(\cdot,t,x),g(\cdot,s,y) \rangle_{\cH^{\otimes (n-1)}}$ using the Fourier transforms in the space variables. Using also Fubini's theorem, we obtain:
\begin{align*}
&E[\langle DG,u \rangle_{\cH}]
=n!\int_{\bR_{+}^{2n}}\gamma(t-s)\prod_{i=1}^{n-1}\gamma(t_i-s_i)\int_{\bR^{(n-1)d}}\mu(d\xi_1)\ldots \mu(d\xi_{n-1})
\Big(\int_{\bR^{2d}}dxdy
f(x-y)  \\
& \cF f_{n-1}(t_1,\cdot,\ldots,t_{n-1},\cdot,t,x)(\xi_1,\ldots,\xi_{n-1})
\overline{\cF g(s_1,\cdot,\ldots,s_{n-1},\cdot,s,y)(\xi_1,\ldots,\xi_{n-1})}\Big) d{\bf t}d{\bf s}dtds\\
&=n!\int_{\bR_{+}^{2n}}\gamma(t-s)\prod_{i=1}^{n-1}\gamma(t_i-s_i)\int_{\bR^{(n-1)d}}
\left(\int_{\bR^{d}}\cF f_{n-1}(t_1,\cdot,\ldots,t_{n-1},\cdot,t,\cdot)(\xi_1,\ldots,\xi_{n-1},\xi)\right.\\
&\quad \quad \quad \left.\overline{\cF g(s_1,\cdot,\ldots,s_{n-1},\cdot,s,\cdot)(\xi_1,\ldots,\xi_{n-1},\xi)}\mu(d\xi)\right)
\mu(d\xi_1)\ldots \mu(d\xi_{n-1})d{\bf t}d{\bf s}dtds\\
&=n!\langle g,f_{n-1}\rangle_{\cH^{\otimes n}}=n!\langle g,\widetilde{f}_{n-1}\rangle_{\cH^{\otimes n}}=E[GI_{n}(\widetilde{f}_{n-1})].
\end{align*}
Note that the second equality above
is justified by \eqref{2-Fourier-equal} and hypothesis {\em (iv)}. $\Box$

\vspace{3mm}

The next result is used to prove the existence and uniqueness of the solution to equation \eqref{wave} in dimension $d\geq 3$, being applied to the process $v^{(t,x)}$ given by Definition \ref{new-def-solution}.(c).
This result gives a correction to Proposition 2.5 of \cite{balan12}, whose proof is incorrect since the second equality on page 12, line 18 (which states that the action of the random distribution $u(\bullet)$ on $\psi_{\varepsilon}* \widetilde{\phi}$ is equal to the series $\sum_{n\geq 0}I_n((f_n(\cdot,\bullet),\psi_{\varepsilon}* \widetilde{\phi}))$) cannot be justified.

\begin{proposition}
\label{corr-prop-25-2}
Let $u=\{u(t,\cdot);t\geq 0\}$ be a process with values in $\cS'(\bR^d)$, defined on
a probability space $(\Omega,\cF,P)$ such that
for any $t \geq 0$ and $\varphi \in \cS(\bR^d)$, $E|(u(t,\cdot),\varphi)|^2<\infty$ and
\begin{equation}
\label{chaos-exp-u-e}
(u(t,\cdot),\varphi)=\sum_{n\geq 0}I_n(S_n^{t,\varphi}) \quad \mbox{in} \ L^2(\Omega),
\end{equation}
where $S_0^{t,\varphi}=E[(u(t,\cdot),\varphi)]$, and for $n \geq 1$, $S_{n}^{t,\varphi} \in \cH^{\otimes n}$ is such that $S_{n}^{t,\varphi}(t_1,\cdot,\ldots,t_n,\cdot)$ is a symmetric distribution in $\cS'(\bR^{nd})$ whose Fourier transform is a function which has a version such that $(t_1,\ldots,t_n,\xi_1,\ldots,\xi_n)\mapsto\cF S_n^{t,\varphi}(t_1,\cdot,\ldots,t_n,\cdot)(\xi_1,\ldots,\xi_n)$ is measurable.

Suppose that $u(t,\cdot)=0$ if $t>T$ and:\\
(i) for any $(\omega,t)\in \Omega \times [0,T]$, the Fourier transform of $u(\omega,t,\cdot)$ is a function which has a version such that $(\omega,t,\xi)\mapsto \cF u(\omega,t,\cdot)(\xi)$ is measurable,
the map $t \mapsto \cF u(\omega,t,\cdot)(\xi)$ is integrable on $[0,T]$ for almost all $(\omega,\xi)\in \Omega \times \bR^d$, and
\begin{equation}
\label{u-cond}
E\left[\int_{\bR^d} \int_{\bR}\left|\int_0^T e^{-i\tau t}\cF u(t,\cdot)(\xi)dt \right|^2 \nu(d\tau)\mu(d\xi)\right]<\infty;
\end{equation}

(ii) for any $n\geq 1$ and $t_1,\ldots,t_n,t \in [0,T]$, there is a distribution $f_n(t_1,\cdot,\ldots,t_n,\cdot,t,\cdot)$ in $\cS'(\bR^{(n+1)d})$ such that for any $\psi \in \cS(\bR^{nd})$ and $\varphi \in \cS(\bR^d)$,
\begin{equation}
\label{relation-fn-Sn}
\big( f_n(t_1,\cdot,\ldots,t_n,\cdot,t,\cdot), \psi \otimes \varphi \big)=\big(S_n^{t,\varphi}(t_1,\cdot,\ldots,t_n,\cdot), \psi \big);
\end{equation}

(iii) $f_n(t_1,\cdot,\ldots,t_n,\cdot,t,\cdot)$ satisfies assumption (iii) of Proposition \ref{corr-prop-25-1};

(iv) for any $t_1,\ldots,t_n,t\in [0,T]$, $\xi_1,\ldots,\xi_n\in \bR^{d}$ and $\phi \in \cS(\bR^d)$, the map $\xi \mapsto \cF f_n(t_1,\cdot,\ldots,t_n,\cdot,t,\cdot)(\xi_1,\ldots,\xi_n,\xi)\phi(\xi)$ is in $\cS(\bR^d)$.

We denote by $\widetilde{f}_n$ the symmetrization of $f_n$ in all $n+1$ variables.
Then $u \in {\rm Dom} \ \delta$ if and only if $V:=\sum_{n\geq 0}I_{n+1}(\widetilde{f}_n)$ converges in $L^2(\Omega)$. In this case, $\delta(u)=V$.
\end{proposition}

\noindent {\bf Proof:} We use the same argument as in the proof of Proposition \ref{corr-prop-25-0}. We only need to show the statement in Step 1, since Steps 2, 3 and 4 remain valid without any modification.

Hypothesis {\em (i)} guarantees that $u\in\cH$ a.s. and $E\|u\|_{\cH}^2<\infty$,
by Theorem \ref{phi-in-H}.b). Similarly, hypothesis {\em (iii)} implies that $f_n \in \cH^{\otimes n}$ by Theorem \ref{phi-in-Hn}.c).

We prove relation \eqref{MalliavinG-u} by regularizing $u$ in space. Let $\phi \in \cD(\bR^d)$ be such that $\phi \geq 0$, the support of $\phi$ is included in the unit ball in $\bR^d$, and $\int_{\bR^d}\phi(x)dx=1$. For any $\varepsilon>0$, let $\phi_{\varepsilon}(x)=\varepsilon^{-d}\phi(x/\varepsilon)$ for all $x \in \bR^d$. Then $\cF \phi_{\varepsilon}(\xi)\to 0$ as $\varepsilon \to 0$ and $|\cF \phi_{\varepsilon}(\xi)|\leq 1$ for all $\xi \in \bR^d$. For any $\omega \in \Omega$, $t \in [0,T]$ and $x \in \bR^d$, let
$u_{\varepsilon}(\omega,t,x)=\big(u(\omega,t,\cdot)*\phi_{\varepsilon}\big)(x)$.
Then $u_{\varepsilon}(\omega, t,\cdot)$ is a $C^{\infty}$-function with polynomial growth (hence a distribution in $\cS'(\bR^d)$),
whose Fourier transform is the function $\cF u(\omega,t,\cdot)(\xi)=\cF u(\omega,t,\cdot)(\xi)\cF \phi_{\varepsilon}(\xi)$ (see Theorem 7.19 of \cite{rudin91}). Hence
$\int_{0}^T e^{-i \tau t}\cF u_{\varepsilon}(t,\cdot)(\xi)dt=\cF \phi_{\varepsilon}(\xi)\int_{0}^T e^{-i \tau t}\cF u(t,\cdot)(\xi)dt$ for any $\tau \in \bR$ and $\xi \in \bR^d$. This implies that
$E\left[\int_{\bR^d}\int_{\bR} \left|\int_0^T e^{-i \tau t}\cF u_{\varepsilon}(t,\cdot)(\xi)d\tau\right|^2 \nu(d\tau)\mu(d\xi)\right]<\infty$ using  \eqref{u-cond} and the fact that $|\cF \phi_{\varepsilon}(\xi)| \leq 1$.
By Theorem \ref{phi-in-H}.b), $u_{\varepsilon} \in \cH$ a.s. and $E\|u_{\varepsilon}\|_{\cH}^2<\infty$. Moreover, by the dominated convergence theorem and \eqref{u-cond}, as $\varepsilon \to 0$, $$E\|u_{\varepsilon}-u\|_{\cH}^2=E\left[\int_{\bR^d} \int_{\bR}|\cF \phi_{\varepsilon}(\xi)-1|^2 \left|\int_0^T e^{-i\tau t}\cF u(t,\cdot)(\xi)dt \right|^2 \nu(d\tau)\mu(d\xi)\right] \to 0.$$

By \eqref{chaos-exp-u-e}, $u_{\varepsilon}(t,x)$ has the Wiener chaos expansion:
$$u_{\varepsilon}(t,x)=\big(u(t,\cdot),\phi_{\varepsilon}(x-\cdot) \big)=\sum_{n\geq 0}I_n(f_{n,\varepsilon}(\cdot,t,x)) \quad \mbox{in} \ L^2(\Omega),$$
where $f_{n,\varepsilon}(\cdot,t,x):=S_{n}^{t,\phi_{\varepsilon}(x-\cdot)}$.
The idea is to write relation \eqref{MalliavinG-u} for the process $u_{\varepsilon}$ and let $\varepsilon \to 0$.
For this, we need to check that $u_{\varepsilon}$ satisfies the hypotheses of Proposition \ref{corr-prop-25-1}. We have already proved that $u_{\varepsilon}$ satisfies hypothesis {\em (i)} of Proposition \ref{corr-prop-25-1}. It remains to check hypotheses {\em (ii)}-{\em (iv)}.

First, note that relation \eqref{relation-fn-Sn} allows us to compute $\cF f_{n,\varepsilon}(t_1,\cdot,\ldots,t_n,\cdot,t,x)$.
For any $\psi \in \cS(\bR^{nd})$ and $\varphi \in \cS(\bR^d)$,
\begin{align*}
&\mbox{LHS of \eqref{relation-fn-Sn}}=\big(\cF f_n(t_1,\cdot,\ldots,t_n,\cdot,t,\cdot), \cF^{-1}\psi \otimes \cF^{-1}\varphi \big)=\\
&\int_{\bR^{nd}} \left(\int_{\bR^d}\cF f_n(t_1,\cdot,\ldots,t_n,\cdot,t,\cdot)(\xi_1,\ldots,\xi_n,\xi)\cF^{-1}\varphi(\xi)d\xi \right)\cF^{-1}\psi(\xi_1,\ldots,\xi_n)d\xi_1 \ldots d\xi_n
\end{align*}
and
\begin{align*}
&\mbox{RHS of \eqref{relation-fn-Sn}}=\big(\cF S_{n}^{t,\varphi}(t_1,\cdot,\ldots,t_n,\cdot),\cF^{-1} \psi \big)=\\
&\int_{\bR^{nd}}\cF S_{n}^{t,\varphi}(t_1,\cdot,\ldots,t_n,\cdot)(\xi_1,\ldots,\xi_n)\cF^{-1}\psi(\xi_1, \ldots,\xi_n)d\xi_1 \ldots d\xi_n.
\end{align*}
Since this happens for any $\psi \in \cS(\bR^{nd})$, we infer that for any $\varphi \in \cS(\bR^d)$,
$$\cF S_{n}^{t,\varphi}(t_1,\cdot,\ldots,t_n,\cdot)(\xi_1,\ldots,\xi_n)=\int_{\bR^d}\cF f_n(t_1,\cdot,\ldots,t_n,\cdot,t,\cdot)(\xi_1,\ldots,\xi_n,\xi)\cF^{-1}\varphi(\xi)d\xi. $$
Recalling the definition of $f_{n,\varepsilon}$, we obtain that:
\begin{align}
\nonumber
&\cF f_{n,\varepsilon}(t_1,\cdot,\ldots,t_n,\cdot,t,x)(\xi_1,\ldots,\xi_n)=
\cF S_{n}^{t,\phi_{\varepsilon}(x-\cdot)}(t_1,\cdot,\ldots,t_n,\cdot)
(\xi_1,\ldots,\xi_n)\\
\nonumber
&\quad \quad \quad =\int_{\bR^d}\cF f_n(t_1,\cdot,\ldots,t_n,\cdot,t,\cdot)(\xi_1,\ldots,\xi_n,\xi)
\cF^{-1}\phi_{\varepsilon}(x-\cdot)(\xi)d\xi\\
\label{expression-F-fne}
& \quad \quad \quad =\frac{1}{(2\pi)^d}\int_{\bR^d}e^{i \xi \cdot x}\cF f_n(t_1,\cdot,\ldots,t_n,\cdot,t,\cdot)(\xi_1,\ldots,\xi_n,\xi)
\cF \phi_{\varepsilon}(\xi)d\xi\\
\label{expression-F-fne2}
& \quad \quad \quad =\cF^{-1}\big(\cF f_n(t_1,\cdot,\ldots,t_n,\cdot,t,\cdot)(\xi_1,\ldots,\xi_n,*)
\cF \phi_{\varepsilon} \big)(x),
\end{align}
where for the last equality we used hypothesis {\em (iv)}.

We show that $f_{n,\varepsilon}$ satisfies hypothesis {\em (ii)} of Proposition \ref{corr-prop-25-1}. For any $t_1,\ldots,t_n,t\in [0,T]$, let $f_{n,\varepsilon}(t_1,\cdot,\ldots,t_n,\cdot,t,\cdot)$ be the distribution in $\cS'(\bR^{(n+1)d})$ whose Fourier transform is the function
\begin{equation}
\label{def-Four-fn-allvar}
\cF f_{n,\varepsilon}(t_1,\cdot,\ldots,t_n,\cdot,t,\cdot)(\xi_1,\ldots,\xi_n,\xi):=\cF f_{n}(t_1,\cdot,\ldots,t_n,\cdot,t,\cdot)(\xi_1,\ldots,\xi_n,\xi)\cF \phi_{\varepsilon}(\xi).
\end{equation}
More precisely, for any $h \in \cS(\bR^{(n+1)d})$,
\begin{align*}
\big(f_{n,\varepsilon}(t_1,\cdot,\ldots,t_n,\cdot,t,\cdot),h\big)&=\int_{\bR^{(n+1)d}}
\cF f_{n}(t_1,\cdot,\ldots,t_n,\cdot,t,\cdot)(\xi_1,\ldots,\xi_n,\xi)\cF \phi_{\varepsilon}(\xi) \\
& \quad \quad \quad \cF^{-1}h(\xi_1,\ldots,\xi_n,\xi)d\xi_1 \ldots d\xi_n d\xi\\
&= \big(f_{n}(t_1,\cdot,\ldots,t_n,\cdot,t,\cdot),H_{h} \big),
\end{align*}
where $H_h \in \cS(\bR^{(n+1)d})$ is such that $\cF^{-1}H_{h}(\xi_1,\ldots,\xi_n,\xi)=\cF \phi_{\varepsilon}(\xi)\cF^{-1}h(\xi_1,\ldots,\xi_n,\xi)$.

We claim that for any $\psi \in \cS(\bR^{nd})$ and $\varphi \in \cS(\bR^d)$,
\begin{equation}
\label{relation-Fne-fne}
\big(f_{n,\varepsilon}(t_1,\cdot,\ldots,t_n,\cdot,t,\cdot),\psi \otimes \varphi\big) = \int_{\bR^d}\big(f_{n,\varepsilon}(t_1,\cdot,\ldots,t_n,\cdot,t,x),\psi\big) \varphi(x) dx.
\end{equation}
To prove this, note that the right-hand side of \eqref{relation-Fne-fne} is equal to
\begin{align*}
& \int_{\bR^d}\big(\cF f_{n,\varepsilon}(t_1,\cdot,\ldots,t_n,\cdot,t,x),\cF^{-1} \psi\big) \varphi(x) dx\\
&=\int_{\bR^d} \left(\int_{\bR^{nd}}\cF f_{n,\varepsilon}(t_1,\cdot,\ldots,t_n,\cdot,t,x)(\xi_1,\ldots,\xi_n) \cF^{-1}\psi(\xi_1,\ldots,\xi_n)d\xi_1 \ldots d\xi_n \right) \varphi(x)dx\\
&=\int_{\bR^{(n+1)d}}\cF f_{n}(t_1,\cdot,\ldots,t_n,\cdot,t,\cdot)(\xi_1,\ldots,\xi_n,\xi)\cF \phi_{\varepsilon}(\xi) \cF^{-1}\psi(\xi_1,\ldots,\xi_n)\\
& \quad \quad \quad \quad \quad
\left(\frac{1}{(2\pi)^d}\int_{\bR^d}e^{i \xi \cdot x} \varphi(x)dx\right)
d\xi_1 \ldots d\xi_n d\xi\\
&=\int_{\bR^{(n+1)d}}\cF f_{n}(t_1,\cdot,\ldots,t_n,\cdot,t,\cdot)(\xi_1,\ldots,\xi_n,\xi)\cF \phi_{\varepsilon}(\xi) \cF^{-1}\psi(\xi_1,\ldots,\xi_n) \cF^{-1}\varphi(\xi)
d\xi_1 \ldots d\xi_n d\xi\\
&=\big(f_{n,\varepsilon}(t_1,\cdot,\ldots,t_n,\cdot,t,\cdot),\psi \otimes \varphi\big),
\end{align*}
where we used \eqref{expression-F-fne} for the second equality.

The fact that $\cF f_{n,\varepsilon}(t_1,\cdot,\ldots,t_n,\cdot,t,\cdot)$ satisfies hypothesis {\em (iii)} of Proposition \ref{corr-prop-25-1} follows from hypothesis {\em (iii)}, \eqref{def-Four-fn-allvar} and the fact that $|\cF \phi_{\varepsilon}(\xi)| \leq 1$ for all $\xi \in \bR^d$.

Finally, $f_{n,\varepsilon}$ satisfies hypothesis {\em (iv)} of Proposition \ref{corr-prop-25-1} since for any ${\bf t}=(t_1,\ldots,t_n)\in [0,T]^n$ and $\underline{\xi}=(\xi_1,\ldots,\xi_n)\in \bR^{nd}$, the map $x \mapsto \cF f_{n,\varepsilon}(t_1,\cdot, \ldots,t_n,\cdot,t,x)=:\varphi_{t,\underline{\xi}}(x)$ is in $\cS(\bR^d)$ by \eqref{expression-F-fne2} and hypothesis {\em (iv)}.
This concludes the verification of hypotheses {\em (i)}-{\em (iv)} of Proposition \ref{corr-prop-25-1} for the process $u_{\varepsilon}$.

By Theorem \ref{phi-in-Hn}.c), we infer that $f_{n,\varepsilon} \in \cH^{\otimes (n+1)}$ for any $n\geq 0$. By the dominated convergence theorem and \eqref{def-Four-fn-allvar}, $\|f_{n,\varepsilon}-f_n\|_{\cH^{\otimes (n+1)}}^2 \to 0$ as $\varepsilon \to 0$.

We are now ready to conclude the proof. We write \eqref{MalliavinG-u} for the process $u_{\varepsilon}$:
$$E[\langle DG,u_{\varepsilon} \rangle_{\cH}]=E[G I_{n}(\widetilde{f}_{n-1,\varepsilon})].$$
We let $\varepsilon \to 0$. Using Cauchy-Schwarz inequality, we see that the left-hand side converges to $E[\langle DG,u\rangle_{\cH}]$, since $E\|u_{\varepsilon}-u\|_{\cH}^2 \to 0$ as $\varepsilon \to 0$. Similarly, the right-hand side converges to $E[G I_{n}(\widetilde{f}_{n-1})]$, since $\|f_{n-1,\varepsilon}-f_{n-1}\|_{\cH^{\otimes n}}^2 \to 0$ as $\varepsilon \to 0$. Relation \eqref{MalliavinG-u} follows.
$\Box$

\end{document}